\documentclass[11pt,reqno]{amsart}
\usepackage[T1]{fontenc}
\usepackage{lmodern}
\usepackage{amsmath,amssymb,mathtools,slashed,mathrsfs}
\usepackage{microtype}
\usepackage[margin=30mm]{geometry}
\usepackage{enumitem}
\usepackage{xcolor}
\usepackage[hidelinks]{hyperref}
\newcommand{\revision}[1]{{\color{black}#1}}
\numberwithin{equation}{section}
\newtheorem{theorem}{Theorem}[section]
\newtheorem{proposition}[theorem]{Proposition}
\newtheorem{lemma}[theorem]{Lemma}
\newtheorem{corollary}[theorem]{Corollary}
\theoremstyle{definition}

\theoremstyle{remark}
\newtheorem{remark}[theorem]{Remark}
\newcommand{\R}{\mathbb R}
\newcommand{\C}{\mathbb C}
\newcommand{\Z}{\mathbb Z}
\newcommand{\D}{\slashed D}
\newcommand{\cP}{\mathcal P}
\newcommand{\cH}{\mathcal H}
\newcommand{\cE}{\mathcal E}
\newcommand{\dRclass}[1]{\operatorname{cl}_{\mathrm{dR}}(#1)}
\newcommand{\cR}{\mathcal R}
\newcommand{\dd}{\,d}
\newcommand{\Id}{\operatorname{Id}}
\newcommand{\Harm}{\operatorname{Harm}}
\newcommand{\rank}{\operatorname{rank}}
\newcommand{\Vol}{\operatorname{Vol}}
\newcommand{\Area}{\operatorname{Area}}
\newcommand{\Span}{\operatorname{span}}
\newcommand{\sech}{\operatorname{sech}}

\newcommand{\Rea}{\operatorname{Re}}
\newcommand{\Ima}{\operatorname{Im}}

\newcommand{\tr}{\operatorname{tr}}
\DeclareMathOperator{\diver}{div}

\allowdisplaybreaks[2]
\setlist[enumerate]{label=\textup{(\roman*)},leftmargin=*,itemsep=3pt}

\title[An existence criterion for coupled Dirac-harmonic maps]{A complete existence criterion for coupled Dirac-harmonic maps from tori to surfaces}
\author[J\"urgen Jost]{J\"urgen Jost}
\address{Max Planck Institute for Mathematics in the Sciences,
Inselstrasse 22, 04103 Leipzig, Germany}
\email{jost@mis.mpg.de}
\author[Jingyong Zhu]{Jingyong Zhu}
\address{Department of Mathematics, Sichuan University,
Chengdu 610065, China}
\email{jzhu@scu.edu.cn}
\date{}

\begin{document}
\begin{abstract}
For a closed Riemannian spin torus $(M,g)$ and a closed oriented
Riemannian surface $(N,h)$, we prove that a homotopy class
$\mathfrak f$ contains coupled Dirac-harmonic maps if and only if
$K_h\not\equiv0$ and $\deg(\mathfrak f)\chi(N)=0$.
The metrics and source spin structure are arbitrary, and \revision{each
such class contains sequences of solutions with} unbounded \revision{map} energies.
\revision{For targets of genus at least two, both conditions are automatic,
so every homotopy class occurs.}
Our main theorem treats all \revision{non-null-homotopic} classes into nonflat
targets of positive genus\revision{, including nonzero-degree classes between
tori,} and yields nowhere-zero spinors and maps of rank two on a
nonempty open set. \revision{Together} with the known existence result for
null-homotopic classes and the genus--degree obstruction, \revision{this gives}
the criterion above. The proof relies on an explicit
description of the twisted Dirac kernel and constrained energy minimization,
identifying the curvature force with the term arising from the
Lagrange multipliers.
To exclude harmonic minimizers, we derive restrictions on their connection periods. We also prove that every harmonic map of nonzero degree between Riemannian tori is an unbranched covering, without any curvature sign assumption.
\end{abstract}
\maketitle
\enlargethispage{2pt}

\section{Introduction}\label{sec:introduction}

We determine when a homotopy class of maps from a closed spin torus
to a closed oriented surface contains a coupled Dirac-harmonic map.
The metrics and the source spin structure are prescribed throughout.
Dirac-harmonic maps were introduced in
\cite{CJLW} as critical points of a conformally invariant functional
that couples map energy to a spinor field. They consist of a map
$\phi:(M,g)\to(N,h)$ and a section $\Psi$ of
$\Sigma M\otimes_{\R}\phi^*TN$ satisfying
\begin{equation}\label{eq:DHintro}
 \tau_g(\phi)=\cR_g(\phi,\Psi),\qquad \D^\phi\Psi=0.
\end{equation}
Here $\tau_g$ is the tension field and $\cR_g$ is the curvature force,
quadratic in $\Psi$ and linear in $d\phi$. We call a solution
\emph{coupled} when $\tau_g(\phi)\not\equiv0$.
The spinor then contributes a nonzero term to the map equation, and
the map is not harmonic. A nonzero spinor alone does not ensure coupling.

Index theory \cite{AGindex} gives nonzero twisted harmonic spinors
along harmonic maps, hence uncoupled solutions. Heat-flow methods
\cite{JZflow,JSZ} also give solutions with nonzero spinors but do not
ensure coupling.
For the $\alpha $-perturbed system, Jost--Zhu \cite{JZalpha} obtain
coupled solutions under a uniform energy hypothesis; their limits
as $\alpha\downarrow1 $ need not remain coupled.

Explicit constructions give coupled solutions for special geometries:
maps from noncompact surfaces into hyperbolic three-space in
\cite{JMZ}, and examples based on immersions in \cite{AGexamples}.
More recently, B\'ek\'esi \cite{Bekesi} constructed coupled
examples from closed product manifolds.
These constructions do not give existence for arbitrary prescribed metrics and
source spin structures. For surface targets,
\cite{Yang,CJSZsurfaces} give genus--degree obstructions to coupling.

Our recent work \cite{JZcurves} gives the first general existence
theorem for coupled Dirac-harmonic maps from closed spin tori:
for every source metric and spin structure and every nonflat target
surface metric, there are null-homotopic coupled maps with unbounded
energies. For strictly negatively curved oriented targets, it also
treats every homotopy class with trivial source spin structure;
the other spin structures require a compatibility condition.
These maps factor through curves and have degree zero.

In the present paper, we \revision{treat arbitrary nonflat metrics on
targets of positive genus and every source spin structure.
This removes the} compatibility condition \revision{in \cite{JZcurves}} and
\revision{includes} nonzero-degree classes \revision{between tori}. We minimize
\[
 E_g(\phi)=\frac12\int_M|d\phi|_{g,h}^2\dd V_g
\]
under two constraints on the harmonic part of the pullback connection.
A homotopy class is called nontrivial if it is not
null-homotopic.

\begin{theorem}[Prescribed nontrivial homotopy classes]\label{thm:main}
Let $(M,g)$ be a closed Riemannian torus with any spin structure, and
let $(N,h)$ be a nonflat closed oriented Riemannian surface of genus
at least one. Every nontrivial homotopy class $\mathfrak f\in[M,N]$
contains a sequence of smooth Dirac-harmonic pairs $(\phi_n,\Psi_n)$
such that
\begin{equation}\label{eq:mainproperties}
 [\phi_n]=\mathfrak f,\qquad
 \tau_g(\phi_n)\not\equiv0,\qquad E_g(\phi_n)\longrightarrow\infty.
\end{equation}
Each $\Psi_n$ is nowhere zero. For each $n$, there is a nonempty open
set on which $\rank d\phi_n=2$ and $K_h\circ\phi_n\neq0$.
\end{theorem}

\revision{For a target of genus at least two with $K_h\leq0$, every nonconstant
harmonic map from a torus maps onto a closed geodesic.
Theorem~\ref{thm:main} gives nonharmonic maps in every nontrivial
homotopy class, with rank two on a nonempty open set, and requires
no curvature sign assumption. The rank-two conclusion also
distinguishes these coupled maps from the rank-one solutions
in \cite{JZcurves}.}
 Combining Theorem~\ref{thm:main} with
\cite[Theorem~1.1]{JZcurves} for null-homotopic classes and
\cite[Theorem~3.3]{CJSZsurfaces} for the genus--degree obstruction
gives the following criterion.

\begin{corollary}[Existence criterion for every spin structure]\label{cor:criterion}
Let $(M,g)$ be a closed Riemannian spin torus and $(N,h)$ a closed
oriented Riemannian surface. A class $\mathfrak f\in[M,N]$ contains
a coupled Dirac-harmonic map if and only if
\begin{equation}\label{eq:criterion}
 K_h\not\equiv0,\qquad \deg(\mathfrak f)\,\chi(N)=0.
\end{equation}
Whenever these conditions hold, the class contains coupled solutions
with unbounded map energies.
\end{corollary}

The criterion is independent of the source spin structure and allows
curvature to vanish or change sign. \revision{For targets of genus at least two, the criterion is automatically
satisfied in every homotopy class. For torus targets, only nonflatness
is required, while for sphere targets only the null-homotopic class
occurs.} 

For nonzero-degree classes between tori with trivial source spin
structure, we exclude harmonic minimizers using the following result.

\enlargethispage{\baselineskip}
\begin{theorem}[Harmonic coverings and their connection periods]\label{thm:harmonicintro}
Let $\phi:(M,g)\to(N,h)$ be a harmonic map of nonzero degree
between closed oriented Riemannian tori. Then $\phi$ is an unbranched
covering. Let $g_0$ be a flat metric conformal to $g$. Write
$h=e^{2w}h_b$ with $h_b$ flat, and set $\omega:=*_{h_b}dw$. Then
\begin{equation*}
 \Harm_{g_0}(\phi^*\omega)=0.
\end{equation*}
No curvature sign assumption is required.
\end{theorem}

\revision{The classical univalence theorem of Sampson
\cite[Theorem~11]{Sampson} and Schoen--Yau
\cite[Theorem~3.1]{SchoenYau} states that a degree-one harmonic
map between closed surfaces of the same genus at least two is a
diffeomorphism when the target curvature is nonpositive.
Theorem~\ref{thm:harmonicintro} concerns nonzero-degree maps
between tori with arbitrary metrics and also determines their
connection periods. Its proof uses the local Jacobian statement
in \cite[Proposition~2.2]{SchoenYau}, which requires no curvature
sign assumption.}

For the covering property, we associate a compact minimal immersion
with $\phi $ and compare it with its vertical translates. To prove
vanishing of the periods, we use
\[
 h _{t }=(1-t)g _{0 }+t\phi ^{* }h,\qquad
 \tau _{g _{0 },h _{t }}(\Id)=0,\qquad 0\leq t\leq1. 
\]
The connection periods along this path are continuous, lie in a
discrete set, and vanish at $t=0 $. Hence they vanish at $t=1 $;
see Section~\ref{sec:harmonic}.

We now outline the proof of Theorem~\ref{thm:main}. By conformal
invariance, we use a flat source metric $g_0$. For a nontrivial class
$\mathfrak f $ in the theorem and $\phi\in\mathfrak f $,
\[
 \mathfrak f\text{ is cyclic}
 \quad\Longleftrightarrow\quad
 \phi _{* }\pi _{1 }(M)\cong\Z
 \quad\Longleftrightarrow\quad \deg(\mathfrak f)=0. 
\]
This includes all nontrivial classes into targets of genus at least
two. The remaining classes have nonzero degree between tori.

To use a global oriented orthonormal frame, take $Y=N $ for a torus
target. For genus at least two, lift maps in $\mathfrak f $ to the
associated cyclic cover $Y$, which is diffeomorphic to a cylinder.
Let $\omega $ be the frame connection form. The map enters the twisted
connection through
\[
 A_\Phi:=\Phi^*\omega,\qquad
 \cP(\Phi):=\Harm_{g_0}(A_\Phi),\qquad \Phi:M\to Y .
\]
Here $\Harm_{g_0}$ is the $L^2$-orthogonal projection onto harmonic
one-forms; $A_\Phi$ need not be closed. In \cite{JZcurves}, the maps
factor through curves and a curvature-area constraint controls one
harmonic direction. Here we prescribe both components of $\cP(\Phi) $,
which allows general maps in both types of class.

For the spin structure $(\sigma,\delta)\in\{0,\tfrac12\}^2$, set
\[
 \cH_{\sigma,\delta}
   =(\Z+\sigma)dx^1+(\Z+\delta)dx^2,
\]
using the marked coordinate forms of Section~\ref{sec:conventions}.
Proposition~\ref{prop:kernel} describes the full twisted Dirac kernel:
\[
 \ker\D^\Phi_{\sigma,\delta}\neq\{0\}
 \quad\Longleftrightarrow\quad
 \cP(\Phi)\in\cH_{\sigma,\delta}.
\]
Under this condition, every $g_0$-parallel real one-form $j$ gives
an explicitly constructed nowhere-zero spinor $\Psi_j$ with
\[
 \D^\Phi\Psi_j=0,\qquad
 \cR(\Phi,\Psi_j)=K_Y(\Phi)Jd\Phi(j^\sharp).
\]
Here $J$ is positive rotation by $\pi/2$ on the target and
$g_0(j^\sharp,X)=j(X)$. The reduction applies to every smooth map
into either framed model, independently of minimization;
see Corollary~\ref{cor:reduction}.

For $\eta\in\cH_{\sigma,\delta}$ satisfying \eqref{eq:regularperiod}, consider
\[
 m _{\mathfrak f }(\eta)
 :=\inf\{E(\Phi):\Phi\in\mathcal B _{\mathfrak f },\ \cP(\Phi)=\eta\}, 
\]
where $\mathcal B _{\mathfrak f }$ is the finite-energy class defined in
Section~\ref{sec:compactness}. Changing the map on two disjoint disks
gives any prescribed $\eta$, and the direct method gives a minimizer
$\Phi _{\eta }$. Stokes' theorem and the degree formula show that the
constraint differentials are independent. Thus, for a parallel real
one-form $\lambda _{\eta }$, we have the weak equation
\[
 \tau(\Phi_{\eta })=-K_Y(\Phi_{\eta })Jd\Phi_{\eta }(\lambda_{\eta }^\sharp). 
\]
The projected equation has antisymmetric coefficients, so
\cite[Theorem~1.1]{SharpTopping} and elliptic estimates give
$\Phi _{\eta }\in C ^{\infty }$. With $j=-\lambda_{\eta }$, the spinor formula gives
\[
 \D ^{\Phi _{\eta }}\Psi _{j }=0,\qquad
 \cR(\Phi _{\eta },\Psi _{j })=\tau(\Phi _{\eta }) .
\]

It remains to exclude harmonic minimizers. By
Corollary~\ref{cor:harmonicperiodcondition},
\[
 \begin{gathered}
 \Phi\text{ nonconstant and harmonic}
 \quad\Longrightarrow\quad \cP(\Phi)\in\cH_{0,0},\\
 ( \sigma,\delta)\neq(0,0)
 \quad\Longrightarrow\quad \eta\notin \cH_{0,0}
 \quad\Longrightarrow\quad \tau(\Phi _{\eta })\not\equiv0 .
 \end{gathered}
\]
For the trivial spin structure, Theorem~\ref{thm:harmonicintro} gives
\[
 \deg (\mathfrak f)\neq0,\quad\tau(\Phi )=0
 \quad\Longrightarrow\quad \cP(\Phi)=0 ,
\]
so we choose $\eta\neq0$. Cyclic classes pass to a twofold quotient
with nontrivial spin structure (Lemma~\ref{lem:doublequotient}).

Finally, the uniform connection bound \eqref{eq:framebounds} and
orthogonal projection give
\[
 \begin{gathered}
 \|\eta\|_{L^2(g_0)}^2\leq 2C_\omega^2 E(\Phi_\eta), \\
 \|\eta_{n }\|_{L^2(g_0)}\to\infty
 \quad\Longrightarrow\quad E(\Phi _{\eta _{n }})\to\infty. 
 \end{gathered}
\]
The admissible lattice points are unbounded, and the twofold pullback
multiplies energy by two. This proves the energy conclusion of
Theorem~\ref{thm:main}.

The method suggests a strategy for other Dirac-type coupled
systems: encode a nonzero Dirac kernel by geometric constraints,
minimize energy under those constraints, and choose a kernel
spinor whose coupling force equals the multiplier force. Such
applications require suitable kernel formulas, compactness and
regularity results, and a way to exclude uncoupled solutions.

Sections~\ref{sec:conventions}--\ref{sec:geometry} fix the conventions
and target models. Section~\ref{sec:spinors} solves the twisted Dirac
equation and derives the period condition for harmonic maps.
Sections~\ref{sec:compactness}--\ref{sec:variation} establish constrained
minimization and regularity. Section~\ref{sec:harmonic} proves
Theorem~\ref{thm:harmonicintro}, and Section~\ref{sec:existence} completes
the existence proofs. Appendix~\ref{app:variations} contains the
variation estimates, weak chain rules, and approximation of
ambient test fields.

\section{Conventions and spin structures}\label{sec:conventions}

All surfaces are connected and oriented. Metrics, maps, and
sections are smooth unless stated otherwise.
In matrix formulas, vectors and spinor components are columns,
one-form coefficients are rows, and $\mathsf T$ denotes transpose.

\subsection{Bundles and curvature}\label{subsec:bundleconventions}

The spinor bundle $\Sigma M$ is complex, with Hermitian product
linear in its first variable. Clifford multiplication satisfies
\[
 X\cdot Y\cdot+Y\cdot X\cdot=-2g(X,Y)\Id,\qquad
 \langle X\cdot\xi,\zeta\rangle=-\langle\xi,X\cdot\zeta\rangle.
\]
We use the equivalent descriptions
\[
 \cE_\phi=\Sigma M\otimes_{\R}\phi^*TN
 \cong\Sigma M\otimes_{\C}(\phi^*TN\otimes_{\R}\C).
\]
The target metric is extended Hermitianly in the second description.
We explicitly use its complex-bilinear extension only when discussing
the Hopf differential. The tensor product connection and twisted
Dirac operator are
\[
 \nabla^{\cE_\phi}_X(\xi\otimes Z)
 =\nabla_X^\Sigma\xi\otimes Z+\xi\otimes\nabla_X^{\phi^*TN}Z,
 \qquad
 \D^\phi\Psi=\sum_{i=1}^2e_i\cdot
 \nabla^{\cE_\phi}_{e_i}\Psi .
\]
Here $(e_1,e_2)$ is a local $g$-orthonormal frame. The curvature
convention is
\[
 R^N(X,Y)Z=\nabla_X\nabla_YZ-\nabla_Y\nabla_XZ-\nabla_{[X,Y]}Z.
\]
If $J$ is positive rotation through $\pi/2$ on the target, then
\begin{equation}\label{eq:curvatureconvention}
 \begin{gathered}
 R^N(X,Y)Z=K\bigl(h(Y,Z)X-h(X,Z)Y\bigr),\qquad
 dA_h(X,Y)=h(JX,Y),\\
 R^N(X,Y)=-K\,dA_h(X,Y)J.
 \end{gathered}
\end{equation}
For smooth vector fields $X_1,X_2$ on $M$, the pullback and
tensor product connections satisfy
\begin{equation}\label{eq:pullbacktensorcurvature}
 \begin{aligned}
 R^{\phi^*TN}(X_1,X_2)
   &=R^N(d\phi(X_1),d\phi(X_2)),\\
 R^{\cE_\phi}(X_1,X_2)
   &=R^\Sigma(X_1,X_2)\otimes\Id
     +\Id\otimes R^{\phi^*TN}(X_1,X_2).
 \end{aligned}
\end{equation}
Here $R^\Sigma$ is the curvature of the spinor connection;
it vanishes when the source metric is flat.
Our normalization of the curvature force is
\begin{equation}\label{eq:force}
 h(\cR_g(\phi,\Psi),Y)
 =\frac12\sum_{i=1}^2
 \Rea\langle e_i\cdot R^N(Y,d\phi(e_i))\Psi,\Psi\rangle .
\end{equation}
\subsection{Source geometry and weak derivatives}\label{subsec:sourceconventions}

By uniformization \cite{JostRS}, the prescribed source metric is
conformal to a flat metric. Choosing an oriented lattice basis and
normalizing its periods to $2\pi$, we write
\begin{equation}\label{eq:sourceflat}
 g=e^{2\varpi}g_0,\qquad M=\R^2/(2\pi\Z)^2,\qquad
 g_0=\sum_{i,j=1}^2G_{ij}\,dx^i\,dx^j,
\end{equation}
where $G$ is constant, symmetric, and positive definite.

\emph{Points, coordinates, and lifts.}
Write $x\in  M$ and
$\widetilde x=(\widetilde x^1,\widetilde x^2)^{\mathsf T}\in\widetilde M $, where
\begin{equation}\label{eq:sourcecover}
 \begin{gathered}
 \widetilde M=\R^2,\qquad
 \operatorname{pr}_M:\widetilde M\longrightarrow M,\qquad
 \operatorname{pr}_M(\widetilde x)=[\widetilde x].
 \end{gathered}
\end{equation}
On a coordinate disk $U$, set
$x^i=\widetilde x^i\circ\operatorname{pr}_M^{-1}$ using a local inverse branch.
Branches differ by constant lattice translations; hence $x^i$ is local,
while $dx^i$ and $\partial_{x^i}$ are global, with
\begin{equation}\label{eq:sourcecoordinatepullback}
 \operatorname{pr}_M^*dx^i=d\widetilde x^i,\qquad
 d\operatorname{pr}_M(\partial_{\widetilde x^i})
   =\partial_{x^i},\qquad
 \widetilde g_0:=\operatorname{pr}_M^*g_0
   =\sum_{i,j=1}^2G_{ij}\,d\widetilde x^i\,d\widetilde x^j.
\end{equation}
For a function $v$ or a differential form $\nu$ on $M$, use
\begin{equation}\label{eq:sourcepullbacknotation}
 \begin{gathered}
 \widetilde v:=v\circ\operatorname{pr}_M,\qquad
 \widetilde\nu:=\operatorname{pr}_M^*\nu,\\
 d\widetilde v=\operatorname{pr}_M^*dv,\qquad
 \partial_{\widetilde x^i}\widetilde v
       =(\partial_{x^i}v)\circ\operatorname{pr}_M.
 \end{gathered}
\end{equation}
The lifts are periodic, and \eqref{eq:sourcepullbacknotation}
also holds weakly for $v\in H^1(M,\R^m)$ and
$\widetilde v\in  H^1_{\mathrm{loc}}(\widetilde M,\R^m)$.
We use the same notation for sections of a specified pullback bundle.

Fix $x_0=[0]$. The coordinate loops are
\[
 \begin{gathered}
 a_1,a_2:[0,2\pi]\longrightarrow M,\\
 a_1(\widetilde x^1)=[(\widetilde x^1,0)^{\mathsf T}],\qquad
 a_2(\widetilde x^2)=[(0,\widetilde x^2)^{\mathsf T}].
 \end{gathered}
\]
We identify $\pi_1(M,x_0)\cong\Z^2$ using the generators
$[a_1],[a_2]$. For a map $\phi$, $[\phi]$ denotes its free homotopy class.

\emph{Parallel frames and local derivatives.}
Let $\sqrt G$ be the symmetric positive-definite square root of $G$.
On each source coordinate disk, define oriented orthonormal
coordinates, their coordinate frame, and dual coframe by
\begin{equation}\label{eq:sourceframe}
 \begin{gathered}
 X^i=\sum_{j=1}^2(\sqrt G)_{ij}x^j,\qquad
 e^i=dX^i=\sum_{j=1}^2(\sqrt G)_{ij}\,dx^j,\\
 E_i=\partial_{X^i}
     =\sum_{j=1}^2((\sqrt G)^{-1})_{ji}\,\partial_{x^j},\qquad
 g_0=(dX^1)^2+(dX^2)^2.
 \end{gathered}
\end{equation}
Since the $X^i$ change by constants, $E_i,e^i$ are global and
parallel, with $e^i(E_j)=\delta^i_j$. Their lifts satisfy
\begin{equation}\label{eq:sourceframelifts}
 \begin{gathered}
 \widetilde e^i=\operatorname{pr}_M^*e^i,\qquad
 d\operatorname{pr}_M(\widetilde E_i)=E_i,\qquad
 \widetilde E_i(\widetilde v)=(E_i v)\circ\operatorname{pr}_M.
 \end{gathered}
\end{equation}
For a scalar one-form $\nu$, write
\begin{equation}\label{eq:covectorrows}
 \begin{gathered}
 dx:=\binom{dx^1}{dx^2},\qquad
 \nu=\nu_1dx^1+\nu_2dx^2=(\nu_1\ \nu_2)\,dx.
 \end{gathered}
\end{equation}
Its metric dual satisfies
\[
 g_0(\nu^\sharp,X)=\nu(X),\qquad
 \nu^\sharp=\sum_{i=1}^2\nu(E_i)E_i.
\]
For a target one-form, $\sharp$ instead uses the target metric.
Write
\[
 dV=e^1\wedge e^2,\qquad V_M=\int_MdV,\qquad
 \widetilde{dV}=\operatorname{pr}_M^*dV.
\]
The constant coefficients in \eqref{eq:sourceframe} give
\begin{equation}\label{eq:sourceintegration}
 \begin{gathered}
 [E_1,E_2]=0,\qquad
 \operatorname{div}_{g_0}E_i=0,\\
 \int_M E_i(\varphi)\,\rho\,dV
       =-\int_M \varphi E_i(\rho)\,dV,
 \qquad\varphi\in H^1(M),\quad\rho\in C^\infty(M).
 \end{gathered}
\end{equation}
The source Laplace operators are
\begin{equation}\label{eq:sourcelaplacian}
 \begin{gathered}
 \Delta_{g_0}:=\operatorname{div}_{g_0}\operatorname{grad}_{g_0}
       =E_1^2+E_2^2=\partial_{X^1}^2+\partial_{X^2}^2,\\
 \Delta_{\widetilde g_0}=\widetilde E_1^2+\widetilde E_2^2,\qquad
 \Delta_{\widetilde g_0}\widetilde\varphi
       =(\Delta_{g_0}\varphi)\circ\operatorname{pr}_M.
 \end{gathered}
\end{equation}
For $\varphi,\varphi_i\in L^2(M,\R^m)$, define $E_i\varphi$ componentwise by
\begin{equation}\label{eq:sourceweakderivative}
 E_i\varphi=\varphi_i
 \quad\Longleftrightarrow\quad
 \int_M \varphi E_i(\rho)\,dV=-\int_M\varphi_i\rho\,dV
 \quad\text{for every }\rho\in C^\infty(M,\R).
\end{equation}
Since $E_1,E_2$ is a global orthonormal frame,
\begin{equation}\label{eq:sourceh1criterion}
 \varphi\in H^1(M,\R^m)
 \quad\Longleftrightarrow\quad
 \varphi\in L^2(M,\R^m),\quad
 E_1\varphi,E_2\varphi\in L^2(M,\R^m),
\end{equation}
with weak derivatives as in \eqref{eq:sourceweakderivative}. Also,
\begin{equation}\label{eq:periodicintegration}
 \int_M \varphi\,dV
 =\int_{[0,2\pi]^2}\widetilde\varphi\,\widetilde{dV}.
\end{equation}
For a smooth complex function $F$ on $\widetilde M$ with
$F(\widetilde x+2\pi k)=\mu(k)F(\widetilde x)$ and $|\mu(k)|=1$,
the functions $|F|^2$ and $|\widetilde E_iF|^2$ are periodic.
If $\Delta_{\widetilde g_0} F=0$, integration by parts on the lattice cell gives
\begin{equation}\label{eq:unitaryharmonicconstant}
 \begin{aligned}
 0&=\int_{[0,2\pi]^2}\Delta_{\widetilde g_0}|F|^2\,\widetilde{dV}\\
  &=2\int_{[0,2\pi]^2}
       \bigl(|\widetilde E_1F|^2+|\widetilde E_2F|^2\bigr)
          \,\widetilde{dV}
 \quad\Longrightarrow\quad F\text{ is constant}.
 \end{aligned}
\end{equation}
Opposite boundary terms cancel by periodicity.

Using the local real coordinates in \eqref{eq:sourceframe}, set
\begin{equation}\label{eq:sourceisothermal}
 \begin{gathered}
 z_M=X^1+iX^2,\qquad dz_M=e^1+ie^2,\qquad
 \partial_{z_M}=\tfrac12(E_1-iE_2),\quad
 \partial_{\bar z_M}=\tfrac12(E_1+iE_2).
 \end{gathered}
\end{equation}
Thus $dz_M,\partial_{z_M},\partial_{\bar z_M}$ are global.
On $\widetilde M$ use the global linear coordinate
\[
 \widetilde z_M=\sum_{j=1}^2
       \bigl((\sqrt G)_{1j}+i(\sqrt G)_{2j}\bigr)\widetilde x^j.
\]
Set
\begin{equation}\label{eq:covercomplexoperators}
 \begin{gathered}
 \partial_+=\widetilde E_1+i\widetilde E_2
     =2\partial_{\overline{\widetilde z_M}},\qquad
 \partial_-=\widetilde E_1-i\widetilde E_2
     =2\partial_{\widetilde z_M},\\
 \partial_-\partial_+=\partial_+\partial_-=\Delta_{\widetilde g_0}.
 \end{gathered}
\end{equation}

\subsection{Target covers and pullback connections}\label{subsec:coverconventions}

\leavevmode We fix the covering notation for the models constructed in
Section~\ref{sec:geometry}. With its coordinate identifications,
$\widetilde Y=\R^2$ and
\begin{equation}\label{eq:coordinatecover}
 \begin{aligned}
 \operatorname{pr}_Y\!\left(\binom{\widetilde r}{\widetilde s}\right)
   &=(\widetilde r,\widetilde s\pmod L),
       &&Y=\R\times(\R/L\Z),\\
 \operatorname{pr}_Y\!\left(\binom{\widetilde y^1}{\widetilde y^2}\right)
   &=\left[\binom{\widetilde y^1}{\widetilde y^2}\right],
       &&Y=\R^2/(2\pi\Z)^2.
 \end{aligned}
\end{equation}
Here $L>0$ is the cylinder period. The translation lattice is
\begin{equation}\label{eq:targetlattice}
 \Lambda_Y=
 \begin{cases}
 \{(0,mL)^{\mathsf T}:m\in\Z\},&\text{cylinder},\\
 (2\pi\Z)^2,&\text{torus}.
 \end{cases}
\end{equation}
The coordinates $(r,s)$ and $(y^1,y^2)$ are local, except for
the global coordinate $r$. Circle-valued map components are written
$\mathsf{s}$. We have
\begin{equation}\label{eq:targetcoordinatepullback}
 \operatorname{pr}_Y^*dr=d\widetilde r,\qquad
 \operatorname{pr}_Y^*ds=d\widetilde s
 \quad\text{on the cylinder},\qquad
 \operatorname{pr}_Y^*dy^j=d\widetilde y^j
 \quad\text{on the torus}.
\end{equation}
A local inverse is
\begin{equation}\label{eq:covercoordinateforms}
 \operatorname{pr}_Y^{-1}(y)=
 \begin{cases}
 \binom{r}{s+mL},&\text{cylinder},\\[3pt]
 \binom{y^1+2\pi k_1}{y^2+2\pi k_2},&\text{torus},
 \end{cases}
\end{equation}
with constant integers on the chart. For
$\widetilde y= \operatorname{pr}_Y^{-1}(y)$ and $V\in T_yY$,
\begin{equation}\label{eq:coverinversedifferential}
 \begin{gathered}
 \bigl((d\operatorname{pr}_Y)_{\widetilde y}\bigr)^{-1}
 =d(\operatorname{pr}_Y^{-1})_y:
 T_yY\longrightarrow T_{\widetilde y}\widetilde Y,\\
 d(\operatorname{pr}_Y^{-1})_y[V]
 =\begin{cases}
 \binom{dr_y(V)}{ds_y(V)},&\text{cylinder},\\[3pt]
 \binom{(dy^1)_y(V)}{(dy^2)_y(V)},&\text{torus}.
 \end{cases}
 \end{gathered}
\end{equation}
Branches differ by elements of $\Lambda _{Y }$, so their differentials
agree in real coordinates. Pullback is injective on covariant tensors:
\begin{equation}\label{eq:coverpullbackinjective}
 \operatorname{pr}_Y^*\mathcal T=0\quad\Longrightarrow\quad\mathcal T=0,
\end{equation}
and the same holds for $\operatorname{pr}_M$.

For a map $\Phi:M\to Y$, distinguish the following two lifts:
\begin{equation}\label{eq:mapliftdomains}
 \begin{gathered}
 \widetilde\Phi:\widetilde M\to\widetilde Y,\qquad
 \operatorname{pr}_Y\circ\widetilde\Phi
       =\Phi\circ\operatorname{pr}_M,\\
 \widetilde\Phi_U:U\to\widetilde Y,\qquad
 \operatorname{pr}_Y\circ\widetilde\Phi_U=\Phi|_U.
 \end{gathered}
\end{equation}
For smooth maps, choose compatible lifts on $\widetilde M$ and
simply connected disks $U$. Then
\begin{equation}\label{eq:mapliftdifferentials}
 \begin{gathered}
 \widetilde\Phi_U
       =\widetilde\Phi\circ\operatorname{pr}_M^{-1},\\
 (d\operatorname{pr}_Y)_{\widetilde\Phi_U(x)}
       \circ d\widetilde\Phi_U|_x=d\Phi|_x,\\
 (d\operatorname{pr}_Y)_{\widetilde\Phi(\widetilde x)}
       \circ d\widetilde\Phi|_{\widetilde x}
       =d\Phi|_{[\widetilde x]}
            \circ(d\operatorname{pr}_M)_{\widetilde x}.
 \end{gathered}
\end{equation}
The explicit real Sobolev representatives below satisfy these
identities almost everywhere. On connected overlaps, local lifts
differ by constant elements of $\Lambda_Y$.

\emph{Pullback connections and complex frames.}
For a smooth map $\Phi:(M,g_0)\to(Y,h_Y)$, write $\nabla^\Phi$
for the pullback of the Levi--Civita connection $\nabla^Y$.
Torsion-freeness gives
\begin{equation}\label{eq:pullbacktorsion}
 \nabla^\Phi_{X_1}\{d\Phi(X_2)\}
 -\nabla^\Phi_{X_2}\{d\Phi(X_1)\}
 =d\Phi([X_1,X_2]).
\end{equation}
Here $X_1,X_2$ are smooth vector fields on $M$.
Since $E_1,E_2$ are parallel and commute, we have
\begin{equation}\label{eq:parallelframetension}
 \begin{gathered}
 \nabla^\Phi_{E_1}\{d\Phi(E_2)\}
   =\nabla^\Phi_{E_2}\{d\Phi(E_1)\},\\
 \tau(\Phi):=\operatorname{tr}_{g_0}\nabla d\Phi
   =\sum_{i=1}^2\nabla^\Phi_{E_i}\{d\Phi(E_i)\}.
 \end{gathered}
\end{equation}

For a global positively oriented $h_Y$-orthonormal frame $P_1,P_2$,
the rotation $J$ and the connection one-form $\omega$ satisfy
\begin{equation}\label{eq:targetframeconnection}
 \begin{gathered}
 JP_1=P_2,\qquad JP_2=-P_1,\\
 \nabla^Y P_1=\omega\otimes P_2,\qquad
 \nabla^Y P_2=-\omega\otimes P_1.
 \end{gathered}
\end{equation}
Along $\Phi$, set
\begin{equation}\label{eq:complexframe}
 \begin{gathered}
 A_\Phi:=\Phi^*\omega,\qquad
 Z:=(P_1+iP_2)\circ\Phi,\qquad
 \overline Z:=(P_1-iP_2)\circ\Phi,\\
 Z,\overline Z\in C^\infty(M,\Phi^*TY\otimes_\R\C).
 \end{gathered}
\end{equation}
Extend the connection and differential complex linearly.
Pulling back \eqref{eq:targetframeconnection} gives
\begin{equation}\label{eq:complexframeconnection}
 \begin{aligned}
 \nabla^\Phi_X Z&=-iA_\Phi(X)Z,\qquad 
 \nabla^\Phi_X\overline Z=iA_\Phi(X)\overline Z.
 \end{aligned}
\end{equation}

\subsection{Harmonic forms and phases}\label{subsec:hodgeconventions}

On source one-forms the Hodge star satisfies
$\alpha\wedge *\beta=\langle\alpha,\beta\rangle_{g_0}dV$; thus
$*e^1=e^2$ and $*e^2=-e^1$. A real one-form $\nu$ is harmonic
when $d\nu=d(*\nu)=0$. Since $*_g=*_{g_0}$ on one-forms,
the harmonic forms for $g$ and $g_0$ coincide and are precisely
the $g_0$-parallel forms:
\[
 \cH^1(M)=\Span_\R\{dx^1,dx^2\}.
\]
Write $\widetilde *$ for the Hodge star of $\widetilde g_0$; then
\[
 \operatorname{pr}_M^*(*\nu)=\widetilde *\,\widetilde\nu.
\]
For a smooth real one-form $A=A_1e^1+A_2e^2$, its normalized
Hodge decomposition is
\begin{equation}\label{eq:sourcehodge}
 A=\eta+df+*du,\qquad \eta=\Harm(A),\qquad
 \int_M f\,dV=\int_M u\,dV=0.
\end{equation}
The potentials are the unique mean-zero solutions of
\begin{equation}\label{eq:hodgepotentials}
 \Delta_{g_0} f=E_1(A_1)+E_2(A_2),\qquad
 \Delta_{g_0} u=E_1(A_2)-E_2(A_1).
\end{equation}
The right-hand sides have mean zero by \eqref{eq:sourceintegration}.
For $\eta=\eta_1dx^1+\eta_2dx^2$, set
\begin{equation}\label{eq:phaseprimitive}
 \begin{gathered}
 \vartheta_0(\widetilde x)=\eta_1\widetilde x^1+\eta_2\widetilde x^2,
 \qquad \vartheta=\vartheta_0+\widetilde f,\\
 d\vartheta=\widetilde\eta+d\widetilde f,\qquad
 \widetilde A-d\vartheta=\widetilde *\,d\widetilde u,\\
 \vartheta(\widetilde x+2\pi k)-\vartheta(\widetilde x)
      =2\pi(\eta_1k_1+\eta_2k_2).
 \end{gathered}
\end{equation}
The phase factors $e^{\pm i\vartheta}$ and amplitudes
$e^{\pm\widetilde u}$ enter the spinor formulas in
Section~\ref{sec:spinors}.

\subsection{Norms}\label{subsec:normconventions}

Single bars denote pointwise norms; tensor norms are Hilbert--Schmidt.
For $\Phi:(M,g_0)\to(Y,h_Y)$ and a section $W$ of $\Phi^*TY$,
\begin{equation}\label{eq:pointwisenorms}
 \begin{gathered}
 |d\Phi|_{g_0,h_Y}^2=\sum_{i=1}^2|d\Phi(E_i)|_{h_Y}^2,\qquad 
 |\nabla^\Phi W|_{g_0,h_Y}^2
     =\sum_{i=1}^2|\nabla^\Phi_{E_i}W|_{h_Y}^2.
 \end{gathered}
\end{equation}
Forms and spinors use the induced exterior-power and Hermitian
metrics. For a linear map $\mathsf A$, set
\[
 \|\mathsf A\|_{\mathrm{op}}
       :=\sup_{|V|=1}|\mathsf A V|,
\]
using the stated domain and range metrics.
For a real $m\times m$ matrix $\mathsf A=(A_{\alpha\beta})$,
the Frobenius norm is
\begin{equation}\label{eq:frobeniusnorm}
 |\mathsf A|_{\mathrm F}^{\,2}
 :=\sum_{\alpha,\beta=1}^m|A_{\alpha\beta}|^2. 
\end{equation}
For real $m\times m$ matrices $\mathsf A,\mathsf B$, we have
\begin{equation}\label{eq:matrixproductbounds}
 |\mathsf A\mathsf B|_{\mathrm F}
       \leq\|\mathsf A\|_{\mathrm{op}}|\mathsf B|_{\mathrm F},\qquad
 |\mathsf B\mathsf A|_{\mathrm F}
       \leq|\mathsf B|_{\mathrm F}\|\mathsf A\|_{\mathrm{op}}.
\end{equation}
The norms $\|\cdot\|_2$ and $\|\cdot\|_\infty$ use the domain of
the given field, with measure $dV $ on $M$. Unless specified otherwise, source
tensors use $g_0$, target tensors use $h_Y$ (or $h$ on $N$), and
real coordinate functions use the Euclidean metric. In particular,
$|d\Phi|=|d\Phi|_{g_0,h_Y}$ and $|d\phi|=|d\phi|_{g_0,h}$.

\subsection{Spin structures and conformal covariance}\label{subsec:spinconventions}

\leavevmode Since $\operatorname{pr}_M^*\Sigma_{g_0}M$ is flat and
$\widetilde M$ is simply connected, choose a parallel unitary
frame in which Clifford multiplication is
\begin{equation}\label{eq:cliffordmatrices}
 \widetilde E_1\cdot=\begin{pmatrix}0&-1\\1&0\end{pmatrix},\qquad
 \widetilde E_2\cdot=\begin{pmatrix}0&i\\i&0\end{pmatrix}.
\end{equation}
The chirality operator is
\[
 i\widetilde E_1\cdot\widetilde E_2\cdot
       =\begin{pmatrix}1&0\\0&-1\end{pmatrix}.
\]
The first and second components have positive and negative chirality,
respectively \cite{LawsonMichelsohn}.

For the constant unit spinors, \eqref{eq:cliffordmatrices} gives
\begin{equation}\label{eq:cliffordbasisaction}
 \begin{gathered}
 \begin{aligned}
 \widetilde E_2\cdot\binom10
     &=i\widetilde E_1\cdot\binom10,\\
 \widetilde E_2\cdot\binom01
     &=-i\widetilde E_1\cdot\binom01,
 \end{aligned}\\
 \left|\widetilde E_i\cdot\binom10\right|
 =\left|\widetilde E_i\cdot\binom01\right|=1,
 \qquad i=1,2.
 \end{gathered}
\end{equation}

For $\xi\in C^\infty(M,\Sigma_{g_0}M)$, denote its pulled-back
section and its component function in this frame by
$\widetilde\xi=(\widetilde\xi_+,\widetilde\xi_-)^{\mathsf T}:
\widetilde M\to\C^2$.
The spin structure $(\sigma,\delta)\in\{0,\tfrac12\}^2$ prescribes
\begin{equation}\label{eq:spintransitions}
 \widetilde\xi(\widetilde x+2\pi k)
   =\chi_{\sigma,\delta}(k)\widetilde\xi(\widetilde x),\qquad
 \chi_{\sigma,\delta}(k)=e^{2\pi i(\sigma k_1+\delta k_2)}.
\end{equation}
The trivial spin structure $(0,0)$ has periodic components and
admits nonzero parallel spinors.

In this parallel frame,
\begin{equation}\label{eq:lifteddirac}
 \begin{gathered}
 \widetilde{\nabla^\Sigma_{E_i}\xi}
      =\widetilde E_i(\widetilde\xi),\\
 \widetilde{\D}_{g_0}\widetilde\xi
 :=\widetilde{\D_{g_0}\xi}
 =\sum_{i=1}^2\widetilde E_i\cdot\widetilde E_i(\widetilde\xi)
 =\begin{pmatrix}0&-\partial_-\\ \partial_+&0\end{pmatrix}
       \widetilde\xi.
 \end{gathered}
\end{equation}
For a real form $A=A_1e^1+A_2e^2$ on $M$, put
$D_A=\D_{g_0}-iA^\sharp\cdot$ and denote its lift by
$\widetilde D_{\widetilde A}$. With $\widetilde A_i=A_i\circ\operatorname{pr}_M$,
\begin{equation}\label{eq:liftedmagneticdirac}
 \begin{aligned}
 \widetilde D_{\widetilde A}
 &=\widetilde{\D}_{g_0}-i\widetilde A^\sharp\cdot\\
 &=\begin{pmatrix}
      0&-\partial_-+i\widetilde A_1+\widetilde A_2\\
      \partial_+-i\widetilde A_1+\widetilde A_2&0
    \end{pmatrix}.
 \end{aligned}
\end{equation}
Thus $D_A\xi=0$ is equivalent to
$\widetilde D_{\widetilde A}\widetilde\xi=0$ with \eqref{eq:spintransitions}.

By \eqref{eq:phaseprimitive}, $e^{i\vartheta_0}$ satisfies
\eqref{eq:spintransitions} precisely when
$\eta_1-\sigma,\eta_2-\delta\in\Z$. Define
\begin{equation}\label{eq:spinlattice}
 \cH_{\sigma,\delta}
 =\{(m+\sigma)dx^1+(n+\delta)dx^2:m,n\in\Z\}.
\end{equation}
Thus $\cH_{0,0}$ consists of parallel forms with periods in $2\pi\Z$,
and $\cH_{\sigma,\delta}=(\sigma dx^1+\delta dx^2)+\cH_{0,0}$.
Since $e^{i\widetilde f}$ is periodic, $e^{i\vartheta}$ has the same
spin transformation law as $e^{i\vartheta_0}$.
In Section~\ref{sec:spinors}, we take $A=A_\Phi$.

Under the canonical spin-bundle isometry
$\mathcal B:\Sigma_{g_0}M\to\Sigma_gM$, conformal covariance
\cite{CJLW,LawsonMichelsohn} reads
\begin{align}
 \Psi_g&=e^{-\varpi/2}(\mathcal B\otimes\Id)\Psi_0,\notag\\
 \D_g^\phi\Psi_g
 &=e^{-3\varpi/2}(\mathcal B\otimes\Id)\D_{g_0}^\phi\Psi_0,\notag\\
 \tau_g(\phi)&=e^{-2\varpi}\tau_{g_0}(\phi),&
 \cR_g(\phi,\Psi_g)&=e^{-2\varpi}\cR_{g_0}(\phi,\Psi_0),\label{eq:conformal}\\
 E_g(\phi)&=E_{g_0}(\phi).\notag
\end{align}
We therefore prove existence for $g_0$ and restore $g$ at the end.
All assertions about parallel twisted spinors refer to $g_0$.

\section{Homotopy classes and framed target models}\label{sec:geometry}

We describe the homotopy classes and the target covers on which
we use a global frame. First fix the cohomology notation.

For closed real one-forms $\alpha,\beta $ on a smooth manifold $X$,
write $\dRclass{\alpha}\in H ^{1 }_{\mathrm{dR} }(X;\R) $ for the de Rham class. Then
\begin{equation}\label{eq:derhamclass}
 \begin{gathered}
 \dRclass{\alpha}=\dRclass{\beta}
 \quad\Longleftrightarrow\quad
 \alpha-\beta=d\varphi
 \quad\text{for some }\varphi\in C^\infty(X,\R).
 \end{gathered}
\end{equation}
For a smooth map $\phi:M\to N$ and closed forms on $N $,
\begin{equation}\label{eq:derhampullback}
 \phi^*\dRclass{\alpha}:=\dRclass{\phi^*\alpha},\qquad
 \phi^*(d\varphi)=d(\varphi\circ\phi).
\end{equation}

Fix a nontrivial homotopy class $\mathfrak f$ from $M$ to a target
$N$ of positive genus. To use a global oriented orthonormal frame,
we work on a covering $\pi:Y\to N$ with
\begin{equation}\label{eq:targetcovermetric}
 \begin{gathered}
 h_Y=\pi^*h,\qquad K_Y=K_h\circ\pi.
 \end{gathered}
\end{equation}
Write $\Phi:M\to Y$ for a lift and $\phi=\pi\circ\Phi$ for its
projection. For a torus target, $Y=N$ and $\pi=\Id$.
For genus at least two, we construct the cyclic cover associated
with $\mathfrak f$. In both cases, $\pi$ is a local isometry.

\subsection{Targets of genus at least two}

Suppose that $N$ has genus at least two. By uniformization
\cite{JostRS}, there is a metric $h_b$ of curvature $-1$ in the
conformal class of $h$, so $h=e^{2w}h_b$ for a smooth function $w$.
On the upper half-plane $\mathbb H^2=\{\zeta\in\C:\Ima  \zeta>0\}$,
choose an oriented universal covering
\[
 p_N:\mathbb H^2\longrightarrow N,\qquad
 p_N^*h_b=g_{\mathrm{hyp}}
       =\frac{|d\zeta|^2}{(\Ima\zeta)^2},\quad \Ima\zeta>0.
\]
Let $\Gamma=\operatorname{Deck}(p_N)$, so $N=\mathbb H^2/\Gamma$.
Each deck transformation preserves $p_N^*h_b$ and $p_N^*h$ because
$p_N\circ T=p_N$ for $T\in  \Gamma$.

Fix a smooth reference map in the prescribed free homotopy class:
\[
 \phi_{\mathrm{ref}}:M\longrightarrow N,\qquad
 [\phi_{\mathrm{ref}}]=\mathfrak f.
\]
Set $y_0=\phi_{\mathrm{ref}}(x_0)$ and choose $z_0\in p_N^{-1}(y_0)$.
For a loop representing $\gamma$, let $T_\gamma(z_0)$ be the endpoint
of its lift from $z_0$. Then \cite[Section 1.3]{Hatcher} gives
\begin{equation}\label{eq:deckidentification}
 \pi_1(N,y_0)\xrightarrow{\;\cong\;}\Gamma,
 \qquad\gamma\longmapsto T_\gamma. 
\end{equation}

A closed surface $N$ of positive genus has contractible universal
cover, so it is a $K(\pi_1(N),1)$. Thus, for any connected CW complex $M$,
\begin{equation}\label{eq:homotopyclassification}
 [M,N]\cong
 \operatorname{Hom}(\pi_1(M),\pi_1(N))/\text{conjugacy}.
\end{equation}
See \cite[Propositions 1B.9 and 4A.2]{Hatcher}. In particular,
$[\phi]$ is null-homotopic exactly when $\phi_*\pi_1(M)=\{1\}$.
Every nonidentity $T  _{\gamma  }\in  \Gamma$ is hyperbolic
\cite[Section 1.1.3]{FarbMargalit}. Its centralizer
$C_\Gamma(T_\gamma)$ preserves its oriented axis.
Signed translation gives an injective homomorphism
$\ell_\gamma:C_\Gamma(T_\gamma)\to\R$ whose image is discrete by
proper discontinuity. Hence the centralizer, and every nontrivial
abelian subgroup of $\Gamma$, is infinite cyclic.

Since $\pi_1(M,x_0)=\Z^2$ is abelian and $\mathfrak f$ is nontrivial,
the induced homomorphism $\mathfrak f_*:=(\phi_{\mathrm{ref}})_*$
has image $\langle c\rangle\subset\pi_1(N,y_0)$. For the coordinate
loops $a_1,a_2$ fixed in Section~\ref{sec:conventions}, write
\begin{equation}\label{eq:cyclicclass}
 \begin{gathered}
 p=(p_1,p_2)\in\Z^2,\qquad
 \mathfrak f_*([a_i])=c^{p_i}\quad(i=1,2),\\
 \gcd(|p_1|,|p_2|)=1.
 \end{gathered}
\end{equation}
Here $c$ generates the image subgroup and may be a proper power
in $\pi_1(N,y_0)$. We fix $c  $ and $p$ from the reference map for
the rest of the argument. No curvature sign assumption on $h$ is used.

Under \eqref{eq:deckidentification}, $\langle c\rangle  $ corresponds
to $H_c:=\langle T_c\rangle\subset\Gamma$. Set
\begin{equation}\label{eq:cycliccoverprojection}
 Y=\mathbb H^2/H_c,\qquad
 q(\zeta)=[\zeta]_{H_c},\qquad
 \pi([\zeta]_{H_c})=[\zeta]_\Gamma,
 \qquad p_N=\pi\circ q.
\end{equation}
Since $H_c\subset\Gamma$, $\pi$ is well defined. The universal
cover $q $ identifies $\pi _{1 }(Y,q(z _{0 })) $ with $H _{c }$, giving
\begin{equation}\label{eq:cycliccovergroup}
 \pi_1(Y,q(z_0))\cong\Z,\qquad
 \pi_*(m)=c^m,\qquad \pi_*\pi_1(Y)=\langle c\rangle.
\end{equation}
For a based map $\phi:(M,x_0)\to(N,y_0)$, the lifting criterion
\cite[Proposition 1.33]{Hatcher} states that
\begin{equation}\label{eq:coverliftingcriterion}
 \begin{gathered}
 \exists\,\Phi:(M,x_0)\to(Y,q(z_0))
 \quad\text{such that }\pi\circ\Phi=\phi\\
 \Longleftrightarrow\quad
 \phi_*\pi_1(M,x_0)\subseteq\pi_*\pi_1(Y,q(z_0)).
 \end{gathered}
\end{equation}
For $\phi_{\mathrm{ref}}$, both subgroups in
\eqref{eq:coverliftingcriterion} equal $\langle c\rangle$.
Thus there is a lift $\Phi_{\mathrm{ref}}:M\to Y$ with
$\Phi  _{\mathrm{ref}  }(x  _{0  })=q(z  _{0  })  $. Homotopy lifting then gives a
lift to the same $Y$ for every $\phi\in\mathfrak f$
\cite[Proposition 1.30]{Hatcher}.

We now construct coordinates on $Y$. Orient the axis of $T_c$
in its translation direction. Choose an orientation-preserving
isometry $\mathscr{A} $ of $\mathbb H^2$ taking this axis to $i\R_{>0}$, directed toward
$\infty$. Then
\begin{equation}\label{eq:normalizeddeck}
 \mathscr{A} (T_c\zeta)=e^L\mathscr{A} (\zeta),\qquad \zeta\in\mathbb H^2,\qquad
 L=\ell_c(T_c)=\inf_{\zeta\in\mathbb H^2}
                  d_{\mathrm{hyp}}(\zeta,T_c\zeta)>0.
\end{equation}
The number $L$ is determined by $c$ and $h_b$.

Write
\[
 \mathscr{A} (\zeta)=X_{\mathbb H}+iY_{\mathbb H},\qquad Y_{\mathbb H}>0.
\]
The global Fermi parametrization $\mathcal F:\R^2\to\mathbb H^2$
and its inverse are
\begin{equation}\label{eq:fermiexplicit}
 \begin{gathered}
 \mathcal F(\widetilde r,\widetilde s)=\mathscr{A} ^{-1}\!\left(e^{\widetilde s}(\tanh \widetilde r+i\sech \widetilde r)\right),\\
 \widetilde r=\operatorname{arsinh}(X_{\mathbb H}/Y_{\mathbb H}),\qquad
 \widetilde s=\log\sqrt{X_{\mathbb H}^2+Y_{\mathbb H}^2}.
 \end{gathered}
\end{equation}
The inverse formulas show that $\mathcal F$ is a diffeomorphism.
Direct differentiation gives
\[
 dX_{\mathbb H}^2+dY_{\mathbb H}^2
 =e^{2\widetilde s}
       (\sech^2\widetilde r\,d\widetilde r^2+d\widetilde s^2).
\]
Since $\mathscr{A} $ is an isometry and
$Y_{\mathbb H}=e^{\widetilde s}\sech\widetilde r$,
\begin{equation}\label{eq:fermiequivariance}
 \mathcal F^*g_{\mathrm{hyp}}
 =(\mathscr{A} \circ\mathcal F)^*g_{\mathrm{hyp}}
 =d\widetilde r^2+\cosh^2\widetilde r\,d\widetilde s^2,
 \qquad
 T_c\circ\mathcal F(\widetilde r,\widetilde s)
 =\mathcal F(\widetilde r,\widetilde s+L).
\end{equation}
Let
\[
 q_L:\R^2\longrightarrow\R\times(\R/L\Z),\qquad
 q_L(\widetilde r,\widetilde s)=(\widetilde r,\widetilde s\pmod L).
\]
The equivariance in \eqref{eq:fermiequivariance} gives a diffeomorphism
\begin{equation}\label{eq:cylinderidentification}
 \overline{\mathcal F}:\R\times(\R/L\Z)\longrightarrow Y,
 \qquad q\circ\mathcal F=\overline{\mathcal F}\circ q_L.
\end{equation}
Put $h_{\mathrm{cyl}}=dr^2+\cosh^2r\,ds^2$. Since
$q^*\pi^*h_b=g_{\mathrm{hyp}}$,
\[
 q_L^*\overline{\mathcal F}^{\,*}\pi^*h_b
 =\mathcal F^*g_{\mathrm{hyp}}=q_L^*h_{\mathrm{cyl}}.
\]
Injectivity of pullback \eqref{eq:coverpullbackinjective} yields
\begin{equation}\label{eq:cylmetric}
 \begin{gathered}
 \overline{\mathcal F}^{\,*}(\pi^*h_b)
       =h_{\mathrm{cyl}}=dr^2+\cosh^2r\,ds^2,\\
 h_Y=\pi^*h=e^{2w_Y}\pi^*h_b,\qquad w_Y=w\circ\pi.
 \end{gathered}
\end{equation}
We henceforth use \eqref{eq:cylinderidentification} to identify $Y$
with the cylinder and write $h_b$ on $Y$ for $\pi^*h_b$.
Before this identification, $\operatorname{pr}  _{Y  }=\overline{\mathcal F}\circ q  _{L  }$;
in cylinder coordinates, $\operatorname{pr}  _{Y  }=q  _{L  }$.
All following coordinate formulas use this identification and the
pulled-back metrics in \eqref{eq:cylmetric}.
Thus $Y$ is a noncompact cylinder, with the global closed angular
one-form $ds$ of period $L$.

In these coordinates, the background and target orthonormal frames are
\begin{equation}\label{eq:cylinderframe}
 \begin{aligned}
 Q_1&=\partial_r,& Q_2&=\sech r\,\partial_s,\\
 P_1&=e^{-w_Y}Q_1,& P_2&=e^{-w_Y}Q_2.
 \end{aligned}
\end{equation}
They are invariant under $s\mapsto s+L$ and hence global.
For the background connection $\nabla^b$, we have
$\omega_b=\sinh r\,ds$. The conformal change formula
\[
 \nabla^Y_XV=\nabla^b_XV+X(w_Y)V+V(w_Y)X
                    -h_b(X,V)\nabla^bw_Y
\]
and \eqref{eq:targetframeconnection}, with
$*_{h_b}dr=\cosh r\,ds$ and $*_{h_b}ds=-\sech r\,dr$, give
\begin{equation}\label{eq:cylconnection}
 \omega=\sinh r\,ds+*_{h_b}dw_Y
 =-\frac{\partial_sw_Y}{\cosh r}\,dr+
   \bigl(\sinh r+\partial_rw_Y\cosh r\bigr)\,ds .
\end{equation}
Using $d(*_{h_b}dw_Y)=(\Delta_{h_b}w_Y)dA_{h_b}$ and
$K_Y=-e^{-2w_Y}(1+\Delta_{h_b}w_Y)$, we obtain
\begin{equation}\label{eq:connectioncurvature}
 d\omega=-K_YdA_{h_Y}.
\end{equation}

Compactness of $N$ and local isometry of $\pi$ give
\[
 \sup_Y\bigl(|w_Y|+|dw_Y|_{h_b}
                  +|(\nabla^b)^2w_Y|_{h_b}\bigr)<\infty.
\]
Moreover,
\[
 \begin{gathered}
 \omega_b=\tanh r\,(\cosh r\,ds),\qquad
 \nabla^b\omega_b
 =\sech^2r\,dr\otimes(\cosh r\,ds)
   -\tanh^2r\,(\cosh r\,ds)\otimes dr,\\
 |\omega_b|_{h_b}\leq1,\qquad |\nabla^b\omega_b|_{h_b}\leq1.
 \end{gathered}
\]
The conformal connection formula therefore bounds $\omega$ and
$\nabla^Y\omega$ uniformly. Differentiating
\eqref{eq:targetframeconnection} gives
\[
 (\nabla^Y)^2P_{i  }
 = (\nabla^Y\omega)\otimes J  P_{i  }-  \omega\otimes\omega\otimes P_{i  },
 \qquad i=  1,2.
\]
Thus there are constants $C_\omega\geq1$ and $C<\infty$,
independent of $y\in Y$, such that
\begin{equation}\label{eq:framebounds}
 \begin{gathered}
 c_h:=\min_N e^{2w}>0,\qquad C_h:=\max_N e^{2w}<\infty,\\
 c_hh_b\leq h_Y\leq C_hh_b,\qquad
 |\omega|_{h_Y}\leq C_\omega,\\
 \begin{aligned}
 &|K_Y|+|\nabla^Y\omega|_{h_Y}
       +|\nabla^YP_1|_{h_Y}+|\nabla^YP_2|_{h_Y}\\
 &\qquad+|(\nabla^Y)^2P_1|_{h_Y}+|(\nabla^Y)^2P_2|_{h_Y}\leq C.
 \end{aligned}
 \end{gathered}
\end{equation}

Write
\[
 \operatorname{pr}_s:Y\longrightarrow\R/L\Z,\qquad
 \operatorname{pr}_s(r,s\pmod L)=s\pmod L.
\]
Write $ds$ also for the pullback of the standard angular form to $Y$.
By \eqref{eq:cycliccovergroup}, $\pi_*(m)=c^m$, so
\[
 \frac1L\int_{a_i}\Phi^*ds=p_i\quad(i=1,2)
 \qquad\text{for }\Phi\simeq\Phi_{\mathrm{ref}}.
\]

Set
\begin{equation}\label{eq:betap}
 \begin{gathered}
 \widetilde s_p:\widetilde M\longrightarrow\R,\qquad
 \widetilde s_p(\widetilde x)=\frac{L}{2\pi}(p_1\widetilde x^1+p_2\widetilde x^2),\\
 \mathsf{s}_p:M\longrightarrow\R/L\Z,\qquad
 \mathsf{s}_p([\widetilde x])=\widetilde s_p(\widetilde x)\pmod L,\\
 \beta_p:=\mathsf{s}_p^*ds
       =\frac{L}{2\pi}(p_1dx^1+p_2dx^2)\in\Omega^1(M),
 \qquad \operatorname{pr}_M^*\beta_p=d\widetilde s_p.
 \end{gathered}
\end{equation}
The increments $L(p_1k_1+p_2k_2)$ make $\mathsf{s}_p$ well defined, and
\[
 \nabla^{g_0}\beta_p=0,\qquad \int_{a_i}\beta_p=Lp_i.
\]

For $\mathsf{s}:=\operatorname{pr}_s\circ\Phi$, a real lift
$\widetilde s:\widetilde M\to\R$ satisfies
\[
 \widetilde s(\widetilde x+2\pi k)-\widetilde s(\widetilde x)
 =L(p_1k_1+p_2k_2),\qquad k\in\Z^2.
\]
Thus $\widetilde s-\widetilde s_p=\widetilde v $ for
$v\in C^\infty(M,\R)$, and the smooth maps in the lifted class are
\[
 \Phi=\bigl(r,\mathsf{s}_p+v\pmod L\bigr),\qquad
 r,v\in C^\infty(M,\R),\qquad \Phi^*ds=\beta_p+dv,
\]
where $v$ is determined modulo $L\Z$.
Conversely, the homotopy to $\Phi_p:=(0,\mathsf{s}_p)$ is
\[
 \Phi_t=\bigl((1-t)r,\mathsf{s}_p+(1-t)v\pmod L\bigr),\qquad
 0\leq t\leq1,\qquad \Phi_0=\Phi,\quad\Phi_1=\Phi_p.
\]
By \eqref{eq:cyclicclass}, \eqref{eq:homotopyclassification}, and this homotopy,
\[
 [\pi\circ\Phi]=[\pi\circ\Phi_p]=\mathfrak f.
\]
The map $\pi\circ\Phi_p$ factors through a circle, and
$H_2(S^1;\Z)=0$, so $\deg(\mathfrak f)=0$.
Null-homotopic maps also have degree zero. Thus every map from a
torus to a surface of genus at least two has degree zero.

\subsection{Torus targets}\label{subsec:torustargets}

For a torus target, take $Y=N$ and $\pi=\Id$.
As in \eqref{eq:sourceflat}, uniformization gives marked coordinates
$y\in\R^2/(2\pi\Z)^2$ with
\[
 h=e^{2w}h_b,\qquad
 h_b=\sum_{i,j}H_{ij}\,dy^i\,dy^j,
\]
where $H=(H_{ij})$ is a constant symmetric positive-definite matrix.
Fix the oriented parallel $h_b$-orthonormal frame $Q_1,Q_2$
and its dual coframe $\alpha^1,\alpha^2$ by
\begin{equation}\label{eq:torusparallelframe}
 \alpha^{i  }=\sum_{j=1}^2(\sqrt H)_{ij  }\,dy^j,\qquad  
 Q_{i  }=\sum_{j=1}^2((\sqrt H)^{-1})_{ji  }\,\partial_{y^j},
 \qquad i  =1,  2.
\end{equation}
Here $\sqrt H$ is the symmetric positive-definite square root,
as in \eqref{eq:sourceframe}. We have
\begin{equation}\label{eq:torusparallelidentities}
 \begin{gathered}
 \alpha^i(Q_j)=\delta^i_j,\qquad
 h_b=\alpha^1\otimes\alpha^1+\alpha^2\otimes\alpha^2,\\
 \nabla^{h_b}\alpha^1=\nabla^{h_b}\alpha^2=0,\qquad
 d\alpha^1=d\alpha^2=0,\\
 *_h\alpha^1=*_{h_b}\alpha^1=\alpha^2,\qquad
 *_h\alpha^2=*_{h_b}\alpha^2=-\alpha^1,\\
 |\alpha^i|_h=e^{-w},\qquad
 (\alpha^i)^{\sharp_h}=e^{-2w}Q_i\quad(i=1,2).
 \end{gathered}
\end{equation}
Set
\begin{equation}\label{eq:torusframe}
 P_1=e^{-w}Q_1,\qquad P_2=e^{-w}Q_2.
\end{equation}
The background frame connection vanishes, so the same conformal
connection calculation gives
\begin{equation}\label{eq:torusconnection}
 \omega=*_{h_b}dw=-Q_2(w)\alpha^1+Q_1(w)\alpha^2,\qquad
 d\omega=-K_hdA_h .
\end{equation}
Here $d\omega=(\Delta_{h_b}w)dA_{h_b}$ and
$K_h=-e^{-2w}\Delta_{h_b}w$. All the bounds in
\eqref{eq:framebounds} hold by compactness. Since
$d\alpha^1=d\alpha^2=0$
and $dw=Q_1(w)\alpha^1+Q_2(w)\alpha^2$, wedging with
$\alpha^1,\alpha^2$ also gives
\begin{equation}\label{eq:torusexact}
 \alpha^1\wedge\omega=d(w\alpha^2),\qquad
 \alpha^2\wedge\omega=-d(w\alpha^1).
\end{equation}

The marked coordinates identify both fundamental groups with
$\Z^2$. Write
\[
 \phi_*\binom{k_1}{k_2}
 =B\binom{k_1}{k_2},\qquad
 B_{ji}=\frac1{2\pi}\int_{a_i}\phi^*dy^j\in\Z
 \quad(i,j=1,2).
\]
The columns of $B$ give the images of $[a_i]$; its rows give the
periods of $\phi^*dy^j$ divided by $2\pi$.
A lift $\widetilde\phi:\widetilde M\to\widetilde Y$ satisfies
\[
 \widetilde\phi(\widetilde x+2\pi k)
 =\widetilde\phi(\widetilde x)+2\pi Bk.
\]
Hence $\widetilde\phi-B\widetilde x=\widetilde v $ for
$v\in C ^{\infty }(M,\R ^{2 }) $, and the smooth maps inducing $B$ are
\begin{equation}\label{eq:torusclass}
 \phi_v([\widetilde x])
 =B\widetilde x+v([\widetilde x])\pmod{(2\pi\Z)^2},
 \qquad v\in C^\infty(M,\R^2).
\end{equation}
The homotopy $\phi_{(1-t)v}$, $0\leq t\leq1$, joins $\phi_v$
to $\phi_0$; thus $B$ determines the homotopy class.

Differentiating the local lifts in \eqref{eq:torusclass} gives
\begin{equation}\label{eq:toruscoordinatepullback}
 \phi_v^*dy^k=\sum_{i=1}^2B_{ki}\,dx^i+dv^k,
 \qquad k=1,2.
\end{equation}
Taking de Rham classes gives
\begin{equation}\label{eq:toruscohomologypullback}
 \phi_v^*\dRclass{dy^k}
 =\sum_{i=1}^2B_{ki}\,\dRclass{dx^i},\qquad k=1,2.
\end{equation}
In these cohomology bases, $\phi_v^*$ acts on coefficient rows
by right multiplication by $B$.
In particular,
\[
 \begin{aligned}
 \phi_v^*(dy^1\wedge dy^2)
 &=(\det B)dx^1\wedge dx^2\\
 &\quad+d\left(v^1\sum_iB_{2i}dx^i
               -v^2\sum_iB_{1i}dx^i+v^1dv^2\right).
 \end{aligned}
\]
Integration over $M$ therefore gives
\begin{equation}\label{eq:torusdegree}
 \deg\phi_v=\frac1{(2\pi)^2}
 \int_M\phi_v^*(dy^1\wedge dy^2)=\det B.
\end{equation}

If $B=0$, the class is null-homotopic. If $\rank B=1$, choose
an integer column vector $\mathbf c\neq0$ generating
$B\Z^2=\Z\mathbf c$. There is a unique integer row $(p_1\ p_2)$
such that
\begin{equation}\label{eq:rankonefactorization}
 B=\mathbf c\,(p_1\ p_2),\qquad \gcd(|p_1|,|p_2|)=1.
\end{equation}
For torus targets, $p$ denotes this pair, fixed by $B$ and
$\mathbf c$. Replacing $\mathbf c$ by $-\mathbf c$ replaces
$p$ by $-p$, and
\[
 \ker(B:\Z^2\to\Z^2)=\Z(-p_2,p_1)^{\mathsf T}.
\]
These are the nontrivial cyclic classes. Finally, $\rank B=2$
is equivalent to nonzero degree.
The fixed parallel parts of the pulled-back target coframe are
\begin{equation}\label{eq:torusparallelpullback}
 \begin{gathered}
 \beta_i:=\phi_0^*\alpha^i
      =\sum_{j=1}^2(\sqrt HB)_{ij}\,dx^j,\qquad
 \phi_v^*\alpha^i
      =\beta_i+d\left(\sum_{j=1}^2(\sqrt H)_{ij}v^j\right),\\
 \nabla^{g_0}\beta_i=0\quad(i=1,2),\qquad
 \dim\Span\{\beta_1,\beta_2\}=\rank B.
 \end{gathered}
\end{equation}
In rank one,
\[
 \beta_i=(\sqrt H\mathbf c)_i(p_1dx^1+p_2dx^2),\qquad
 \Span\{\beta_1,\beta_2\}
       =\R(p_1dx^1+p_2dx^2).
\]
This is the same line as that generated by $\beta_p$ in the
cylindrical model.

\subsection{Connection periods}\label{subsec:connectionperiods}

\emph{Curvature area and parallel transport.}
Let $\gamma:[0,2\pi]\to Y$ be a smooth closed curve.
For a unit vector field $V  $ along $\gamma$, let
$\theta_{\mathrm{pt}}:[0,2\pi]\to\R$ be its angle in $(P  _{1  },P  _{2  })  $.
By \eqref{eq:targetframeconnection},
\[
 \nabla^\gamma V=(d\theta_{\mathrm{pt}}+\gamma^*\omega)\otimes JV.
\]
For parallel $V$, if $\gamma$ positively parametrizes
$\partial\mathcal D $ for a disk $\mathcal D\subset Y$, Stokes' theorem gives
\begin{equation}\label{eq:curvatureareaangle}
 \theta_{\mathrm{pt}}(2\pi)-\theta_{\mathrm{pt}}(0)
 =-\int_\gamma\omega
 =\int_{\mathcal D}K_YdA_{h_Y}.
\end{equation}
The rotation angle is understood modulo $2\pi$.

\emph{The harmonic part of the connection.}
For a smooth map $\Phi:M\to Y$, the pullback connection form
$A_\Phi=\Phi^*\omega$ in \eqref{eq:complexframe} satisfies
\begin{equation}\label{eq:connectionformmeaning}
 (A_\Phi)_x(X)
 =\omega_{\Phi(x)}(d\Phi_xX)
 =h_Y\bigl(\nabla^\Phi_X(P_1\circ\Phi),P_2\circ\Phi\bigr)_x.
\end{equation}
Since $(e^i,e^j)_{L^2}=V_M\delta_{ij}$, the harmonic projection is
\begin{equation}\label{eq:period}
 \begin{aligned}
 \cP(\Phi):=\Harm(A_\Phi)
 &=\frac1{V_M}\sum_{i=1}^2
       \left(\int_M A_\Phi(E_i)\dd V\right)e^i.
 \end{aligned}
\end{equation}
Here $\Harm=\Harm_{g_0}$. Set
$I_i(\Phi):=\int_M A_\Phi(E_i)\dd V$, $i=1,2$, and $\eta=\cP(\Phi)$.
We call $\int_{a_{i  }}\eta$ the connection periods in the chosen target
frame. They determine $\eta$ by
\[
 \eta=\frac1{2\pi}\left(\int_{a_1}\eta\right)dx^1
      +\frac1{2\pi}\left(\int_{a_2}\eta\right)dx^2.
\]
These periods need not equal $\int  _{a  _{i  }}A_\Phi$:
by \eqref{eq:sourcehodge},
\[
 \int_{a_i}A_\Phi=\int_{a_i}\eta+\int_{a_i}*du,
 \qquad i=1,2.
\]
The pullback curvature is
\begin{equation}\label{eq:pullbackconnectioncurvature}
 dA_\Phi=-\Phi^*(K_YdA_{h_Y}).
\end{equation}
For a source disk $U\subset M$, the rotation angle for parallel
transport in $\Phi^*TY$ around the positively oriented boundary
$\partial U$ is, modulo $2\pi$,
\[
 -\int_{\partial U}A_\Phi
 =-\int_{\partial U}*du
 =\int_U\Phi^*(K_YdA_{h_Y}). 
\]
The spin condition $\eta\in\cH_{\sigma,\delta}$ constrains only the
harmonic part and does not force $dA_\Phi=0$.
If $dA_\Phi\not\equiv0$, a sufficiently small disk $U$ on which
$*dA_\Phi$ has a fixed nonzero sign satisfies
\[
 0<\left|\int_UdA_\Phi\right|
 \leq\|*dA_\Phi\|_{L^\infty(U)}\Vol_{g_0}(U)<2\pi.
\]
Thus there are contractible loops around which parallel transport
has nonzero rotation modulo $2\pi$.

\emph{Maps through curves.}
If $dA_\Phi=0$, then $A_\Phi=\eta+df$.
For $\eta\in\cH_{\sigma,\delta}$, this gives parallel twisted spinors
as in Remark~\ref{rem:coexact}.

To compare with \cite[Section 5]{JZcurves}, write
$\eta_{\mathrm{curv}}$ for its curvature-area primitive on a target
coordinate disk $D$. We have
\[
 d\eta_{\mathrm{curv}}=K_YdA_{h_Y},\qquad
 \omega+\eta_{\mathrm{curv}}=d\chi\quad\text{on }D
\]
for a smooth real $\chi$, since $D$ is contractible and
$d\omega=-K _{Y }dA _{h _{Y }}$. Thus, for $\Phi(M)\subset D$,
\[
 A_\Phi=-\Phi^*\eta_{\mathrm{curv}}+d(\chi\circ\Phi),\qquad
 \cP(\Phi)=-\Harm(\Phi^*\eta_{\mathrm{curv}}).
\]

For a curve $\gamma:\R/(2\pi\Z)\to D$, use the angular coordinate
$\theta _{\gamma }$ of period $2\pi $. Choose $(n_1,n_2)\in\Z^2$ and put
\[
 t([\widetilde x])=n_1\widetilde x^1+n_2\widetilde x^2\pmod{2\pi},
 \qquad dt:=t^*d\theta _{\gamma }=n_1dx^1+n_2dx^2.
\]
Writing $\gamma^*\omega=\bar a\,d\theta _{\gamma }+dF$ with $F$ periodic gives
\[
 \bar a=\frac1{2\pi}\int_\gamma\omega
       =-\frac{\mathcal A_K(\gamma)}{2\pi},\qquad
 \mathcal A_K(\gamma):=\int_\gamma\eta_{\mathrm{curv}},
\]
and hence
\[
 A_{\gamma\circ t}=\bar a\,dt+d(F\circ t),\qquad
 \cP(\gamma\circ t)
 =-\frac{\mathcal A_K(\gamma)}{2\pi}\,dt.
\]
The earlier curvature-area constraint fixes the harmonic part in
the direction $dt  $. Here we prescribe both harmonic components
and allow $dA_\Phi\neq0$.

\emph{Change of target frame.}
Let $\theta:Y\to\R/(2\pi\Z)$ be smooth and rotate the frame by
\[
 P'_1=\cos\theta\,P_1+\sin\theta\,P_2,\qquad
 P'_2=-\sin\theta\,P_1+\cos\theta\,P_2.
\]
Local real lifts define the closed form $d\theta$, with periods in
$2\pi\Z$. Differentiating gives
\[
 \omega'=\omega+d\theta,\qquad
 A'_\Phi=A_\Phi+\Phi^*d\theta.
\]
The pullback $\Phi^*d\theta$ is closed with periods in $2\pi\Z$,
and harmonic projection preserves these periods. Hence
\[
 \cP'(\Phi)-\cP(\Phi)=\Harm(\Phi^*d\theta)\in\cH_{0,0},\qquad
 \cP'(\Phi)\in\cH_{\sigma,\delta}
 \ \Longleftrightarrow\ \cP(\Phi)\in\cH_{\sigma,\delta}.
\]
The value zero in Theorem~\ref{thm:harmonicintro} refers to the
fixed conformal frame \eqref{eq:torusframe}.

\section{The complete twisted Dirac kernel}\label{sec:spinors}

Let $\Phi $ be a smooth map from the flat spin torus to either
framed model $Y$. Use $Z,\overline Z$ from \eqref{eq:complexframe}.
Their pullbacks $\widetilde Z=\operatorname{pr}_M^*Z$ and
$\overline{\widetilde Z}=\operatorname{pr}_M^*\overline Z$ satisfy
\[
 \widetilde Z(\widetilde x)
   =(P_1+iP_2)_{\Phi([\widetilde x])},\qquad
 \widetilde Z(\widetilde x+2\pi k)=\widetilde Z(\widetilde x).
\]
The complexified target metric satisfies
\begin{equation}\label{eq:complexframepairings}
 h_{\C}(Z,Z)=h_{\C}(\overline Z,\overline Z)=2,
 \qquad h_{\C}(Z,\overline Z)=0,
\end{equation}
where $h_{\C}$ denotes the Hermitian extension in this section.
Equation~\eqref{eq:complexframeconnection} shows that the connection
preserves the two complex line subbundles spanned by $Z$ and $\overline Z$.
Every twisted spinor on $M$ has a unique expression
\[
 \Psi=\xi\otimes Z+\zeta\otimes\overline Z,\qquad
 \xi,\zeta\in C^\infty(M,\Sigma_{g_0}M).
\]
Their component functions in the parallel spinor trivialization
of Section~\ref{sec:conventions} satisfy \eqref{eq:spintransitions}.
For the operators defined in Section~\ref{sec:conventions},
\begin{equation}\label{eq:splitdirac}
 \D^\Phi\Psi=0
 \quad\Longleftrightarrow\quad
 D_{A_\Phi}\xi=0,\qquad D_{-A_\Phi}\zeta=0.
\end{equation}

Apply \eqref{eq:sourcehodge} to $A _{\Phi }$:
\begin{equation}\label{eq:hodge}
 A_\Phi=\eta+df+*du,\qquad \eta=\cP(\Phi),\qquad
 \int_Mf\dd V=\int_Mu\dd V=0.
\end{equation}
Use the potentials and phase from
\eqref{eq:hodgepotentials}--\eqref{eq:phaseprimitive}.

\begin{proposition}[Kernel and curvature current]\label{prop:kernel}
The twisted Dirac kernel is nonzero if and only if
$\eta\in\cH_{\sigma,\delta}$. In that case it has complex dimension
four, and all its elements have the following expressions in the
parallel trivialization fixed in Section~\ref{sec:conventions}:
\begin{equation}\label{eq:completekernel}
 \Psi=\xi\otimes Z+\zeta\otimes\overline Z,\qquad
 \widetilde\xi=e^{i\vartheta}
       \begin{pmatrix}c_+e^{-\widetilde u}\\c_-e^{\widetilde u}\end{pmatrix},
 \qquad
 \widetilde\zeta=e^{-i\vartheta}
       \begin{pmatrix}d_+e^{\widetilde u}\\d_-e^{-\widetilde u}\end{pmatrix},
\end{equation}
where $c_\pm,d_\pm\in\C$ are arbitrary constants.
For $\eta\in\cH_{\sigma,\delta}$ and a spinor $\Psi$ given by
\eqref{eq:completekernel}, define its associated parallel real
one-form $j$ by
\begin{equation}\label{eq:current}
 j(E_1)=2\Ima(c_+\overline{c_-}-d_+\overline{d_-}),\qquad
 j(E_2)=2\Rea(c_+\overline{c_-}-d_+\overline{d_-}).
\end{equation}
Then
\begin{equation}\label{eq:currentforce}
 \cR(\Phi,\Psi)=K_Y(\Phi)Jd\Phi(j^\sharp).
\end{equation}
Conversely, assume $\eta\in\cH_{\sigma,\delta}$. For any
$g_0$-parallel real one-form $j$, including $j=0$, choose
\begin{equation}\label{eq:chooseconstants}
 d_+=d_-=0,\qquad c_+=1,\qquad
 c_-=\frac{j(E_2)-ij(E_1)}2.
\end{equation}
Substituting these constants into \eqref{eq:completekernel} gives
a nowhere-zero spinor $\Psi_j\in\ker\D^\Phi$ whose associated
one-form \eqref{eq:current} is precisely $j$.
\end{proposition}

\begin{proof}
\emph{Necessity.}
Suppose $\D^\Phi\Psi=0$ and $\Psi\not\equiv0$.
Write $A=A_\Phi$ and
$\Psi=\xi\otimes Z+\zeta\otimes\overline Z$.
On $\widetilde M$, put $\widetilde A=\operatorname{pr}_M^*A$ and define
\[
 \widehat\xi:=e^{-i\vartheta}\widetilde\xi
   =\begin{pmatrix}\widehat\xi_+\\\widehat\xi_-\end{pmatrix},
 \qquad
 \widehat\zeta:=e^{i\vartheta}\widetilde\zeta
   =\begin{pmatrix}\widehat\zeta_+\\\widehat\zeta_-\end{pmatrix}.
\]
By \eqref{eq:phaseprimitive},
$\widetilde A-d\vartheta=\widetilde *\,d\widetilde u$.
The product rule therefore gives
\begin{equation}\label{eq:spinorgaugereduction}
 \begin{aligned}
 \widetilde D_{\widetilde A}\widetilde\xi
 &=e^{i\vartheta}
   \bigl(\widetilde{\D}_{g_0}\widehat\xi
       -i(\widetilde A-d\vartheta)^\sharp\cdot\widehat\xi\bigr)
   =e^{i\vartheta}
       \widetilde D_{\widetilde *\,d\widetilde u}\widehat\xi,\\
 \widetilde D_{-\widetilde A}\widetilde\zeta
 &=e^{-i\vartheta}
   \bigl(\widetilde{\D}_{g_0}\widehat\zeta
       +i(\widetilde A-d\vartheta)^\sharp\cdot\widehat\zeta\bigr)
   =e^{-i\vartheta}
       \widetilde D_{-\widetilde *\,d\widetilde u}\widehat\zeta.
 \end{aligned}
\end{equation}
Put $\widetilde u_i=\widetilde E_i(\widetilde u)$. Substituting
$\widetilde *\,d\widetilde u
 =-\widetilde u_2\widetilde e^1+\widetilde u_1\widetilde e^2$
into \eqref{eq:liftedmagneticdirac}, and using \eqref{eq:spinorgaugereduction}, gives
\[
 \begin{aligned}
 \widetilde D_{\widetilde *\,d\widetilde u}\widehat\xi
 &=\begin{pmatrix}
 -\partial_-\widehat\xi_-+(\partial_-\widetilde u)\widehat\xi_-\\
 \partial_+\widehat\xi_++(\partial_+\widetilde u)\widehat\xi_+
 \end{pmatrix}
 =\begin{pmatrix}
 -e^{\widetilde u}\partial_-(e^{-\widetilde u}\widehat\xi_-)\\
 e^{-\widetilde u}\partial_+(e^{\widetilde u}\widehat\xi_+)
 \end{pmatrix}=0,\\
 \widetilde D_{-\widetilde *\,d\widetilde u}\widehat\zeta
 &=\begin{pmatrix}
 -\partial_-\widehat\zeta_--(\partial_-\widetilde u)\widehat\zeta_-\\
 \partial_+\widehat\zeta_+-(\partial_+\widetilde u)\widehat\zeta_+
 \end{pmatrix}
 =\begin{pmatrix}
 -e^{-\widetilde u}\partial_-(e^{\widetilde u}\widehat\zeta_-)\\
 e^{\widetilde u}\partial_+(e^{-\widetilde u}\widehat\zeta_+)
 \end{pmatrix}=0.
 \end{aligned}
\]
Define functions on $\widetilde M$ by
\begin{equation}\label{eq:weightedspinorcomponents}
 F^\xi_+=e^{\widetilde u}\widehat\xi_+,\qquad
 F^\xi_-=e^{-\widetilde u}\widehat\xi_-,\qquad
 F^\zeta_+=e^{-\widetilde u}\widehat\zeta_+,\qquad
 F^\zeta_-=e^{\widetilde u}\widehat\zeta_-.
\end{equation}
Then
\begin{equation}\label{eq:holospin}
 \begin{aligned}
 \partial_+F^\xi_+&=0,& \partial_-F^\xi_-&=0,\\
 \partial_+F^\zeta_+&=0,& \partial_-F^\zeta_-&=0.
 \end{aligned}
\end{equation}
By \eqref{eq:covercomplexoperators}, all four functions are
$\widetilde g_0$-harmonic. By \eqref{eq:phaseprimitive} and
\eqref{eq:spintransitions},
\[
 \begin{aligned}
 F^\xi_\pm(\widetilde x+2\pi k)
 &=\chi_{\sigma,\delta}(k)e^{-2\pi i(\eta_1k_1+\eta_2k_2)}
                         F^\xi_\pm(\widetilde x),\\
 F^\zeta_\pm(\widetilde x+2\pi k)
 &=\chi_{\sigma,\delta}(k)e^{2\pi i(\eta_1k_1+\eta_2k_2)}
                         F^\zeta_\pm(\widetilde x).
 \end{aligned}
\]
These multipliers have absolute value one. Hence
\eqref{eq:unitaryharmonicconstant} gives
\[
 F^\xi_+=c_+,\qquad F^\xi_-=c_-,\qquad
 F^\zeta_+=d_+,\qquad F^\zeta_-=d_-.
\]
This proves \eqref{eq:completekernel}. Since $\Psi\not\equiv0$,
at least one constant is nonzero, and its transformation law gives
\[
 \begin{cases}
 \eta_1-\sigma,\ \eta_2-\delta\in\Z,
       &\text{if }(c_+,c_-)\neq(0,0),\\
 \eta_1+\sigma,\ \eta_2+\delta\in\Z,
       &\text{if }(d_+,d_-)\neq(0,0).
 \end{cases}
\]
Because $2\sigma,2\delta\in\Z$, either condition is equivalent to
$\eta\in\cH_{\sigma,\delta}$.

\emph{Sufficiency and dimension.}
Assume $\eta\in\cH_{\sigma,\delta}$. The increment formula
\eqref{eq:phaseprimitive} gives
\[
 e^{\pm i\vartheta(\widetilde x+2\pi k)}
   =\chi_{\sigma,\delta}(k)e^{\pm i\vartheta(\widetilde x)},
\]
since $\chi_{\sigma,\delta}(k)^{-1}=\chi_{\sigma,\delta}(k)$.
For any constants in \eqref{eq:completekernel}, the periodicity of
$\widetilde u$ shows that $\widetilde\xi,\widetilde\zeta$ satisfy
\eqref{eq:spintransitions}. They therefore define spinors on $M$. Reversing the calculation
leading to \eqref{eq:holospin} gives
$D_A\xi=D_{-A}\zeta=0$, hence $\D^\Phi\Psi=0$.
The four constants in \eqref{eq:completekernel} are independent,
so $\dim_\C\ker\D^\Phi=4$.

\emph{Curvature and the current.}
Since $JZ=-iZ$ and $J\overline Z=i\overline Z$,
\eqref{eq:curvatureconvention} gives, for $V,W\in T_{\Phi(x)}Y$,
\[
 R^Y(V,W)Z=iK_Y(\Phi)dA_{h_Y}(V,W)Z,\qquad
 R^Y(V,W)\overline Z=-iK_Y(\Phi)dA_{h_Y}(V,W)\overline Z.
\]

Compute the spinor products on $\widetilde M$.
The Hermitian product is linear in its first entry. The phase
and amplitude factors cancel in the mixed products:
\begin{equation}\label{eq:spinormixedproducts}
 \widetilde\xi_+\overline{\widetilde\xi_-}
     =c_+\overline{c_-},\qquad
 \widetilde\zeta_+\overline{\widetilde\zeta_-}
     =d_+\overline{d_-}.
\end{equation}
Equation~\eqref{eq:cliffordmatrices} gives
\[
 \begin{aligned}
 \langle\widetilde E_1\cdot\widetilde\xi,\widetilde\xi\rangle
 &=-\widetilde\xi_-\overline{\widetilde\xi_+}
    +\widetilde\xi_+\overline{\widetilde\xi_-}
   =2i\Ima(c_+\overline{c_-}),\\
 \langle\widetilde E_2\cdot\widetilde\xi,\widetilde\xi\rangle
 &=i\widetilde\xi_-\overline{\widetilde\xi_+}
    +i\widetilde\xi_+\overline{\widetilde\xi_-}
   =2i\Rea(c_+\overline{c_-}),\\
 \langle\widetilde E_1\cdot\widetilde\zeta,\widetilde\zeta\rangle
 &=2i\Ima(d_+\overline{d_-}),\\
 \langle\widetilde E_2\cdot\widetilde\zeta,\widetilde\zeta\rangle
 &=2i\Rea(d_+\overline{d_-}).
 \end{aligned}
\]
On $M$, these identities and \eqref{eq:current} give
\[
 \langle E_i\cdot\xi,\xi\rangle
 -\langle E_i\cdot\zeta,\zeta\rangle=ij(E_i),\qquad
 \nabla^{g_0}j=0.
\]
For a real test vector $V$ along $\Phi$, \eqref{eq:complexframepairings}
shows that the mixed target terms vanish because $h_\C(Z,\overline Z)=0$, while
$h_\C(Z,Z)=h_\C(\overline Z,\overline Z)=2$.
Thus \eqref{eq:force} becomes
\[
 \begin{aligned}
 h_Y(\cR(\Phi,\Psi),V)
 &=\frac12\sum_{i=1}^2\Rea\Bigl[
       2iK_Y(\Phi)dA_{h_Y}(V,d\Phi E_i)\\
 &\hspace{35mm}\cdot
       \bigl(\langle E_i\cdot\xi,\xi\rangle
             -\langle E_i\cdot\zeta,\zeta\rangle\bigr)\Bigr]\\
 &=-K_Y(\Phi)\sum_{i=1}^2j(E_i)dA_{h_Y}(V,d\Phi E_i)\\
 &=K_Y(\Phi)\sum_{i=1}^2j(E_i)h_Y(V,Jd\Phi E_i)\\
 &=K_Y(\Phi)h_Y\bigl(V,Jd\Phi(j^\sharp)\bigr).
 \end{aligned}
\]
Here $dA_{h_Y}(V,W)=-h_Y(V,JW)$ and
$j^\sharp=j(E_1)E_1+j(E_2)E_2$. This proves
\eqref{eq:currentforce}.

Finally, for a prescribed parallel real form $j$, the constants
in \eqref{eq:chooseconstants} satisfy
\[
 c_+\overline{c_-}-d_+\overline{d_-}
   =\frac{j(E_2)+ij(E_1)}2.
\]
Thus
\[
 |\Psi_j|^2=2\bigl(e^{-2u}+|c_-|^2e^{2u}\bigr)
   =2e^{-2u}+\frac{|j|_{g_0}^2}{2}e^{2u}>0
 \quad\text{on }M.
\]

\end{proof}

\begin{corollary}[Equivalent formulation of the coupled system]\label{cor:reduction}
Let $Y$ be either framed target model of Section~\ref{sec:geometry},
and let $\Phi:(M,g_0)\to Y$ be smooth. It is the map component of
a coupled Dirac-harmonic pair for the spin structure
$(\sigma,\delta)$ if and only if there is a $g_0$-parallel real
one-form $j$ such that
\begin{equation}\label{eq:completereduction}
 \cP(\Phi)\in\cH_{\sigma,\delta},\qquad
 \tau(\Phi)=K_Y(\Phi)Jd\Phi(j^\sharp)\not\equiv0.
\end{equation}
For each such $j$, the constants \eqref{eq:chooseconstants}
in \eqref{eq:completekernel} define a nowhere-zero spinor
completing $\Phi$ to a coupled Dirac-harmonic pair.
\end{corollary}

\begin{proof}
For a coupled pair, Proposition~\ref{prop:kernel} gives the period
condition and \eqref{eq:currentforce}, hence
\eqref{eq:completereduction}. Conversely,
\eqref{eq:chooseconstants} gives a nowhere-zero spinor with
$\D^\Phi\Psi=0$ and $\cR(\Phi,\Psi)=\tau(\Phi)\not\equiv0$.
\end{proof}

\begin{corollary}[A period condition for harmonic maps]\label{cor:harmonicperiodcondition}
If $\Phi:(M,g_0)\to Y$ is smooth, nonconstant, and harmonic, then
\begin{equation}\label{eq:harmonicperiodcondition}
 \cP(\Phi)\in\cH_{0,0}.
\end{equation}
Consequently, for a nontrivial source spin structure, no nonconstant
harmonic map satisfies $\cP(\Phi)\in\cH_{\sigma,\delta}$.
\end{corollary}

\begin{proof}
For this proof, equip $(M,g_0)$ with the trivial spin structure
and write $\D^\Phi_{0,0}$ for the corresponding twisted Dirac
operator. By \eqref{eq:spintransitions}, the constant spinor
$\widetilde\xi_0=\binom01$ descends to a parallel unit spinor
$\xi_0$ on $M$. Set
\begin{equation}\label{eq:harmonicspinorconstruction}
 \Psi_\Phi:=\sum_{i=1}^2 E_i\cdot\xi_0\otimes d\Phi(E_i).
\end{equation}
Equation~\eqref{eq:cliffordbasisaction} gives
\[
 \Psi_\Phi=(E_1\cdot\xi_0)\otimes
       \bigl(d\Phi(E_1)-i\,d\Phi(E_2)\bigr),
 \qquad |E_1\cdot\xi_0|=1.
\]
Since $d\Phi(E_1),d\Phi(E_2)$ are real, the Hermitian extension
$h_\C$ of the target metric gives
\[
 \begin{aligned}
 |\Psi_\Phi|^2
 &=\bigl|d\Phi(E_1)-i\,d\Phi(E_2)\bigr|_{h_\C}^2\\
 &=|d\Phi(E_1)|_{h_Y}^2+|d\Phi(E_2)|_{h_Y}^2
  =|d\Phi|_{g_0,h_Y}^2.
 \end{aligned}
\]
Since $E_1,E_2,\xi_0$ are parallel, \eqref{eq:parallelframetension}
gives
\begin{equation}\label{eq:harmonicspinordirac}
 \begin{aligned}
 \D^\Phi_{0,0}\Psi_\Phi
 &=\sum_{i,j=1}^2 E_j\cdot E_i\cdot\xi_0
       \otimes\nabla^\Phi_{E_j}\{d\Phi(E_i)\}\\
 &=-\xi_0\otimes\tau(\Phi)
   +(E_1\cdot E_2\cdot\xi_0)\otimes
       \bigl(\nabla^\Phi_{E_1}\{d\Phi(E_2)\}
              -\nabla^\Phi_{E_2}\{d\Phi(E_1)\}\bigr)\\
 &=-\xi_0\otimes\tau(\Phi)=0.
 \end{aligned}
\end{equation}
Since $\Phi$ is nonconstant, $\Psi_\Phi\not\equiv0$.
Proposition~\ref{prop:kernel}, applied with $(\sigma,\delta)=(0,0)$,
therefore gives
\[
 \ker\D^\Phi_{0,0}\neq\{0\}
 \quad\Longrightarrow\quad
 \cP(\Phi)\in\cH_{0,0}.
\]
The last assertion follows from
\[
 (\sigma,\delta)\neq(0,0)
 \quad\Longrightarrow\quad
 \cH_{0,0}\cap\cH_{\sigma,\delta}=\varnothing.
\]
\end{proof}

\begin{remark}[Dirac solutions and parallel spinors]\label{rem:coexact}
The potential $u$ of the coexact part satisfies
\begin{equation}\label{eq:poissonu}
 (\Delta_{g_0} u)dV=dA_\Phi=-\Phi^*(K_YdA_{h_Y}).
\end{equation}
The mean-zero normalization gives
\begin{equation}\label{eq:flatpullbackpotential}
 u=0\quad\Longleftrightarrow\quad dA_\Phi=0
 \quad\Longleftrightarrow\quad\Phi^*(K_YdA_{h_Y})=0.
\end{equation}
Thus $\rank d\Phi\leq1$ implies $u=0$; nonconstant amplitudes
occur precisely when $dA_\Phi\not\equiv0$.

For $\eta\in\cH_{\sigma,\delta}$ and $c_\pm\in\C$ not both zero,
the spinor $\widetilde\xi=e^{i\vartheta}\binom{c_+}{c_-}$ descends to $M $.
By \eqref{eq:phaseprimitive} and \eqref{eq:complexframeconnection},
\begin{equation}\label{eq:phaseconnectionremainder}
 \begin{aligned}
 \nabla^\Sigma\xi&=i(\eta+df)\otimes\xi,\\
 \nabla^{\cE_\Phi}(\xi\otimes Z)
 &=i(\eta+df-A_\Phi)\otimes(\xi\otimes Z)\\
 &=-i*du\otimes(\xi\otimes Z).
 \end{aligned}
\end{equation}
Thus $\xi\otimes Z$ is parallel precisely when $u=0$.

For every smooth $\Phi:(M,g_0)\to Y$, with $g_0$ flat,
\begin{equation}\label{eq:diracparallelcriteria}
 \begin{aligned}
 \ker\D^\Phi\neq\{0\}
 &\quad\Longleftrightarrow\quad\eta\in\cH_{\sigma,\delta},\\
 \ker\nabla^{\cE_\Phi}\neq\{0\}
 &\quad\Longleftrightarrow\quad
       dA_\Phi=0\ \text{and}\ \eta\in\cH_{\sigma,\delta}.
 \end{aligned}
\end{equation}
Here $\ker\nabla^{\cE_\Phi}$ is the space of global parallel
twisted spinors. The first equivalence is
Proposition~\ref{prop:kernel}; sufficiency in the second follows
from \eqref{eq:phaseconnectionremainder}.
Conversely, a nonzero parallel spinor has constant positive norm.
By flatness of $g_0$, \eqref{eq:pullbacktensorcurvature}, and
\eqref{eq:curvatureconvention},
\begin{equation}\label{eq:twistedconnectioncurvature}
 \begin{gathered}
 R^{\cE_\Phi}(E_1,E_2)
 =dA_\Phi(E_1,E_2)(\Id\otimes J),\qquad
 (\Id\otimes J)^2=-\Id,\\
 \nabla^{\cE_\Phi}\Psi=0
 \quad\Longrightarrow\quad
 dA_\Phi(E_1,E_2)(\Id\otimes J)\Psi=0
 \quad\Longrightarrow\quad dA_\Phi=0.
 \end{gathered}
\end{equation}
Also $\D^\Phi\Psi=0$, so Proposition~\ref{prop:kernel} gives
$\eta\in\cH_{\sigma,\delta}$. For
$\eta\in\cH_{\sigma,\delta}$ and $dA_\Phi\not\equiv0$,
every nonzero spinor in \eqref{eq:completekernel} is a global
Dirac solution and is not parallel.
\end{remark}

\section{Finite-energy classes, compactness, and period identities}
\label{sec:compactness}

In this section and Section~\ref{sec:variation}, assume
$K_h\not\equiv0$. This provides the curved target disk used to
prescribe arbitrary values of $I_1,I_2$.
Corollary~\ref{cor:harmonicperiodcondition} and the harmonic-map
results of Section~\ref{sec:harmonic} do not require this assumption.

\subsection{Finite-energy classes}

We use $H^1=W^{1,2}$ with the weak derivatives of
\eqref{eq:sourceweakderivative}--\eqref{eq:sourceh1criterion}.
We use Poincare's inequality and Rellich compactness on the closed
torus:
\begin{equation}\label{eq:sobolevcompactfacts}
 \|v-\overline v\|_2\leq C_P\|dv\|_2,
 \qquad \overline v=V_M^{-1}\int_Mv\dd V,
 \qquad H^1(M)\Subset L^2(M).
\end{equation}
Also,
\begin{equation}\label{eq:aeextraction}
 v_n\to v\text{ in }L^2
 \quad\Longrightarrow\quad
 v_{n_k}\to v\text{ a.e. for some }n_k\uparrow\infty.
\end{equation}
Together with weak compactness, this gives
\begin{equation}\label{eq:sobolevsubsequence}
 \sup_n\|v_n\|_{H^1}<\infty
 \quad\Longrightarrow\quad
 \exists\,v\in H^1,\ n_k\uparrow\infty:
 \begin{cases}
 v_{n_k}\rightharpoonup v&\text{in }H^1,\\
 v_{n_k}\to v&\text{in }L^2\text{ and a.e.}
 \end{cases}
\end{equation}

Choose $\rho\in C_c^\infty(\R^2)$ with $\rho\geq0$ and $\int\rho=1$,
and put $\rho_\epsilon=\epsilon^{-2}\rho(\,\cdot\,/\epsilon)$.
For $\varphi\in H^1(M,\R^m)$, convolution of its periodic lift
defines smooth functions on $M$ satisfying, as $\epsilon\downarrow0$,
\begin{equation}\label{eq:periodicmollification}
 \widetilde{\varphi_\epsilon}:=\rho_\epsilon*\widetilde\varphi,
 \qquad \varphi_\epsilon\to\varphi\text{ in }H^1,
 \qquad \|\varphi_\epsilon\|_\infty\leq\|\varphi\|_\infty
       \quad\text{if }\varphi\in L^\infty.
\end{equation}

We also need product estimates on fixed metric vector bundles over
$M$. Suppose $b_n\to b$ a.e. and $|b_n|\leq C$, with $b _{n }$ real.
For every fixed $T\in L^2$, dominated convergence gives
$\|(b_n-b)T\|_2\to0$. Thus, if $D_n\rightharpoonup D$ in $L^2$,
\begin{equation}\label{eq:weakproductrule}
 \begin{gathered}
 \langle b_nD_n-bD,T\rangle_{L^2}
 =\langle D_n,(b_n-b)T\rangle_{L^2}
       +\langle D_n-D,bT\rangle_{L^2}\longrightarrow0,\\
 b_nD_n\rightharpoonup bD\quad\text{in }L^2.
 \end{gathered}
\end{equation}
For $\mathcal U_n\to\mathcal U$ in $L^2$,
\begin{equation}\label{eq:strongproductestimate}
 \|b_n\mathcal U_n-b\mathcal U\|_2
 \leq C\|\mathcal U_n-\mathcal U\|_2
       +\|(b_n-b)\mathcal U\|_2\longrightarrow0.
\end{equation}
Both rules hold for uniformly bounded matrix coefficients, using
their adjoints in \eqref{eq:weakproductrule}. For one-forms
$\mu_n\to\mu$ and $\nu_n\to\nu$ in $L^2$,
\begin{equation}\label{eq:strongwedgeestimate}
 \|\mu_n\wedge\nu_n-\mu\wedge\nu\|_1
 \leq\|\mu_n-\mu\|_2\|\nu_n\|_2
       +\|\mu\|_2\|\nu_n-\nu\|_2\longrightarrow0.
\end{equation}

For the fixed cylindrical data $c,p,\mathsf{s}_p,\beta_p$ of
Section~\ref{sec:geometry}, write
\begin{equation}\label{eq:weakcylinder}
 \Phi_{r,v}=(r,\mathsf{s}_p+v\pmod L),\qquad
 \Phi_{r,v}^*ds=\beta_p+dv.
\end{equation}
Define
\begin{equation}\label{eq:weakcylinderclass}
 \mathcal B_{\mathfrak f}
 :=\left\{\Phi_{r,v}:\ r,v\in H^1(M,\R),\quad
       \cosh r\,(\beta_p+dv)\in L^2(T^*M)\right\}.
\end{equation}
For $\Phi=\Phi_{r,v}$, the energy is
\begin{equation}\label{eq:cylinderenergy}
 E(\Phi)=\frac12\int_M e^{2(w_Y\circ\Phi)}
       \bigl(|dr|_{g_0}^2+\cosh^2r\,|\Phi^*ds|_{g_0}^2\bigr)\dd V.
\end{equation}
In the torus model, keep the matrix $B$ in
\eqref{eq:torusclass} fixed and define
\begin{equation}\label{eq:weaktorusclass}
 \mathcal B_{\mathfrak f}
 :=\left\{\phi_v([\widetilde x])=B\widetilde x+v([\widetilde x])\pmod{(2\pi\Z)^2}:\quad
                         v\in H^1(M,\R^2)\right\}.
\end{equation}
The energy is
\begin{equation}\label{eq:torusenergy}
 E(\phi_v)=\frac12\int_Me^{2(w\circ\phi_v)}
                 |B\,dx+dv|_{g_0,h_b}^2\dd V .
\end{equation}
In both models, maps are defined almost everywhere, derivatives
are weak, and the class data $p$ or $B$ are fixed.
For smooth elements, \eqref{eq:toruscoordinatepullback} or
\eqref{eq:weakcylinder} shows that the weak differentials of the
real correction terms are smooth; hence these terms are smooth.
The coordinate homotopies in Section~\ref{sec:geometry} therefore give
\begin{equation}\label{eq:smoothclasshomotopy}
 \Phi\in\mathcal B_{\mathfrak f}\cap C^\infty(M,Y)
       \quad\Longrightarrow\quad [\pi\circ\Phi]=\mathfrak f.
\end{equation}
In the torus model $Y=N$, $\pi=\Id$, and $\Phi=\phi_v$.

For $\Phi\in\mathcal B_{\mathfrak f}$ in either model,
the weak pullback $A_\Phi:=\Phi^*\omega$ satisfies, by
\eqref{eq:framebounds},
\[
 |A_\Phi|_{g_0}
 \leq |\omega|_{h_Y}\circ\Phi\,|d\Phi|_{g_0,h_Y}
 \leq C_\omega|d\Phi|_{g_0,h_Y}\quad\text{a.e. on }M.
\]
Thus $A_\Phi\in L^2(T^*M)$, and \eqref{eq:period} defines
$\cP(\Phi)$ on $\mathcal B_{\mathfrak f}$. Since $\Harm$ is an
orthogonal projection and $\|d\Phi\|_2^2=2E(\Phi)$,
\begin{equation}\label{eq:energyperiod}
 \|\cP(\Phi)\|_2^2\leq\|A_\Phi\|_2^2
 \leq 2C_\omega^2 E(\Phi).
\end{equation}
For $\eta\in\cH^1(M)$ set
\begin{equation}\label{eq:constraintset}
 \mathcal A_\eta=\{\Phi\in\mathcal B_{\mathfrak f}:
                 \cP(\Phi)=\eta\},\qquad
 m_{\mathfrak f}(\eta)=\inf_{\mathcal A_\eta}E.
\end{equation}

\subsection{Prescribing the periods}

\begin{lemma}[Independent local adjustment of the periods]\label{lem:feasible}
For every $\eta\in\cH^1(M)$, there is a smooth
$\Phi\in\mathcal B_{\mathfrak f}$ with $\cP(\Phi)=\eta$.
In particular,
$m_{\mathfrak f}(\eta)<\infty$.
\end{lemma}

\begin{proof}
By \eqref{eq:period}, it suffices
to construct a smooth $\Phi\in\mathcal B_{\mathfrak f}$ with
\[
 I_i(\Phi)=V_M\eta(E_i),\qquad i=1,2.
\]

\emph{Local target coordinates.}
Choose $y_*\in Y$ with $K_Y(y_*)\neq0$ and any positively oriented
local coordinates $z=(z^1,z^2)$ with $z(y_*)=0$.
Restrict to a coordinate rectangle centered at $0$ on which $K_Y\neq0$.
Writing $h_Y=\sum_{i,j=1}^2h_{ij}(z)\,dz^i\otimes dz^j$,
we obtain from \eqref{eq:connectioncurvature} and \eqref{eq:torusconnection}
\[
 d\omega=\varrho(z)\,dz^1\wedge dz^2,\qquad
 \varrho(z):=-K_Y(z)\sqrt{\det(h_{ij}(z))}\neq0.
\]
Set
\[
 a=z^1,\qquad
 b=\int_0^{z^2}\varrho(z^1,u)\dd u.
\]
Then
\[
 \begin{aligned}
 db&=\left(\int_0^{z^2}\partial_1\varrho(z^1,u)\dd u\right)dz^1
       +\varrho(z)\,dz^2,\\
 \det\frac{\partial(a,b)}{\partial(z^1,z^2)}
   &=\varrho(z)\neq0,\\
 da\wedge db&=\varrho(z)\,dz^1\wedge dz^2=d\omega.
 \end{aligned}
\]
Restrict to a neighborhood $\mathcal O\subset Y$ of $y_*$
whose $(a,b)$-image is a Euclidean ball $B_{2\varepsilon}(0)$
and whose closure lies in the coordinate rectangle. Since
$d\{\tfrac12(a\,db-b\,da)\}=da\wedge db$, the Poincare lemma gives
$\chi_{\mathrm{loc}}\in C^\infty(\mathcal O,\R)$ such that
\begin{equation}\label{eq:darboux}
 d\omega=da\wedge db,\qquad
 \omega=\tfrac12(a\,db-b\,da)+d\chi_{\mathrm{loc}},
 \qquad \chi_{\mathrm{loc}}(y_*)=0.
\end{equation}
\emph{Two constant source disks.}
Choose a smooth $\Phi_{\mathrm{ref}}\in\mathcal B_{\mathfrak f}$,
and two disjoint coordinate disks $U_1,U_2\subset M$.
On each $U_j$, use coordinates $x^1,x^2$ given by a fixed local
inverse of $\operatorname{pr}_M$. Choose smaller disks
$U_j'\Subset U_j$, $j=1,2$, and cutoff functions
\[
 \varsigma_j\in C_c^\infty(U_j,[0,1]),\qquad
 \varsigma_j=1\text{ on }\overline{U_j'}.
\]
Choose smooth local lifts
$\widetilde\Phi_{U_j}:U_j\to\R^2$ as in
\eqref{eq:mapliftdomains}.
Fix $\widetilde y_*\in\operatorname{pr}_Y^{-1}(y_*)$ and define
\[
 \Phi_{\mathrm{ref},t}(x)=
 \begin{cases}
 \operatorname{pr}_Y\bigl((1-t\varsigma_j(x))\widetilde\Phi_{U_j}(x)
                   +t\varsigma_j(x)\widetilde y_*\bigr),&x\in U_j,\\
 \Phi_{\mathrm{ref}}(x),&x\notin U_1\cup U_2,
 \end{cases}
 \quad 0\leq t\leq1,
\]
where $j=1,2$ in the first case. Since
$\varsigma_j\in C_c^\infty(U_j)$, the formulas agree near
$\partial U_j$ and define a smooth homotopy. In particular,
\[
 \Phi_{\mathrm{base}}:=\Phi_{\mathrm{ref},1}
       \in\mathcal B_{\mathfrak f},\qquad
 \Phi_{\mathrm{base}}=y_*
       \quad\text{on }\overline{U_1'}\cup\overline{U_2'}.
\]

\emph{Independent changes of the constraints.}
Using these coordinates on $U_j$, set
\[
 X^{j }=(\sqrt G)_{j1}x^1+(\sqrt G)_{j2}x^2.
\]
For the frame fixed in Section~\ref{sec:conventions}, this gives
\[
 dX^{j }=\sum_{k=1}^2(\sqrt G)_{jk}\,dx^k=e^j,\qquad
 E_iX^{j }=\delta_{ij}\quad\text{on }U_j.
\]
Choose
\[
 \rho_j\in C_c^\infty(U_j'),\qquad
 0\leq\rho_j<\varepsilon,\qquad
 \int_{U_j'}\rho_j^2\dd V>0.
\]
For $\nu=(\nu_1,\nu_2)\in\R^2$, set
$\Phi_\nu=\Phi_{\mathrm{base}}$ outside $U_1'\cup U_2'$ and, on $U_j'$,
\begin{equation}\label{eq:diskinsertion}
 a\circ\Phi_\nu=\rho_j\cos(\nu_jX^{j }),\qquad
 b\circ\Phi_\nu=\rho_j\sin(\nu_jX^{j }).
\end{equation}
Since
\[
 (a\circ\Phi_\nu)^2+(b\circ\Phi_\nu)^2=\rho_j^2<\varepsilon^2,
 \qquad \Phi_\nu(U_j')\subset\mathcal O,
\]
the formulas are well defined, and
\[
 \Phi_\nu=\Phi_{\mathrm{base}}=y_*
 \quad\text{on }U_j'\setminus\operatorname{supp}\rho_j.
\]
Hence $\Phi_{\nu }$ is smooth on $M$.
Replacing these coordinates by $t(a\circ\Phi_\nu,b\circ\Phi_\nu)$, $0\leq t\leq1$,
gives a homotopy relative to the complement of the disks. Hence
$\Phi_\nu\simeq\Phi_{\mathrm{base}}\simeq\Phi_{\mathrm{ref}}$ in $Y$.
Differentiation gives
\[
 \begin{aligned}
 &d(a\circ\Phi_\nu)=\cos(\nu_jX^{j })d\rho_j
        -\nu_j\rho_j\sin(\nu_jX^{j })e^j,\\
 &d(b\circ\Phi_\nu)=\sin(\nu_jX^{j })d\rho_j
        +\nu_j\rho_j\cos(\nu_jX^{j })e^j,\\
 &(a\circ\Phi_\nu)\,d(b\circ\Phi_\nu)-(b\circ\Phi_\nu)\,d(a\circ\Phi_\nu)=\nu_j\rho_j^2e^j.
 \end{aligned}
\]
For this fixed $\nu$, define
\[
 f_j(x)=
 \begin{cases}
 \chi_{\mathrm{loc}}(\Phi_\nu(x)),&x\in U_j',\\
 0,&x\in M\setminus U_j'.
 \end{cases}
\]
Since $\chi_{\mathrm{loc}}(y_*)=0$,
\[
 \operatorname{supp}f_j\subset\operatorname{supp}\rho_j\Subset U_j',
 \qquad f_j\in C^\infty(M,\R).
\]
Extending $\rho_j$ by zero as well gives
\[
 A_{\Phi_\nu}-A_{\Phi_{\mathrm{base}}}
       =\frac12\sum_{j=1}^2\nu_j\rho_j^2e^j+d(f_1+f_2).
\]
Equation~\eqref{eq:sourceintegration}, with $\rho=1$, gives
\[
 \int_M E_i(f_j)\dd V=0.
\]
Consequently
\begin{equation}\label{eq:momentinsertion}
 \begin{aligned}
 I_i(\Phi_\nu)-I_i(\Phi_{\mathrm{base}})
 &=\frac12\sum_{j=1}^2\nu_j
       \int_{U_j'}\rho_j^2e^j(E_i)\dd V\\
 &=\frac{\nu_i}{2}\int_{U_i'}\rho_i^2\dd V,\qquad i=1,2.
 \end{aligned}
\end{equation}
For the prescribed $\eta$, choose
\[
 \nu_i=\frac{2\{V_M\eta(E_i)-I_i(\Phi_{\mathrm{base}})\}}{\int_{U_i'}\rho_i^2\dd V}
 \quad\Longrightarrow\quad
 I_i(\Phi_\nu)=V_M\eta(E_i),\qquad \cP(\Phi_\nu)=\eta.
\]

\end{proof}

\subsection{Attainment of the minimum}

\begin{proposition}[Attainment of the constrained minimum]\label{prop:attainment}
For every $\eta\in\cH^1(M)$, the minimum
$m_{\mathfrak f}(\eta)$ is attained in $\mathcal A_\eta$.
Moreover,
\begin{equation}\label{eq:minlower}
 m_{\mathfrak f}(\eta)\geq\frac{\|\eta\|_2^2}{2C_\omega^2}.
\end{equation}
\end{proposition}

\begin{proof}
Choose $\Phi_n\in\mathcal A_\eta$ with
$E(\Phi_n)\to m_{\mathfrak f}(\eta)<\infty$.
We keep the index $n $ for subsequences. Write
$(G^{ij})=G^{-1}$ and, for a scalar or vector-valued function,
$\overline v_n=V_M^{-1}\int_Mv_n\dd V$.

\emph{Torus target.}
Write
\[
 \phi_n([\widetilde x])=[B\widetilde x+v_n([\widetilde x])],\qquad v_n\in H^1(M,\R^2).
\]
Set
\[
 D_n=B\,dx+dv_n,\qquad
 D_n^k=\phi_n^*dy^k
   =\sum_{i=1}^2(B_{ki}+\partial_{x^i}v_n^k)\,dx^i,
 \quad k=1,2.
\]
Since $h=e^{2w}h_b$,
\[
 2E(\phi_n)=\|e^{w\circ\phi_n}D_n\|_{L^2(g_0,h_b)}^2.
\]
The form $B\,dx$ has constant coefficients. Integration by parts
gives $\int_M\partial_{x^i}v_n^k\,dV=0$, and therefore
\begin{equation}\label{eq:torusenergyorthogonality}
 \int_M\langle B\,dx,dv_n\rangle_{g_0,h_b}\,dV=0,
 \qquad
 \|D_n\|_{L^2(g_0,h_b)}^2
 =V_M|B\,dx|_{g_0,h_b}^2+\|dv_n\|_{L^2(g_0,h_b)}^2.
\end{equation}
By \eqref{eq:framebounds} and \eqref{eq:torusenergyorthogonality},
\[
 2E(\phi_n)\geq c_h\|D_n\|_{L^2(g_0,h_b)}^2
 =c_h\bigl(V_M|B\,dx|_{g_0,h_b}^2
                   +\|dv_n\|_{L^2(g_0,h_b)}^2\bigr).
\]
Since $H$ is positive definite, $\|dv_n\|_2\leq C $ in the Euclidean norm.
Adding $2\pi k_n$, $k_n\in\Z^2$, to $v_n$ leaves $\phi_n$ unchanged.
Choose $\overline v_n\in[0,2\pi)^2$. Then
\[
 \begin{aligned}
 \|v_n\|_2
 &\leq\|v_n-\overline v_n\|_2+\|\overline v_n\|_2\\
 &\leq C_P\|dv_n\|_2+V_M^{1/2}|\overline v_n|\leq C.
 \end{aligned}
\]
Thus $v_n$ is bounded in $H^1(M,\R^2)$. By
\eqref{eq:sobolevsubsequence}, a subsequence satisfies
\[
 v_n\rightharpoonup v\text{ in }H^1,\qquad
 v_n\to v\text{ in }L^2\text{ and almost everywhere}.
\]
Define
\[
 \phi_v([\widetilde x])=[B\widetilde x+v([\widetilde x])],\qquad D=B\,dx+dv.
\]
Then $\phi_v\in\mathcal B_{\mathfrak f}$,
$\phi_n\to\phi_v$ almost everywhere, and
$D_n\rightharpoonup D$ in $L^2(g_0,h_b)$. Moreover,
\[
 e^{w\circ\phi_n}\to e^{w\circ\phi_v}\quad\text{a.e.},\qquad
 e^{w\circ\phi_n}\leq e^{\|w\|_{L ^{\infty }(N) }}.
\]
Applying \eqref{eq:weakproductrule} gives
\begin{equation}\label{eq:torusweightedlimit}
 e^{w\circ\phi_n}D_n\rightharpoonup e^{w\circ\phi_v}D
       \quad\text{in }L^2(g_0,h_b).
\end{equation}
Weak lower semicontinuity applied to \eqref{eq:torusweightedlimit} gives
\[
 2E(\phi_v)=\|e^{w\circ\phi_v}D\|_{L^2(g_0,h_b)}^2
 \leq\liminf_{n\to\infty}
        \|e^{w\circ\phi_n}D_n\|_{L^2(g_0,h_b)}^2
 =\liminf_{n\to\infty}2E(\phi_n).
\]

To pass the constraint to the limit, write the fixed connection form as
\[
 \omega=\sum_{k=1}^2\omega_{y^k}\,dy^k.
\]
The coefficients $\omega_{y^k}$ are smooth on the compact torus
$N$ and have bounded periodic lifts. Thus
\[
 \omega_{y^k}\circ\phi_n\to\omega_{y^k}\circ\phi_v
       \quad\text{a.e.},\qquad
 |\omega_{y^k}\circ\phi_n|\leq\|\omega_{y^k}\|_\infty.
\]
Since $D_n^k\rightharpoonup D^k$ in $L^2(T^*M)$,
\eqref{eq:weakproductrule} yields
\begin{equation}\label{eq:torusconnectionlimit}
 A_{\phi_n}=\sum_{k=1}^2(\omega_{y^k}\circ\phi_n)D_n^k
 \rightharpoonup
 \sum_{k=1}^2(\omega_{y^k}\circ\phi_v)D^k=A_{\phi_v}
       \quad\text{in }L^2(T^*M).
\end{equation}
Testing \eqref{eq:torusconnectionlimit} against $e^1,e^2$ gives
\[
 \begin{aligned}
 I_i(\phi_n)-I_i(\phi_v)
 &=\int_M\langle A_{\phi_n}-A_{\phi_v},e^i\rangle_{g_0}\dd V
       \longrightarrow0,\\
 I_i(\phi_v)&=V_M\eta(E_i),\qquad i=1,2.
 \end{aligned}
\]
Hence $\phi_v\in\mathcal A_\eta$.

\emph{Cylindrical target.}
Write
\[
 \Phi_n=(r_n,\mathsf{s}_n),\qquad
 \mathsf{s}_n=\mathsf{s}_p+v_n\pmod L,\qquad
 W_n=\cosh r_n\,(\beta_p+dv_n).
\]
By \eqref{eq:cylinderenergy} and metric comparison,
\[
 \begin{gathered}
 c_h\bigl(\|dr_n\|_2^2+\|W_n\|_2^2\bigr)
       \leq2E(\Phi_n)\leq C,\\
 \Phi_n^*ds=\beta_p+dv_n=\sech r_n\,W_n,\qquad
 \|dv_n\|_2\leq\|W_n\|_2+\|\beta_p\|_2\leq C.
 \end{gathered}
\]
Replacing $v_n$ by $v_n+m_nL$, $m_n\in\Z$, preserves $\Phi_n$.
Choose $m_n$ so that $\overline v_n\in[0,L)$. Poincare's
inequality then gives
\[
 \|v_n\|_2\leq C_P\|dv_n\|_2
                    +V_M^{1/2}|\overline v_n|\leq C.
\]

It remains to bound $\overline r_n=V_M^{-1}\int_Mr_n\dd V$.
Since
\[
 \|r_n-\overline r_n\|_2\leq C_P\|dr_n\|_2\leq C,
\]
an unbounded sequence of means would have a subsequence with
$|\overline r_n|\to\infty$. For every fixed $R>0$ and all
sufficiently large $n$,
\begin{equation}\label{eq:cylinderescapeestimate}
 \begin{gathered}
 \{|r_n|\leq R\}
   \subset\{|r_n-\overline r_n|\geq|\overline r_n|-R\},\\
 \Vol\{|r_n|\leq R\}
 \leq\frac{\|r_n-\overline r_n\|_2^2}
              {(|\overline r_n|-R)^2}\longrightarrow0,\\
 \int_M\sech^2r_n\dd V
 \leq\Vol\{|r_n|\leq R\}+V_M\sech^2R.
 \end{gathered}
\end{equation}
In \eqref{eq:cylinderescapeestimate}, first let $n\to\infty $ and
then $R\to\infty$. We would get
$\|\sech r_n\|_2\to0$.

The fixed nontrivial cyclic class has $(p_1,p_2)\neq(0,0)$.
Define the vector field on $M$
\[
 X=\frac{\beta_p^\sharp}{|\beta_p|_{g_0}}
  =\frac{\displaystyle\sum_{i,j=1}^2G^{ij}p_j\partial_{x^i}}
         {\displaystyle\sqrt{\sum_{i,j=1}^2G^{ij}p_ip_j}}.
\]
Its coordinate coefficients are constant, and
\[
 \begin{gathered}
 \nabla^{g_0}X=0,\qquad |X|_{g_0}=1,\qquad
 \diver_{g_0}X=0,\\
 \beta_p(X)=|\beta_p|_{g_0}
  =\frac{L}{2\pi}\sqrt{\sum_{i,j=1}^2G^{ij}p_ip_j}>0.
 \end{gathered}
\]
Integration by parts gives
\[
 \int_Mdv_n(X)\dd V
   =-\int_Mv_n\diver_{g_0}X\dd V=0.
\]
Consequently
\begin{equation}\label{eq:noescape}
 \begin{aligned}
 0<V_M|\beta_p|_{g_0}
 &=\left|\int_M(\beta_p+dv_n)(X)\dd V\right|\\
 &=\left|\int_M\sech r_n\,W_n(X)\dd V\right|\\
 &\leq\|\sech r_n\|_2\|W_n\|_2\longrightarrow0,
 \end{aligned}
\end{equation}
a contradiction. Thus $\overline r_n$ is bounded, and
$r_n,v_n$ are bounded in $H^1(M)$.

By weak compactness and \eqref{eq:sobolevsubsequence}, a common
subsequence satisfies
\[
 \begin{gathered}
 r_n\rightharpoonup r,\quad v_n\rightharpoonup v
       \quad\text{in }H^1,\\
 r_n\to r,\quad v_n\to v
       \quad\text{in }L^2\text{ and almost everywhere},\\
 W_n\rightharpoonup W\quad\text{in }L^2(T^*M).
 \end{gathered}
\]
Set $\mathsf{s}=\mathsf{s}_p+v\pmod L$ and
$\Phi=(r,\mathsf{s})$. Then
$\Phi_n\to\Phi$ almost everywhere. Since
$0<\sech r_n\leq1$, \eqref{eq:weakproductrule} gives
\[
 \beta_p+dv_n\rightharpoonup\beta_p+dv,\qquad
 \sech r_n\,W_n\rightharpoonup\sech r\,W
       \quad\text{in }L^2(T^*M).
\]
Since $r$ is finite a.e., the weak limits give
\[
 \beta_p+dv=\sech r\,W,\qquad
 W=\cosh r\,(\beta_p+dv)\quad\text{a.e.}
\]
Thus $\Phi\in\mathcal B_{\mathfrak f}$.

By \eqref{eq:framebounds},
\[
 e^{w_Y\circ\Phi_n}\to e^{w_Y\circ\Phi}\quad\text{a.e.},\qquad
 e^{w_Y\circ\Phi_n}\leq e^{\|w\|_\infty}.
\]
Applying \eqref{eq:weakproductrule} to both components gives
\[
 e^{w_Y\circ\Phi_n}(dr_n,W_n)
 \rightharpoonup e^{w_Y\circ\Phi}(dr,W)
       \quad\text{in }L^2(T^*M)\oplus L^2(T^*M).
\]
Weak lower semicontinuity in this Hilbert space gives
\[
 \begin{aligned}
 2E(\Phi)
 &=\int_Me^{2(w_Y\circ\Phi)}
       \bigl(|dr|_{g_0}^2+\cosh^2r\,|\beta_p+dv|_{g_0}^2\bigr)\dd V\\
 &=\|e^{w_Y\circ\Phi}dr\|_2^2+\|e^{w_Y\circ\Phi}W\|_2^2\\
 &\leq\liminf_{n\to\infty}
    \bigl(\|e^{w_Y\circ\Phi_n}dr_n\|_2^2
         +\|e^{w_Y\circ\Phi_n}W_n\|_2^2\bigr)\\
 &=\liminf_{n\to\infty}2E(\Phi_n).
 \end{aligned}
\]

Finally, define the background-frame coefficients of $\omega$ by
\[
 \omega_1:=\omega(Q_1)=-\frac{\partial_sw_Y}{\cosh r},\qquad
 \omega_2:=\omega(Q_2)=\tanh r+\partial_rw_Y.
\]
Then
\begin{equation}\label{eq:cylboundedcoeff}
 \begin{gathered}
 \omega=\omega_1\,dr+\omega_2\cosh r\,ds,\qquad
 |\omega_1|+|\omega_2|\leq1+2\|dw_Y\|_{L^\infty(h_b)},\\
 A_{\Phi_n}=(\omega_1\circ\Phi_n)\,dr_n
                   +(\omega_2\circ\Phi_n)W_n,\qquad
 A_\Phi=(\omega_1\circ\Phi)\,dr+(\omega_2\circ\Phi)W.
 \end{gathered}
\end{equation}
Since $\omega_i\circ\Phi_n\to\omega_i\circ\Phi$ almost everywhere,
\eqref{eq:weakproductrule} gives
\[
 A_{\Phi_n}\rightharpoonup A_\Phi\quad\text{in }L^2(T^*M).
\]
Therefore
\[
 I_i(\Phi)=\langle A_\Phi,e^i\rangle_{L^2}
   =\lim_{n\to\infty}\langle A_{\Phi_n},e^i\rangle_{L^2}
   =V_M\eta(E_i),\qquad i=1,2,
\]
and $\Phi\in\mathcal A_\eta$.

In either model, writing $\Phi$ for the limit,
\[
 m_{\mathfrak f}(\eta)\leq E(\Phi)
 \leq\liminf_{n\to\infty}E(\Phi_n)=m_{\mathfrak f}(\eta).
\]
Thus the minimum is attained. The lower bound follows from
\eqref{eq:energyperiod}.
\end{proof}

\subsection{Approximation and period identities}

For $\Phi\in\mathcal B_{\mathfrak f}$, define
\begin{equation}\label{eq:framecomponentforms}
 \Theta_i(\Phi):=h_Y(d\Phi,P_i\circ\Phi)
              =\Phi^*\bigl(h_Y(P_i,\cdot)\bigr)
              \in L^2(T^*M),\qquad i=1,2.
\end{equation}

Orthonormality and orientation give, almost everywhere,
\begin{equation}\label{eq:weakframeexpansion}
 d\Phi(V)=\sum_{i=1}^2\Theta_i(\Phi)(V)P_i\circ\Phi,
 \qquad V\in TM,
\end{equation}
and
\begin{equation}\label{eq:weakframeidentities}
 \begin{gathered}
 |d\Phi|_{g_0,h_Y}^2=|\Theta_1(\Phi)|_{g_0}^2
                         +|\Theta_2(\Phi)|_{g_0}^2,\\
 \Phi^*dA_{h_Y}=\Theta_1(\Phi)\wedge\Theta_2(\Phi),\qquad
 A_\Phi=\sum_{i=1}^2(\omega(P_i)\circ\Phi)\Theta_i(\Phi).
 \end{gathered}
\end{equation}

\begin{lemma}[Strong approximation and pullback identities]
\label{lem:weakpullback}
For every $\Phi\in\mathcal B_{\mathfrak f}$ there are
$\Phi_n\in\mathcal B_{\mathfrak f}\cap C^\infty(M,Y)$ such that
\[
 \Phi_n\to\Phi\text{ almost everywhere},\qquad
 \Theta_i(\Phi_n)\longrightarrow\Theta_i(\Phi)
      \text{ in }L^2(T^*M),\quad i=1,2.
\]
The approximations need not lie in $\mathcal A_\eta$. For every such $\Phi$,
\begin{equation}\label{eq:weakcurvature}
 dA_\Phi=-\Phi^*(K_YdA_{h_Y})
 \quad\text{in distributions}.
\end{equation}
The two-form on the right belongs to $L^1$.
In the cylinder model,
\begin{equation}\label{eq:cylstokes}
 \int_M\Phi^*ds\wedge A_\Phi=0.
\end{equation}
In the torus model, let $\alpha^1,\alpha^2$ be the parallel
$h_b$-orthonormal coframe on $N$ fixed in
Section~\ref{subsec:torustargets}. Then
\begin{equation}\label{eq:torusstokes}
 \int_M\phi^*\alpha^1\wedge A_\phi
   =\int_M\phi^*\alpha^2\wedge A_\phi=0.
\end{equation}
\end{lemma}

\begin{proof}
\emph{Smooth approximation for the torus.}
Write $\phi=\phi_v$. Standard mollification on the flat torus gives
$v_n\in C^\infty(M,\R^2)$ with $v_n\to v$ in $H^1$.
Apply \eqref{eq:aeextraction} and relabel the subsequence by $n$,
so that
\[
 v_n\to v\quad\text{in }H^1\text{ and almost everywhere},\qquad
 \phi_n([\widetilde x])=[B\widetilde x+v_n([\widetilde x])].
\]
Then $\phi_n\in\mathcal B_{\mathfrak f}\cap C^\infty(M,N)$
and $\phi_n\to\phi$ almost everywhere. By
\eqref{eq:torusparallelpullback}, with $\beta_i$ fixed
independently of $n$,
\begin{equation}\label{eq:torusapproxcoframe}
 \phi_n^*\alpha^i
 =\beta_i+d\left(\sum_{j=1}^2(\sqrt H)_{ij}v_n^j\right)
 \longrightarrow\phi^*\alpha^i\quad\text{in }L^2,
 \qquad i=1,2.
\end{equation}
Since $e^{w\circ\phi_n}\to e^{w\circ\phi}$ almost everywhere
and $e^{w\circ\phi_n}\leq e^{\|w\|_\infty}$, applying
\eqref{eq:strongproductestimate} to \eqref{eq:torusapproxcoframe}
gives
\begin{equation}\label{eq:torusapproxframe}
 \Theta_i(\phi_n)=e^{w\circ\phi_n}\phi_n^*\alpha^i
 \longrightarrow e^{w\circ\phi}\phi^*\alpha^i
 =\Theta_i(\phi)\quad\text{in }L^2,
 \qquad i=1,2.
\end{equation}

\emph{Smooth approximation for the cylinder.}
For $\Phi=(r,\mathsf{s}_p+v\pmod L)$, the frame
\eqref{eq:cylinderframe} gives
\begin{equation}\label{eq:cylinderframecomponents}
 \Theta_1(\Phi)=e^{w_Y\circ\Phi}dr,\qquad
 \Theta_2(\Phi)=e^{w_Y\circ\Phi}\cosh r\,(\beta_p+dv).
\end{equation}
For $R>0$, set
\[
 r^{(R)}=\max\{-R,\min\{r,R\}\},\qquad
 \Phi^{(R)}=(r^{(R)},\mathsf{s}).
\]
The Sobolev chain rule and dominated convergence give
\begin{equation}\label{eq:cylindertruncation}
 \begin{gathered}
 dr^{(R)}=\mathbf1_{\{|r|<R\}}dr\quad\text{a.e.},\qquad
 |r^{(R)}|\leq|r|,\qquad \cosh r^{(R)}\leq\cosh r,\\
 (\Phi^{(R)})^*ds=\beta_p+dv,\qquad
 r^{(R)}\longrightarrow r\quad\text{in }H^1\quad(R\to\infty).
 \end{gathered}
\end{equation}
Indeed,
\[
 \|r^{(R)}-r\|_{H^1}^2
 \leq\int_{\{|r|\geq R\}}(|r|^2+|dr|_{g_0}^2)\dd V
 \longrightarrow0.
\]
Also $\Phi^{(R)}\in\mathcal B_{\mathfrak f}$. By
\eqref{eq:cylindertruncation} and \eqref{eq:cylinderframecomponents},
\[
 \Theta_i(\Phi^{(R)})\to\Theta_i(\Phi)\quad\text{a.e.},
 \qquad i=1,2.
\]
Moreover, $|w_Y|\leq\|w\|_\infty$ gives the bounds
\[
 \begin{aligned}
 |\Theta_1(\Phi^{(R)})-\Theta_1(\Phi)|_{g_0}^2
 &\leq4e^{2\|w\|_\infty}|dr|_{g_0}^2,\\
 |\Theta_2(\Phi^{(R)})-\Theta_2(\Phi)|_{g_0}^2
 &\leq4e^{2\|w\|_\infty}\cosh^2r\,|\beta_p+dv|_{g_0}^2.
 \end{aligned}
\]
These bounds are in $L^1$ and independent of $R$, so dominated convergence gives
\begin{equation}\label{eq:truncatedframeconvergence}
 \sum_{i=1}^2\|\Theta_i(\Phi^{(R)})-\Theta_i(\Phi)\|_2
 \longrightarrow0\quad(R\to\infty).
\end{equation}

For fixed $R$, apply \eqref{eq:periodicmollification} to
$r^{(R)},v$ and extract a common subsequence using
\eqref{eq:aeextraction}. This gives
$r_{R,m},v_{R,m}\in C^\infty(M,\R)$ with
\begin{equation}\label{eq:cylindermollifiedcoordinates}
 \begin{gathered}
 |r_{R,m}|\leq R,\qquad
 r_{R,m}\to r^{(R)},\quad v_{R,m}\to v
       \quad\text{in }H^1\text{ and almost everywhere},\\
 \Phi_{R,m}=(r_{R,m},\mathsf{s}_p+v_{R,m}\pmod L)
       \in\mathcal B_{\mathfrak f}\cap C^\infty(M,Y).
 \end{gathered}
\end{equation}
In \eqref{eq:strongproductestimate}, take respectively
\[
 (b_m,\mathcal U_m)=
 \begin{cases}
 \bigl(e^{w_Y\circ\Phi_{R,m}},\,dr_{R,m}\bigr),&i=1,\\
 \bigl(e^{w_Y\circ\Phi_{R,m}}\cosh r_{R,m},\,
                         \beta_p+dv_{R,m}\bigr),&i=2.
 \end{cases}
\]
By \eqref{eq:cylindermollifiedcoordinates}, $\mathcal U_m$ converges in
$L^2$, $b_m$ converges almost everywhere, and
\[
 e^{w_Y\circ\Phi_{R,m}}\leq e^{\|w\|_\infty},\qquad
 e^{w_Y\circ\Phi_{R,m}}\cosh r_{R,m}
                 \leq e^{\|w\|_\infty}\cosh R.
\]
For fixed $R$, \eqref{eq:strongproductestimate} gives
\begin{equation}\label{eq:cylindermollifiedframe}
 \Theta_i(\Phi_{R,m})\longrightarrow\Theta_i(\Phi^{(R)})
 \quad\text{in }L^2\quad(m\to\infty),\qquad i=1,2.
\end{equation}

For the diagonal sequence, set $R=n$. By
\eqref{eq:cylindermollifiedcoordinates} and
\eqref{eq:cylindermollifiedframe}, choose $m(n)\geq n$ such that
\begin{equation}\label{eq:cylinderdiagonalchoice}
 \begin{aligned}
 &\|r_{n,m(n)}-r^{(n)}\|_{H^1}
  +\|v_{n,m(n)}-v\|_{H^1}\\
 &\qquad+\sum_{i=1}^2
   \|\Theta_i(\Phi_{n,m(n)})-\Theta_i(\Phi^{(n)})\|_2<\frac1n.
 \end{aligned}
\end{equation}
Set $r_n=r_{n,m(n)}$, $v_n=v_{n,m(n)}$, and
$\Phi_n=\Phi_{n,m(n)}$. Equations~\eqref{eq:cylindertruncation},
\eqref{eq:truncatedframeconvergence}, and
\eqref{eq:cylinderdiagonalchoice} give
\begin{equation}\label{eq:cylinderapproxframe}
 \begin{gathered}
 r_n\to r,\quad v_n\to v\quad\text{in }H^1,\\
 \Theta_i(\Phi_n)\to\Theta_i(\Phi)\quad\text{in }L^2,
 \qquad i=1,2.
 \end{gathered}
\end{equation}
In particular,
\begin{equation}\label{eq:cylinderapproxangular}
 \Phi_n^*ds=\beta_p+dv_n\longrightarrow\beta_p+dv=\Phi^*ds
 \quad\text{in }L^2.
\end{equation}
After a further subsequence, \eqref{eq:aeextraction} also gives
$\Phi_n\to\Phi$ a.e.

\emph{Passing the curvature formula to the limit.}
By \eqref{eq:framebounds}, the smooth functions $\omega(P_i)$
and $K_Y$ are bounded on $Y$. Apply
\eqref{eq:strongproductestimate} with
$b_n=\omega(P_i)\circ\Phi_n$ and $\mathcal U_n=\Theta_i(\Phi_n)$,
using \eqref{eq:torusapproxframe} or
\eqref{eq:cylinderapproxframe}. The expression for $A_\Phi$ in
\eqref{eq:weakframeidentities} then gives
\begin{equation}\label{eq:approxconnectionstrong}
 A_{\Phi_n}
 =\sum_{i=1}^2(\omega(P_i)\circ\Phi_n)\Theta_i(\Phi_n)
 \longrightarrow A_\Phi\quad\text{in }L^2.
\end{equation}
Similarly, \eqref{eq:strongwedgeestimate},
\eqref{eq:torusapproxframe}, and \eqref{eq:cylinderapproxframe}
give
\begin{equation}\label{eq:approxareastrong}
 \Theta_1(\Phi_n)\wedge\Theta_2(\Phi_n)
 \longrightarrow\Theta_1(\Phi)\wedge\Theta_2(\Phi)
 \quad\text{in }L^1.
\end{equation}
By \eqref{eq:weakframeidentities} and \eqref{eq:approxareastrong},
\begin{equation}\label{eq:approxcurvaturestrong}
 \begin{aligned}
 &\|\Phi_n^*(K_YdA_{h_Y})-\Phi^*(K_YdA_{h_Y})\|_1\\
 &\quad\leq\|K_Y\|_\infty
   \|\Theta_1(\Phi_n)\wedge\Theta_2(\Phi_n)
                -\Theta_1(\Phi)\wedge\Theta_2(\Phi)\|_1\\
 &\qquad+\int_M|K_Y\circ\Phi_n-K_Y\circ\Phi|\,
                |\Theta_1(\Phi)\wedge\Theta_2(\Phi)|_{g_0}\dd V
 \longrightarrow0.
 \end{aligned}
\end{equation}
Here dominated convergence applies with the bound
\[
 2\|K_Y\|_\infty|\Theta_1(\Phi)|_{g_0}|\Theta_2(\Phi)|_{g_0}
 \in L^1.
\]
For smooth $\Phi_n$, the ordinary pullback identity gives
\[
 dA_{\Phi_n}=d(\Phi_n^*\omega)
       =\Phi_n^*(d\omega)=-\Phi_n^*(K_YdA_{h_Y}).
\]
By \eqref{eq:approxconnectionstrong} and
\eqref{eq:approxcurvaturestrong}, for every
$\rho\in C^\infty(M,\R)$,
\[
 \begin{aligned}
 -\int_Md\rho\wedge A_\Phi
 &=\lim_{n\to\infty}\left(-\int_Md\rho\wedge A_{\Phi_n}\right)\\
 &=-\lim_{n\to\infty}\int_M\rho\,\Phi_n^*(K_YdA_{h_Y})\\
 &=-\int_M\rho\,\Phi^*(K_YdA_{h_Y}).
 \end{aligned}
\]
This proves \eqref{eq:weakcurvature}.

\emph{Passing the Stokes identities to the limit.}
For the torus, \eqref{eq:torusexact} gives
\[
 \begin{aligned}
 \phi_n^*\alpha^1\wedge A_{\phi_n}
 &=d\bigl((w\circ\phi_n)\phi_n^*\alpha^2\bigr),\\
 \phi_n^*\alpha^2\wedge A_{\phi_n}
 &=-d\bigl((w\circ\phi_n)\phi_n^*\alpha^1\bigr).
 \end{aligned}
\]
For the cylinder, define the smooth function
\[
 F(r,s)=\int_0^r\frac{\partial_sw_Y(t,s)}{\cosh t}\dd t,
       \qquad F(r,s+L)=F(r,s).
\]
Then $F\,ds$ is a global one-form on $Y$, and
\eqref{eq:cylconnection} gives
\[
 \begin{aligned}
 ds\wedge\omega
 &=\frac{\partial_sw_Y}{\cosh r}\,dr\wedge ds
   =\partial_rF\,dr\wedge ds=d(F\,ds),\\
 \Phi_n^*ds\wedge A_{\Phi_n}
 &=d\bigl((F\circ\Phi_n)\Phi_n^*ds\bigr).
 \end{aligned}
\]
Since $\partial M=\varnothing$, Stokes' theorem yields, in the
respective models,
\begin{equation}\label{eq:smoothapproxstokes}
 \int_M\phi_n^*\alpha^i\wedge A_{\phi_n}=0\quad(i=1,2),
 \qquad
 \int_M\Phi_n^*ds\wedge A_{\Phi_n}=0.
\end{equation}
Apply \eqref{eq:strongwedgeestimate} to
\eqref{eq:torusapproxcoframe} and \eqref{eq:approxconnectionstrong}
in the torus case, and to \eqref{eq:cylinderapproxangular} and
\eqref{eq:approxconnectionstrong} in the cylinder case. This gives
\[
 \begin{aligned}
 \phi_n^*\alpha^i\wedge A_{\phi_n}
 &\longrightarrow\phi^*\alpha^i\wedge A_\phi
       &&\text{in }L^1,\quad i=1,2,\\
 \Phi_n^*ds\wedge A_{\Phi_n}
 &\longrightarrow\Phi^*ds\wedge A_\Phi
       &&\text{in }L^1.
 \end{aligned}
\]
Passing to the limit in \eqref{eq:smoothapproxstokes} proves
\eqref{eq:torusstokes} and \eqref{eq:cylstokes}.
\end{proof}

\begin{lemma}[Period restrictions when $dA_\Phi=0$]\label{lem:stokesobstruction}
Let $\Phi\in\mathcal B_{\mathfrak f}$. Suppose $dA_\Phi=0$ in
distributions and $\cP(\Phi)=\eta$. Then in the cylinder model
\begin{equation}\label{eq:cylperiodline}
 \int_M\beta_p\wedge\eta=0,\qquad \eta\in\R\beta_p.
\end{equation}
In the torus model, we have
\begin{equation}\label{eq:torusperiodline}
 \int_M\beta_1\wedge\eta=\int_M\beta_2\wedge\eta=0.
\end{equation}
Thus $\eta\in\Span\{\beta_1,\beta_2\}$ when $\rank B=1$,
and $\eta=0$ when $\rank B=2$.
\end{lemma}

\begin{proof}
The $L^2$ Hodge decomposition gives $A_\Phi=\eta+df$ with
$f\in H^1(M,\R)$. For a parallel one-form $\beta$ and
$f_{\mathrm{aux}}\in H^1(M,\R)$, smooth approximation, Stokes' theorem,
and \eqref{eq:strongwedgeestimate} give
\begin{equation}\label{eq:weakexactintegrals}
 \int_M\beta\wedge df_{\mathrm{aux}}
 =\int_Mdf_{\mathrm{aux}}\wedge df=0.
\end{equation}
Using \eqref{eq:weakcylinder}, \eqref{eq:torusparallelpullback},
and the Stokes identities of Lemma~\ref{lem:weakpullback}, we obtain
\[
 \begin{aligned}
 0&=\int_M\Phi^*ds\wedge A_\Phi
   =\int_M\beta_p\wedge\eta&&\text{(cylinder)},\\
 0&=\int_M\phi^*\alpha^i\wedge A_\phi
   =\int_M\beta_i\wedge\eta,\quad i=1,2&&\text{(torus)}.
 \end{aligned}
\]
For parallel $\beta,\eta$,
\begin{equation}\label{eq:parallelwedgecriterion}
 \int_M\beta\wedge\eta
 =V_M\{\beta(E_1)\eta(E_2)-\beta(E_2)\eta(E_1)\}.
\end{equation}
If $\beta\neq0$, this vanishes exactly when $\eta\in\R\beta$.
Now use $\beta_p\neq0$ and
$\dim\Span\{\beta_1,\beta_2\}=\rank B$.
\end{proof}

For nonzero degree, the next result applies to every value of
$\cP(\phi_v)$.

\begin{lemma}[Degree identity at finite energy]\label{lem:degree}
In the torus model, for $\phi_v\in\mathcal B_{\mathfrak f}$,
define its signed Jacobian relative to $g_0$ and $h$ by
\begin{equation}\label{eq:signedjacobiandef}
 \operatorname{Jac}(\phi_v)(x)
 :=(dA_h)_{\phi_v(x)}
       \bigl((d\phi_v)_xE_1(x),(d\phi_v)_xE_2(x)\bigr)
       \quad\text{for a.e. }x\in M.
\end{equation}
Then, for every smooth two-form $\Omega_{N }$ on $N$,
\begin{equation}\label{eq:weakdegree}
 \int_M\phi_v^*\Omega_{N }=(\det B)\int_N\Omega_{N } .
\end{equation}
If $\det B\neq0$ and $K_h\not\equiv0$, then
\begin{equation}\label{eq:Ksquaredegree}
 \int_M K_h(\phi_v)^2\operatorname{Jac}(\phi_v)\dd V
     =(\det B)\int_N K_h^2\dd A_h\neq0.
\end{equation}
Consequently $K_h(\phi_v)\neq0$ and $\rank d\phi_v=2$
on a set of positive measure.
\end{lemma}

\begin{proof}
The oriented area forms and \eqref{eq:toruscoordinatepullback} give
\[
 \begin{gathered}
 dV=\sqrt{\det G}\,dx^1\wedge dx^2,\qquad
 dA_h=e^{2w}\sqrt{\det H}\,dy^1\wedge dy^2,\\
 \mathcal J:=\det\begin{pmatrix}
 B_{11}+\partial_{x^1}v^1&B_{12}+\partial_{x^2}v^1\\
 B_{21}+\partial_{x^1}v^2&B_{22}+\partial_{x^2}v^2
 \end{pmatrix},\\
 \phi_v^*dA_h=e^{2(w\circ\phi_v)}\sqrt{\det H}\,
           \mathcal J\,dx^1\wedge dx^2
       =\operatorname{Jac}(\phi_v)dV.
 \end{gathered}
\]
Thus
\begin{equation}\label{eq:metricjacobianformula}
 \operatorname{Jac}(\phi_v)
 =e^{2(w\circ\phi_v)}\sqrt{\frac{\det H}{\det G}}\,\mathcal J
 \quad\text{a.e.}
\end{equation}
Moreover, \eqref{eq:weakframeidentities} gives
$|\operatorname{Jac}(\phi_v)|\leq\tfrac12|d\phi_v|_{g_0,h}^2\in L^1$.

Write $\Omega_{N }=F\,dA_h$, where $F\in C^\infty(N)$, and choose
the smooth approximations $\phi_n$ from Lemma~\ref{lem:weakpullback}.
Set $\Xi_n=\Theta_1(\phi_n)\wedge\Theta_2(\phi_n)$ and
$\Xi=\Theta_1(\phi_v)\wedge\Theta_2(\phi_v)$.
By \eqref{eq:strongwedgeestimate}, $\Xi_n\to\Xi$ in $L^1$.
Since $\phi_n\to\phi_v$ a.e. and $F$ is bounded,
\[
 \begin{aligned}
 \|\phi_n^*\Omega_{N }- \phi_v^*\Omega_{N }\|_1
 &\leq\|F\|_\infty\|\Xi_n-\Xi\|_1
   +\int_M|F\circ\phi_n-F\circ\phi_v|\,|\Xi|\,dV
 \longrightarrow0.
 \end{aligned}
\]
The last integral tends to zero by dominated convergence.
Every $\phi_n$ has degree $\det B$ by \eqref{eq:torusdegree},
so the smooth degree formula \cite[Theorem~17.35(a)]{LeeSM} gives
\begin{equation}\label{eq:smoothdegreeformula}
 \int_M\phi_n^*\Omega_{N }=(\det B)\int_N\Omega_{N }.
\end{equation}
Passing to the limit proves \eqref{eq:weakdegree}. Taking
$\Omega_{N }=K_h^2dA_h$ proves \eqref{eq:Ksquaredegree}.
The last assertion follows from the nonzero integral.
\end{proof}

\section{Constrained minimizers and regularity}
\label{sec:variation}

Fix $\eta\in\cH^1(M)$ satisfying
\begin{equation}\label{eq:regularperiod}
 \begin{cases}
 \eta\notin\R\beta_p,&Y\text{ is the cyclic cylinder},\\
 \eta\notin\Span\{\beta_1,\beta_2\},&Y=N,\ \rank B=1,\\
 \eta\in\cH^1(M),&Y=N,\ \det B\neq0.
 \end{cases}
\end{equation}
For every $\Phi\in\mathcal A_\eta$,
Lemmas~\ref{lem:stokesobstruction} and \ref{lem:degree} give
\begin{equation}\label{eq:nonzeropullback}
 \Phi^*(K_YdA_{h_Y})\not\equiv0.
\end{equation}

\subsection{First variations and independent constraints}

\begin{lemma}[First variations and constraint independence]\label{lem:firstvariation}
For $\Phi\in\mathcal B_{\mathfrak f}$, recall
\begin{equation}\label{eq:constraintfunctionals}
 I_i(\Phi):=\int_M A_\Phi(E_i)\dd V,\qquad
 \cP(\Phi)=\frac1{V_M}\bigl(I_1(\Phi)e^1+I_2(\Phi)e^2\bigr).
\end{equation}
For $\mathbf{f} =(f_1,f_2)\in C^\infty(M,\R^2)$, define
\begin{equation}\label{eq:expsinglevariation}
 \begin{gathered}
 W_{\mathbf{f} }(x):=\sum_{j=1}^2f_j(x)P_j(\Phi(x))\in T_{\Phi(x)}Y,\\
 \Phi_t(x):=\exp^Y_{\Phi(x)}\bigl(tW_{\mathbf{f} }(x)\bigr)
       \quad\text{for a.e. }x\in M.
 \end{gathered}
\end{equation}
Here $P_1,P_2$ are given by \eqref{eq:cylinderframe} in the
cylinder model and by \eqref{eq:torusframe} in the torus model;
$\exp^Y$ is the exponential map of $h_Y$.
\begin{enumerate}
\item For sufficiently small $|t|$, $\Phi_t\in\mathcal B_{\mathfrak f}$,
and $E(\Phi_t),I_1(\Phi_t),I_2(\Phi_t)$ are $C^1$ in $t$, with
\begin{align}
 dE_\Phi(W_{\mathbf{f} })
 &:=\left.\frac d{dt}E(\Phi_t)\right|_{t=0}
 =\int_M\sum_{i=1}^2
 h_Y(d\Phi(E_i),\nabla^{\Phi}_{E_i}W_{\mathbf{f} })\dd V,
 \label{eq:firstE}\\
 dI_i|_\Phi(W_{\mathbf{f} })
 &:=\left.\frac d{dt}I_i(\Phi_t)\right|_{t=0}
 =\int_M(K_Y\circ\Phi)
 h_Y(W_{\mathbf{f} },Jd\Phi(E_i))\dd V,\quad i=1,2.
 \label{eq:firstI}
\end{align}
For any fixed $\mathbf{f} ^{(1)}=(f_1^{(1)},f_2^{(1)})$ and
$\mathbf{f} ^{(2)}=(f_1^{(2)},f_2^{(2)})$ in $C^\infty(M,\R^2)$,
define the variation with real parameters $t,z_1,z_2$ by
\begin{equation}\label{eq:threeparamvariation}
 \begin{gathered}
 z=(z_1,z_2)\in\R^2,\\
 \Phi_{t,z}(x)
 :=\exp^Y_{\Phi(x)}\!\left(
       tW_{\mathbf{f} }(x)+z_1W_{\mathbf{f} ^{(1)}}(x)+z_2W_{\mathbf{f} ^{(2)}}(x)\right)
       \quad\text{for a.e. }x\in M.
 \end{gathered}
\end{equation}
We have $\Phi_{t,0}=\Phi_t$ and
$\Phi_{t,z}\in\mathcal B_{\mathfrak f}$ near $(t,z)=(0,0)$, with
\[
 \begin{gathered}
 (t,z)\longmapsto
 \bigl(E(\Phi_{t,z}),I_1(\Phi_{t,z}),I_2(\Phi_{t,z})\bigr)
 \quad\text{is }C^1,\\
 \left.\partial_t\Phi_{t,z}\right|_{(0,0)}=W_{\mathbf{f} },\qquad
 \left.\partial_{z_j}\Phi_{t,z}\right|_{(0,0)}=W_{\mathbf{f} ^{(j)}},
 \quad j=1,2.
 \end{gathered}
\]
\item If $\cP(\Phi)=\eta$ satisfies \eqref{eq:regularperiod}, then
$dI_1|_\Phi,dI_2|_\Phi$ are linearly independent on the fields $W_{\mathbf{f} }$:
\begin{equation}\label{eq:constraintindependence}
 \left[
 \sum_{i=1}^2\lambda_i\,dI_i|_\Phi(W_{\mathbf{f} })=0
 \quad\text{for every }\mathbf{f} \in C^\infty(M,\R^2)
 \right]
 \quad\Longrightarrow\quad \lambda_1=\lambda_2=0,
\end{equation}
where $\lambda_1,\lambda_2\in\R$.
\end{enumerate}
\end{lemma}

\begin{proof}
\emph{Variations and first derivatives.}
It suffices to treat the three-parameter family. Set
\begin{equation}\label{eq:expvariationfamily}
 \begin{gathered}
 \mathbf{f} _{t,z}:=t\mathbf{f} +z_{1 }\mathbf{f} ^{(1)}+z_{2 }\mathbf{f} ^{(2)},\qquad
 \mathscr S_{t,z}(x,y):=
       \exp_y^Y\!\left(\sum_{j=1}^2(\mathbf{f} _{t,z})_j(x)P_j(y)\right),\\
 \Phi_{t,z}=\mathscr S_{t,z}\circ(\operatorname{id}_M,\Phi).
 \end{gathered}
\end{equation}
Write $\partial\in\{\partial_t,\partial_{z_1},\partial_{z_2}\}$
and restrict the parameters to $|t|,|z_1|,|z_2|\leq1$.
The estimates in Appendix~\ref{app:variations} give
\begin{equation}\label{eq:expweakchainrule}
 d\Phi_{t,z}(E_i)
 =d_x\mathscr S_{t,z}(x,\Phi(x))[E_i(x)]
       +d_y\mathscr S_{t,z}(x,\Phi(x))[d\Phi(E_i)]
 \quad\text{a.e.}
\end{equation}
For the lifted displacement $\Delta_{t,z}$ defined in
\eqref{eq:expdisplacement}, the global real correction terms are
\begin{equation}\label{eq:exprealcorrections}
 \begin{aligned}
 v_{t,z}&=v+\Delta_{t,z}(\,\cdot\,,\Phi),&&\text{torus},\\
 \binom{r_{t,z}}{v_{t,z}}
 &=\binom rv+\Delta_{t,z}(\,\cdot\,,\Phi),&&\text{cylinder}.
 \end{aligned}
\end{equation}
By \eqref{eq:boundedcompositionchain}, these functions belong to
$H^1$. Equations~\eqref{eq:expfamilybounds} and
\eqref{eq:expweakchainrule} give
\begin{equation}\label{eq:expenergybounds}
 |d\Phi_{t,z}|_{g_0,h_Y}\leq C(1+|d\Phi|_{g_0,h_Y}),\qquad
 E(\Phi_{t,z})\leq C\bigl(V_M+E(\Phi)\bigr).
\end{equation}
In the cylinder case this also bounds
$\|\cosh r_{t,z}(\beta_p+dv_{t,z})\|_2^2$
by $2E(\Phi_{t,z})/c_h$. Thus the real corrections preserve
$B$ or $p$, and $\Phi_{t,z}\in\mathcal B_{\mathfrak f}$.

For almost every $x$, put
\begin{equation}\label{eq:expparametervector}
 W_\partial(x):=\partial\Phi_{t,z}(x)
       =(\partial\mathscr S_{t,z})(x,\Phi(x)).
\end{equation}
The covariant chain rule and mixed-derivative identity in
Appendix~\ref{app:weakvariationcalculus} yield
\begin{equation}\label{eq:expvariationbounds}
 |W_\partial|_{h_Y}\leq C,\qquad
 |\nabla^{\Phi_{t,z}}W_\partial|_{g_0,h_Y}
       \leq C(1+|d\Phi|_{g_0,h_Y}).
\end{equation}
Here $C$ may change from line to line. It depends on
$\mathbf{f} ,\mathbf{f} ^{(1)},\mathbf{f} ^{(2)}$, the metrics, and the target frame, but not on
$\Phi$ or the parameters.
By metric compatibility and \eqref{eq:expmixedderivatives},
\begin{equation}\label{eq:energyderivativedomination}
 \begin{aligned}
 \partial\!\left(\tfrac12|d\Phi_{t,z}|_{g_0,h_Y}^2\right)
 &=\sum_{i=1}^2h_Y\bigl(d\Phi_{t,z}(E_i),
                    \nabla^{\Phi_{t,z}}_{E_i}W_\partial\bigr),\\
 \left|\partial\!\left(\tfrac12|d\Phi_{t,z}|_{g_0,h_Y}^2\right)\right|
 &\leq C(1+|d\Phi|_{g_0,h_Y}^2)\in L^1(M).
 \end{aligned}
\end{equation}
The density and its parameter derivatives are continuous a.e.
The bound \eqref{eq:energyderivativedomination} also controls the difference
quotients, so dominated convergence gives $C^1$ dependence and
\begin{equation}\label{eq:energyparameterderivative}
 \partial E(\Phi_{t,z})
 =\int_M\sum_{i=1}^2h_Y\bigl(d\Phi_{t,z}(E_i),
             \nabla^{\Phi_{t,z}}_{E_i}W_\partial\bigr)\dd V.
\end{equation}

For the connection, the bounded composition rule
\eqref{eq:boundedcompositionchain} and the covariant product rule give
\begin{equation}\label{eq:connectioncontractionderivative}
 E_i\{\omega_{\Phi_{t,z}}(W_\partial)\}
 =(\nabla^Y_{d\Phi_{t,z}(E_i)}\omega)_{\Phi_{t,z}}(W_\partial)
       +\omega_{\Phi_{t,z}}(\nabla^{\Phi_{t,z}}_{E_i}W_\partial).
\end{equation}
Together with \eqref{eq:expmixedderivatives}, this gives
\begin{equation}\label{eq:connectionparameterderivative}
 \begin{aligned}
 \partial\{A_{\Phi_{t,z}}(E_i)\}
 &=(\nabla^Y_{W_\partial}\omega)_{\Phi_{t,z}}(d\Phi_{t,z}(E_i))
       +\omega_{\Phi_{t,z}}(\nabla^{\Phi_{t,z}}_{E_i}W_\partial)\\
 &=E_i\{\omega_{\Phi_{t,z}}(W_\partial)\}
       +(d\omega)_{\Phi_{t,z}}(W_\partial,d\Phi_{t,z}(E_i)).
 \end{aligned}
\end{equation}
The bounds \eqref{eq:framebounds}, \eqref{eq:expenergybounds}, and
\eqref{eq:expvariationbounds} imply
\begin{equation}\label{eq:constraintparameterbounds}
 \begin{gathered}
 |\omega_{\Phi_{t,z}}(W_\partial)|\leq C,\\
 |\partial\{A_{\Phi_{t,z}}(E_i)\}|
       +|E_i\{\omega_{\Phi_{t,z}}(W_\partial)\}|
 \leq C(1+|d\Phi|_{g_0,h_Y})\in L^2(M).
 \end{gathered}
\end{equation}
Thus $\omega_{\Phi_{t,z}}(W_\partial)\in H^1(M)$ and, for each $i$,
\begin{equation}\label{eq:connectiondivergencezero}
 \int_M E_i\{\omega_{\Phi_{t,z}}(W_\partial)\}\,dV=0
\end{equation}
by \eqref{eq:sourceintegration}. Since $M$ is compact,
\eqref{eq:constraintparameterbounds} also gives an $L^1$ bound.
Using $d\omega=-K_YdA_{h_Y}$, the same argument gives
\begin{equation}\label{eq:constraintparameterderivative}
 \partial I_i(\Phi_{t,z})
 =\int_M(K_Y\circ\Phi_{t,z})
       h_Y(W_\partial,Jd\Phi_{t,z}(E_i))\dd V.
\end{equation}
At $(t,z)=(0,0)$, $d\exp_y^Y|_0=\Id_{T_yY}$, so the parameter
directions are $W_{\mathbf{f} },W_{\mathbf{f} ^{(1)}},W_{\mathbf{f} ^{(2)}}$.
Equations~\eqref{eq:energyparameterderivative} and
\eqref{eq:constraintparameterderivative} prove the first assertions.

\emph{Independence.}
If $dI_1|_\Phi,dI_2|_\Phi$ were dependent, \eqref{eq:firstI}
would give a nonzero parallel form $\lambda=\lambda_1e^1+\lambda_2e^2$ with
\[
 \int_M(K_Y\circ\Phi)h_Y(W_{\mathbf{f} },Jd\Phi(\lambda^\sharp))\dd V=0
 \quad\text{for every }\mathbf{f} \in C^\infty(M,\R^2).
\]
Taking $\mathbf{f} =(\rho,0)$ and $\mathbf{f} =(0,\rho)$ separately gives
\[
 \int_M\rho\,(K_Y\circ\Phi)
       h_Y(P_j\circ\Phi,Jd\Phi(\lambda^\sharp))\dd V=0
 \quad\text{for all }\rho\in C^\infty(M,\R),\quad j=1,2.
\]
The two tested functions lie in $L^2(M)$, hence
\[
 (K_Y\circ\Phi)h_Y(P_j\circ\Phi,Jd\Phi(\lambda^\sharp))=0
 \quad\text{a.e.},\quad j=1,2,
 \qquad
 (K_Y\circ\Phi)Jd\Phi(\lambda^\sharp)=0\quad\text{a.e.}
\]
Set $\mathcal{K} _\Phi:=\Phi^*(K_YdA_{h_Y})\in L^1(\Lambda^2T^*M)$.
Then
\[
 \mathcal{K} _\Phi(\lambda^\sharp,E_i)
 =(K_Y\circ\Phi)h_Y(Jd\Phi(\lambda^\sharp),d\Phi(E_i))=0,
 \qquad i=1,2,
\]
so
\[
 \lambda_2\mathcal{K} _\Phi(E_1,E_2)=0,\qquad
 \lambda_1\mathcal{K} _\Phi(E_1,E_2)=0\quad\text{a.e.}
\]
Since $\lambda\neq0$, these identities give $\mathcal{K} _\Phi=0$
almost everywhere, contradicting \eqref{eq:nonzeropullback}.
\end{proof}

\subsection{Multiplier equation and regularity}

\begin{proposition}[Smooth constrained minimizers]\label{prop:smoothmin}
For every $\eta$ satisfying \eqref{eq:regularperiod}, every
minimizer in Proposition~\ref{prop:attainment} has a smooth
representative. Denote any such representative by $\Phi_\eta$.
For this fixed $\Phi_\eta$, there is a unique parallel real
one-form $\lambda_\eta$ such that
\begin{equation}\label{eq:multiplier}
 \tau(\Phi_\eta)
     =-K_Y(\Phi_\eta)Jd\Phi_\eta(\lambda_\eta^\sharp).
\end{equation}
The lower bound \eqref{eq:minlower} holds, and
$\Phi_\eta^*(K_YdA_{h_Y})$ is nonzero on a nonempty open subset of $M$.
\end{proposition}

\begin{proof}
\emph{The multiplier equation.}
Let $\Phi\in\mathcal A_\eta$ be a minimizer. By
\eqref{eq:constraintindependence}, choose
$V_j=W_{\mathbf{f} ^{(j)}}$, $j=1,2$, such that the matrix below is invertible.
For arbitrary $\mathbf{f} \in C^\infty(M,\R^2)$, use
\eqref{eq:threeparamvariation} and set
\begin{equation}\label{eq:multiplierparameterfunctions}
 \begin{gathered}
 Q(t,z):=\binom{I_1(\Phi_{t,z})}{I_2(\Phi_{t,z})},\\
 D_zQ(0,0)=\begin{pmatrix}
 dI_1|_\Phi(V_1)&dI_1|_\Phi(V_2)\\
 dI_2|_\Phi(V_1)&dI_2|_\Phi(V_2)
 \end{pmatrix},\qquad \det D_zQ(0,0)\neq0.
 \end{gathered}
\end{equation}
The map $Q$ is $C^1$ by Lemma~\ref{lem:firstvariation}.
The implicit function theorem gives a $C^1$ map
$z:(-\varepsilon,\varepsilon)\to\R^2$ satisfying
\begin{equation}\label{eq:constraintcorrection}
 \begin{gathered}
 z(0)=0,\qquad Q(t,z(t))=Q(0,0),\qquad
 \Phi_{t,z(t)}\in\mathcal A_\eta,\\
 z'(0)=-(D_zQ(0,0))^{-1}\partial_tQ(0,0).
 \end{gathered}
\end{equation}
Minimality, the parameter chain rule, and
\eqref{eq:constraintcorrection} give
\begin{equation}\label{eq:constrainedstationarity}
 \begin{aligned}
 0&=\left.\frac d{dt}E(\Phi_{t,z(t)})\right|_{t=0}\\
  &=dE_\Phi(W_{\mathbf{f} })+\sum_{j=1}^2z_j'(0)\,dE_\Phi(V_j)\\
  &=dE_\Phi(W_{\mathbf{f} })
       -\begin{pmatrix}dE_\Phi(V_1)&dE_\Phi(V_2)\end{pmatrix}
       (D_zQ(0,0))^{-1}
       \begin{pmatrix}dI_1|_\Phi(W_{\mathbf{f} })\\dI_2|_\Phi(W_{\mathbf{f} })\end{pmatrix}.
 \end{aligned}
\end{equation}
Define the constant row vector and the parallel one-form
\begin{equation}\label{eq:multipliercoefficients}
 \begin{gathered}
 (\lambda_1,\lambda_2)
 :=\begin{pmatrix}dE_\Phi(V_1)&dE_\Phi(V_2)\end{pmatrix}
       (D_zQ(0,0))^{-1},\\
 \lambda:=\lambda_1e^1+\lambda_2e^2.
 \end{gathered}
\end{equation}
The form $\lambda$ in \eqref{eq:multipliercoefficients} is independent
of $\mathbf{f} $, and \eqref{eq:constrainedstationarity} gives
\begin{equation}\label{eq:multiplieridentity}
 dE_\Phi(W_{\mathbf{f} })=\lambda_1dI_1|_\Phi(W_{\mathbf{f} })+\lambda_2dI_2|_\Phi(W_{\mathbf{f} })
 \quad\text{for every }\mathbf{f} \in C^\infty(M,\R^2).
\end{equation}
Subtracting two identities of the form
\eqref{eq:multiplieridentity} and using
\eqref{eq:constraintindependence} proves uniqueness of $\lambda$,
independently of the choice of $\mathbf{f} ^{(1)},\mathbf{f} ^{(2)}$.
Equations \eqref{eq:firstE}--\eqref{eq:firstI} yield
\begin{equation}\label{eq:weakmultiplierexpanded}
 \int_M\sum_{i=1}^2
 h_Y(d\Phi(E_i),\nabla^\Phi_{E_i}W_{\mathbf{f} })\dd V
 =\int_M(K_Y\circ\Phi)
 h_Y(W_{\mathbf{f} },Jd\Phi(\lambda^\sharp))\dd V.
\end{equation}
For smooth $\Phi$, integration by parts gives
\[
 -\int_M h_Y(\tau(\Phi),W_{\mathbf{f} })\dd V
 =\int_M(K_Y\circ\Phi)h_Y(W_{\mathbf{f} },Jd\Phi(\lambda^\sharp))\dd V,
\]
which fixes the sign in \eqref{eq:multiplier}.

\emph{The projected weak equation.}
To prove smoothness, we first project \eqref{eq:weakmultiplierexpanded}
to Euclidean space. Fix the multiplier $\lambda$ and choose a smooth isometric embedding $\iota:(N,h)\hookrightarrow\R^m$
and put
\[
 U:=\iota\circ\pi\circ\Phi.
\]
Since $\pi:(Y,h_Y)\to(N,h)$ is a local isometry,
Lemma~\ref{lem:weakpullback} and \eqref{eq:strongproductestimate} give
\[
 U\in H^1(M,\R^m),\qquad U(x)\in\iota(N)\text{ a.e.},\qquad
 dU=d(\iota\circ\pi)_\Phi\circ d\Phi,\qquad |dU|_{g_0,\R^m}=|d\Phi|_{g_0,h_Y}.
\]
For $y\in Y$, let $\Pi(\iota(\pi(y)))$ be the Euclidean
orthogonal projection onto
\[
 d\iota_{\pi(y)}(T_{\pi(y)}N)
 =d(\iota\circ\pi)_y(T_yY).
\]
Extend $J$ by zero on the normal space:
\begin{equation}\label{eq:ambientJdefinition}
 \begin{gathered}
 \mathsf J(\iota(\pi(y))):\R^m\longrightarrow\R^m,\\
 \mathsf J(\iota(\pi(y)))
       \bigl(d\iota_{\pi(y)}(V)+V^\perp\bigr)
       :=d\iota_{\pi(y)}(JV),\\
 V\in T_{\pi(y)}N,\qquad
 V^\perp\perp d\iota_{\pi(y)}(T_{\pi(y)}N).
 \end{gathered}
\end{equation}

Set $\widehat K:=K_h\circ\iota^{-1}$ on $\iota(N)$.
The tangent--normal decomposition gives
\begin{equation}\label{eq:projectionalgebra}
 \begin{gathered}
 \Pi^{\mathsf T}=\Pi,\quad \Pi^2=\Pi,\quad
 \mathsf J^{\mathsf T}=-\mathsf J,\quad
 \Pi\mathsf J=\mathsf J\Pi=\mathsf J,\\
 \|\Pi\|_{\mathrm{op}}=1,\qquad |\mathsf J|_{\mathrm F}=\sqrt2.
 \end{gathered}
\end{equation}
Let $\mathbf b\in C^\infty(M,\R^m)$ and set
\begin{equation}\label{eq:ambienttestcoefficients}
 \begin{gathered}
 f_j:=\bigl\langle\mathbf b,
     d(\iota\circ\pi)_\Phi(P_j\circ\Phi)\bigr\rangle,
 \qquad j=1,2,\\
 W:=f_1(P_1\circ\Phi)+f_2(P_2\circ\Phi),\qquad
 d(\iota\circ\pi)_\Phi W=\Pi(U)\mathbf b.
 \end{gathered}
\end{equation}

Lemma~\ref{lem:ambienttests} gives smooth coefficients $f_{j,n}$
such that $W_n:=\sum_{j=1}^2f_{j,n}(P_j\circ\Phi)$ satisfies
\eqref{eq:ambienttestconvergence}. Passing to the limit in
\eqref{eq:weakmultiplierexpanded} by $L^2$ duality and using
\eqref{eq:ambienttestpairing} gives
\begin{equation}\label{eq:projectedtestidentity}
 \int_M\sum_{i=1}^2
 \bigl\langle E_i(U),E_i\bigl(\Pi(U)\mathbf b\bigr)\bigr\rangle\dd V
 =\int_M\widehat K(U)
 \left\langle\mathbf b,\mathsf J(U)\sum_{i=1}^2\lambda_iE_i(U)\right\rangle\dd V.
\end{equation}
Write $\Pi=\Pi(U)$. Then
\[
 \Pi E_i(U)=E_i(U),\qquad
 (E_i(\Pi))^{\mathsf T}=E_i(\Pi^{\mathsf T})=E_i(\Pi)
 \quad\text{a.e.}
\]
Consequently,
\[
 \bigl\langle E_i(U),E_i(\Pi\mathbf b)\bigr\rangle
 =\langle E_i(U),E_i(\mathbf b)\rangle
       +\langle (E_i(\Pi))E_i(U),\mathbf b\rangle.
\]
Using \eqref{eq:sourceintegration} and \eqref{eq:sourcelaplacian},
\eqref{eq:projectedtestidentity} for arbitrary $\mathbf b$ gives
the distributional equation
\begin{equation}\label{eq:projectedequation}
 -\Delta_{g_0}U
 =-\sum_{i=1}^2(E_i(\Pi))E_i(U)
       +\widehat K(U)\mathsf J(U)\sum_{i=1}^2\lambda_iE_i(U).
\end{equation}
The right-hand side is in $L^1$. Differentiating $\Pi^2=\Pi$ gives
\begin{equation}\label{eq:projectionderivative}
 (E_i(\Pi))\Pi+\Pi E_i(\Pi)=E_i(\Pi),\qquad
 \Pi(E_i(\Pi))\Pi=0,
\end{equation}
Thus \eqref{eq:projectionderivative}, together with $\Pi E _{i }(U)=E _{i }(U) $, gives
\[
 \bigl(\Pi E_i(\Pi)-(E_i(\Pi))\Pi\bigr)E_i(U)
       =-(E_i(\Pi))E_i(U).
\]
Equation \eqref{eq:projectedequation} therefore takes the form
\begin{equation}\label{eq:antisymmetric}
 \begin{gathered}
 -\Delta_{g_0}U=\sum_{i=1}^2\Omega_iE_i(U),\\
 \Omega_i:=\Pi E_i(\Pi)-(E_i(\Pi))\Pi
                 +\lambda_i\widehat K(U)\mathsf J(U),\qquad i=1,2.
 \end{gathered}
\end{equation}
The symmetry of $\Pi,E_i(\Pi)$ and skew symmetry of $\mathsf J$ give
\begin{equation}\label{eq:omegaskewsymmetry}
 \Omega_i^{\mathsf T}
 =E_i(\Pi)\Pi-\Pi E_i(\Pi)
       -\lambda_i\widehat K(U)\mathsf J(U)
 =-\Omega_i.
\end{equation}
Smoothness of $\Pi$ on the compact submanifold $\iota(N)$ and
\eqref{eq:matrixproductbounds} give
\begin{equation}\label{eq:omegacoefficientbound}
 \begin{aligned}
 |\Omega_i|_{\mathrm F}
 &\leq2|D\Pi(U)[E_i(U)]|_{\mathrm F}
       +\sqrt2\,\|K_h\|_\infty|\lambda_i|\\
 &\leq C\bigl(|E_i(U)|_{\R^m}+|\lambda_i|\bigr).
 \end{aligned}
\end{equation}
Here $\Pi$ may be smoothly extended to a tubular neighborhood
of $\iota(N)$ when taking ambient derivatives. Since $dU\in L^2$
and $\lambda$ is parallel,
\begin{equation}\label{eq:omegabound}
 \begin{gathered}
 \Omega:=\Omega_1e^1+\Omega_2e^2
       \in L^2(T^*M\otimes\mathfrak{so}(m)),\\
 |\Omega|_{g_0,\mathrm F}^2
       =|\Omega_1|_{\mathrm F}^2+|\Omega_2|_{\mathrm F}^2,\qquad
 |\Omega|_{g_0,\mathrm F}
       \leq C(|dU|_{g_0,\R^m}+|\lambda|_{g_0}).
 \end{gathered}
\end{equation}

\emph{Regularity of the projection.}
In the local coordinates of \eqref{eq:sourceframe},
\eqref{eq:sourcelaplacian} allows us to apply
\cite[Theorem~1.1]{SharpTopping}. The distributional equation
\eqref{eq:antisymmetric} with \eqref{eq:omegaskewsymmetry},
\eqref{eq:omegabound}, and zero inhomogeneous term gives
$U\in W^{2,3/2}_{\mathrm{loc}}$. Sobolev embedding and
\eqref{eq:omegabound} then give
\begin{equation}\label{eq:regularityLthree}
 \begin{gathered}
 dU\in L^6_{\mathrm{loc}},\\
 \left|\sum_{i=1}^2\Omega_iE_i(U)\right|_{\R^m}
 \leq C(|dU|_{g_0,\R^m}^2+|\lambda|_{g_0}\,|dU|_{g_0,\R^m})\in L^3_{\mathrm{loc}}.
 \end{gathered}
\end{equation}
Using \eqref{eq:regularityLthree} in \eqref{eq:antisymmetric},
interior elliptic estimates give
\begin{equation}\label{eq:regularityCone}
 U\in W^{2,3}_{\mathrm{loc}}
 \subset C^{1,\alpha}_{\mathrm{loc}},\qquad 0<\alpha<\tfrac13.
\end{equation}
Since $E_i(\Pi)=D\Pi(U)[E_i(U)]$, the right-hand side of
\eqref{eq:projectedequation} is a smooth function of $U,dU$.
Schauder estimates and repeated differentiation now give
$U\in C^\infty(M,\R^m)$ \cite{GT}. Since $\iota(N)$ is closed
and $U\in\iota(N)$ almost everywhere, continuity gives
$U(M)\subset\iota(N)$.

\emph{Regularity of $\Phi$.}
For $Y=N$, $\iota^{-1}\circ U$ is a smooth representative of
$\Phi$. In the cylinder model, write
$\Phi=(r,\mathsf{s}_p+v\pmod L)$. By
\eqref{eq:cylinderenergy}, \eqref{eq:framebounds}, and
$\cosh^2r\geq1$,
\begin{equation}\label{eq:liftderivativebound}
 \begin{gathered}
 \begin{aligned}
 |d\Phi|_{g_0,h_Y}^2
 &=e^{2(w_Y\circ\Phi)}
       \bigl(|dr|_{g_0}^2+\cosh^2r\,|\beta_p+dv|_{g_0}^2\bigr)\\
 &\geq c_h\bigl(|dr|_{g_0}^2+|\beta_p+dv|_{g_0}^2\bigr),
 \end{aligned}\\
 |dr|_{g_0}+|\beta_p+dv|_{g_0}
 \leq\sqrt{\frac2{c_h}}\,|d\Phi|_{g_0,h_Y}
 =\sqrt{\frac2{c_h}}\,|dU|_{g_0,\R^m}\in L^\infty(M).
 \end{gathered}
\end{equation}
The constant is independent of $r$. Since $\beta_p$ is parallel,
\[
 r,v\in H^1(M),\quad dr,dv\in L^\infty
 \quad\Longrightarrow\quad r,v\in W^{1,\infty}(M).
\]
Use the resulting continuous representative of $\Phi$.
The identity $\pi\circ\Phi=\iota^{-1}\circ U$, initially almost
everywhere, now holds everywhere. For each $x\in M$, choose a
smooth local inverse branch $\pi^{-1}$ with
$\pi^{-1}(\pi(\Phi(x)))=\Phi(x)$. By continuity, on a sufficiently
small neighborhood $\mathcal U$ of $x$,
\[
 \Phi|_{\mathcal U}
 =\pi^{-1}\circ\iota^{-1}\circ U|_{\mathcal U}
 \in C^\infty(\mathcal U,Y).
\]
In either case, denote the smooth representative by $\Phi_\eta$;
thus $\Phi_\eta=\Phi$ almost everywhere. Set $\lambda_\eta=\lambda$.
Equation \eqref{eq:weakmultiplierexpanded} now gives
\eqref{eq:multiplier} pointwise. The smooth two-form
$\Phi_\eta^*(K_YdA_{h_Y})$ is not identically zero by
\eqref{eq:nonzeropullback}; hence
\[
 \bigl\{x\in M:(K_Y\circ\Phi_\eta)(x)\neq0,
                  \ \rank(d\Phi_\eta)_x=2\bigr\}
\]
is nonempty and open.
\end{proof}

For the smooth minimizer, \eqref{eq:multiplieridentity} and
\eqref{eq:constraintindependence} also give
\begin{equation}\label{eq:zeromultiplier}
 \lambda_\eta=0\quad\Longleftrightarrow\quad\tau(\Phi_\eta)=0.
\end{equation}

The same regularity proof applies to any $\Phi\in\mathcal B_{\mathfrak f}$
satisfying \eqref{eq:weakmultiplierexpanded} with a fixed parallel real
form $\lambda$. No minimality, period, or spin assumption is needed.

For comparison, Wang--Xu \cite[Definition~1.1 and Theorem~1.5]{WangXu}
prove smoothness of two-dimensional weak Dirac-harmonic pairs
with $\phi\in H^1$ and ambient spinor components in $L^4$ whose
weak first derivatives lie in $L^{4/3}$.
Here \eqref{eq:antisymmetric} and \eqref{eq:omegabound} give
regularity of the map before a spinor is constructed.
Once $\cP(\Phi)\in\cH_{\sigma,\delta}$,
Proposition~\ref{prop:kernel} supplies a smooth spinor with the
required curvature force.

\section{Harmonic coverings and connection periods}\label{sec:harmonic}

By \eqref{eq:zeromultiplier}, excluding zero multipliers means excluding
harmonic minimizers. Corollary~\ref{cor:harmonicperiodcondition} gives
$\cP(\Phi)\in\cH_{0,0}$ for every nonconstant harmonic map to
either framed model. For maps of nonzero degree between tori,
we prove the stronger conclusion $\cP(\phi)=0$ in the conformal
target frame.

\subsection{Harmonic coverings}

\begin{proposition}[Harmonic covering property]\label{prop:harmoniccover}
A smooth harmonic map of nonzero degree between closed oriented
Riemannian tori is an unbranched covering, for arbitrary smooth
source and target metrics.
\end{proposition}

\begin{proof}
By \eqref{eq:conformal}, use the flat source metric $g_0$.
In the Hopf differential formulas below, $h$ denotes its
complex-bilinear extension. For any smooth map $\phi$, put
\begin{equation}\label{eq:hopfcoefficient}
 \begin{gathered}
 \phi_{z_M}:=d\phi(\partial_{z_M})
       =\tfrac12\bigl(d\phi(E_1)-i\,d\phi(E_2)\bigr),\\
 \mathcal Q_\phi=\mathfrak q\,dz_M^2,\qquad
 \mathfrak q=h(\phi_{z_M},\phi_{z_M})
 =\tfrac14\bigl(|d\phi(E_1)|_h^2-|d\phi(E_2)|_h^2
          -2i\,h(d\phi(E_1),d\phi(E_2))\bigr).
 \end{gathered}
\end{equation}
For harmonic $\phi$, the Hopf differential $\mathcal Q_\phi$ is
holomorphic \cite[Section~10.1]{JostRGA}. Since $dz_M$ in
\eqref{eq:sourceisothermal} is a global nowhere-vanishing
holomorphic one-form, $\mathfrak q$ is a holomorphic function
on the compact torus $M$ and hence constant.

If $\mathfrak q=0$, $\phi$ is a nonconstant weakly conformal
harmonic map, hence holomorphic or antiholomorphic.
The Riemann--Hurwitz formula for maps between tori
\cite{JostRS} then shows that $\phi$ is an unbranched covering.

We now assume $\mathfrak q\neq0$.

\emph{A compact minimal immersion.}
Choose a constant $\theta\in\R$ with $e^{2i\theta}\mathfrak q=|\mathfrak q|$, and put
\[
 \begin{aligned}
 \widehat E_1&=\cos\theta\,E_1+\sin\theta\,E_2,&
 \widehat E_2&=-\sin\theta\,E_1+\cos\theta\,E_2,\\
 \widehat e^1&=\cos\theta\,e^1+\sin\theta\,e^2,&
 \widehat e^2&=-\sin\theta\,e^1+\cos\theta\,e^2.
 \end{aligned}
\]
These are a parallel orthonormal frame and its dual. Since
$\widehat E_1-i\widehat E_2=e^{i\theta}(E_1-iE_2)$,
\begin{equation}\label{eq:principalharmonic}
 \begin{gathered}
 h(d\phi(\widehat E_1),d\phi(\widehat E_2))=0,\\
 |d\phi(\widehat E_1)|_h^2-|d\phi(\widehat E_2)|_h^2
       =4|\mathfrak q|>0.
 \end{gathered}
\end{equation}
To cancel the Hopf differential and make the augmented map conformal, set
\begin{equation}\label{eq:heightform}
 \begin{gathered}
 \beta_{\mathrm{ht}}=2\sqrt{|\mathfrak q|}\,\widehat e^2,
 \qquad \nabla^{g_0}\beta_{\mathrm{ht}}=0,\\
 \phi^*h+\beta_{\mathrm{ht}}\otimes\beta_{\mathrm{ht}}
       =|d\phi(\widehat E_1)|_h^2g_0,
 \qquad |d\phi(\widehat E_1)|_h^2\geq4|\mathfrak q|>0.
 \end{gathered}
\end{equation}
The local primitives of $\beta_{\mathrm{ht}}$ are harmonic.
To match the periods of $\beta_{\mathrm{ht}}$, use
\eqref{eq:toruscohomologypullback} and \eqref{eq:torusdegree}:
\begin{equation}\label{eq:heightcohomology}
 \det B=\deg\phi\neq0,\qquad
 \phi^*:H^1_{\mathrm{dR}}(N;\R)
       \xrightarrow{\;\cong\;}H^1_{\mathrm{dR}}(M;\R).
\end{equation}
Let $\alpha_{\mathrm{ht}}$ be the $h_b$-parallel representative of
$(\phi^*)^{-1}\dRclass{\beta_{\mathrm{ht}}}$. Equality of de Rham
classes gives a unique normalized function
\begin{equation}\label{eq:heightpotential}
 \mathsf t_{\mathrm{ht}}\in C^\infty(M,\R),\qquad
 d\mathsf t_{\mathrm{ht}}=\beta_{\mathrm{ht}}-\phi^*\alpha_{\mathrm{ht}},
 \qquad \int_M\mathsf t_{\mathrm{ht}}\,dV=0.
\end{equation}

On $\mathcal N:=N\times\R$, with real coordinate $\mathsf t$ and
projection $\operatorname{pr}_N(y,\mathsf t)=y$, set
\begin{equation}\label{eq:heightambientmetric}
 \begin{gathered}
 h_{\mathrm{ht}}:=\operatorname{pr}_N^*h+
       (d\mathsf t+\operatorname{pr}_N^*\alpha_{\mathrm{ht}})
          \otimes(d\mathsf t+\operatorname{pr}_N^*\alpha_{\mathrm{ht}}),\\
 F:M\longrightarrow(\mathcal N,h_{\mathrm{ht}}),\qquad
 F(x)=(\phi(x),\mathsf t_{\mathrm{ht}}(x)).
 \end{gathered}
\end{equation}
Equations~\eqref{eq:heightform} and \eqref{eq:heightpotential} imply
\begin{equation}\label{eq:minimalaugmentation}
 \begin{aligned}
 F^*h_{\mathrm{ht}}
 &=\phi^*h+(d\mathsf t_{\mathrm{ht}}+\phi^*\alpha_{\mathrm{ht}})
                   \otimes(d\mathsf t_{\mathrm{ht}}+\phi^*\alpha_{\mathrm{ht}})\\
 &=|d\phi(\widehat E_1)|_h^2g_0
       \geq4|\mathfrak q|g_0.
 \end{aligned}
\end{equation}
To verify harmonicity, choose a simply connected coordinate disk
$W\subset N$ and $\mathsf t_W\in C^\infty(W,\R)$ with
$d\mathsf t_W=\alpha_{\mathrm{ht}}|_W$. The map
\[
 \mathscr I_W(y,\mathsf t)=(y,\mathsf t+\mathsf t_W(y))
 \quad\text{satisfies}\quad
 \mathscr I_W^*(h+d\mathsf t^2)=h_{\mathrm{ht}}|_{W\times\R}.
\]
On $\phi^{-1}(W)$,
\begin{equation}\label{eq:heightlocalharmonicity}
 \begin{gathered}
 \mathscr I_W\circ F=(\phi,\mathsf t_{\mathrm{ht}}+\mathsf t_W\circ\phi),\qquad
 d(\mathsf t_{\mathrm{ht}}+\mathsf t_W\circ\phi)=\beta_{\mathrm{ht}},\\
 \Delta_{g_0}(\mathsf t_{\mathrm{ht}}+\mathsf t_W\circ\phi)
   =\sum_{i=1}^2(\nabla^{g_0}_{E_i}\beta_{\mathrm{ht}})(E_i)=0,\\
 \tau_{g_0}(\mathscr I_W\circ F)=(\tau_{g_0}(\phi),0)=0.
 \end{gathered}
\end{equation}
Since $\mathscr I_W$ is a local isometry, $\tau_{g_0}(F)=0$.
Together with \eqref{eq:minimalaugmentation}, this gives
\[
 \frac12\tau_{F^*h_{\mathrm{ht}}}(F)
       =\frac{\tau_{g_0}(F)}{2|d\phi(\widehat E_1)|_h^2}=0.
\]
Thus $F$ is a minimal immersion with compact image $S:=F(M)$.

\emph{Equality of heights over each target point.}
The metric \eqref{eq:heightambientmetric} is independent of
$\mathsf t$, so
\begin{equation}\label{eq:heighttranslations}
 \mathsf T_\ell(y,\mathsf t):=(y,\mathsf t+\ell),\qquad
 \mathsf T_\ell^*h_{\mathrm{ht}}=h_{\mathrm{ht}},\qquad \ell\in\R.
\end{equation}
Thus $\mathsf T_\ell\circ F$ is also a minimal immersion.
The pairs $(x,x')$ with $\phi(x)=\phi(x')$ form a compact subset
of $M\times M$. Define their largest height difference by
\begin{equation}\label{eq:lastheightcontact}
 \begin{gathered}
 \ell_*:=\max_{\substack{x,x'\in M\\\phi(x)=\phi(x')}}
                    \bigl(\mathsf t_{\mathrm{ht}}(x)-\mathsf t_{\mathrm{ht}}(x')\bigr)\geq0,\\
 S\cap\mathsf T_{\ell_*}S\neq\varnothing,\qquad
 S\cap\mathsf T_\ell S=\varnothing\quad(\ell>\ell_*).
 \end{gathered}
\end{equation}
Suppose $\ell_*>0$. At any contact
$F(x)=\mathsf T_{\ell_*}F(x')$, restrict the two immersions to
embedded disks. Distinct tangent planes would give a transverse
intersection, which persists after a small further translation,
contrary to \eqref{eq:lastheightcontact}. Thus the tangent planes
agree. Write the disks as graphs over the common tangent plane.
If their difference changed sign, the signs at two fixed points
would persist after a small further translation, again giving an
intersection. Hence one disk lies locally on one side of the other. The strong maximum
principle \cite[Corollary~1.27 and Chapter~7, Section~1]{ColdingMinicozzi}
implies that they coincide near the contact point.

Local coincidence makes $\{x\in M:F(x)\in\mathsf T_{\ell_*}S\}$ open.
It is nonempty by \eqref{eq:lastheightcontact} and closed because
$\mathsf T _{\ell _{* }}S $ is compact. Since $M $ is connected,
\[
 S\subset\mathsf T_{\ell_*}S.
\]
Taking the minimum height yields
\[
 \min_S\mathsf t\geq\min_{\mathsf T_{\ell_*}S}\mathsf t
        =\min_S\mathsf t+\ell_*,
\]
a contradiction. Thus $\ell_*=0$, and interchanging $x,x'$ gives
\begin{equation}\label{eq:equalfiberheights}
 \phi(x)=\phi(x')\quad\Longrightarrow\quad
 \mathsf t_{\mathrm{ht}}(x)=\mathsf t_{\mathrm{ht}}(x')\quad\Longrightarrow\quad F(x)=F(x').
\end{equation}

\emph{Local injectivity and nonvanishing Jacobian.}
For each $x\in M$, choose a neighborhood $V$ on which $F$ is an
embedding. By \eqref{eq:equalfiberheights},
\[
 x_1,x_2\in V,\quad\phi(x_1)=\phi(x_2)
 \quad\Longrightarrow\quad F(x_1)=F(x_2)
 \quad\Longrightarrow\quad x_1=x_2.
\]
Invariance of domain \cite[Theorem 2B.3]{Hatcher} makes $\phi|_V$
a homeomorphism onto an open subset of $N$. Thus $\phi$ is a local
homeomorphism. Its local degree $\deg_x\phi\in\{-1,1\}$ is locally
constant, hence constant on connected $M$. At regular points it
equals the sign of $\operatorname{Jac}(\phi)$; at critical points
the Jacobian is zero. Thus the Jacobian has one sign, and
\eqref{eq:weakdegree} gives
\[
 (\deg\phi)\operatorname{Jac}(\phi)\geq0,\qquad
 \int_M\operatorname{Jac}(\phi)\dd V
      =\deg\phi\,\Area_h(N)\neq0.
\]
The first assertion in the proof of
\cite[Proposition 2.2, p.~270]{SchoenYau}, applied with the opposite
target orientation when $\deg\phi<0$, gives
\begin{equation}\label{eq:harmonicjacobianzero}
 \operatorname{Jac}(\phi)(x)=0\quad\Longrightarrow\quad d\phi_x=0.
\end{equation}
But \eqref{eq:principalharmonic} gives
$|d\phi_x(\widehat E_1(x))|_h^2\geq4|\mathfrak q|>0$.
Therefore $\operatorname{Jac}(\phi)\neq0$ everywhere.
The local diffeomorphism $\phi$ is proper because $M$ is compact,
and its image is open and closed in connected $N$. Hence $\phi$
is a finite unbranched covering.
\end{proof}

\subsection{Vanishing in the conformal target frame}

\begin{lemma}[Vanishing connection periods]\label{lem:coverperiod}
Let $\phi:(M,g_0)\to(N,h)$ be a harmonic unbranched covering between tori,
with the target frame and $\omega$ as in
\eqref{eq:torusconnection}. Then $\cP(\phi)=0$.
\end{lemma}

\begin{proof}
Since $d\phi$ is invertible, $\phi^*h$ is a smooth metric on $M$.
Set
\begin{equation}\label{eq:harmonicmetricpath}
 h_t=(1-t)g_0+t\phi^*h,\qquad 0\leq t\leq1.
\end{equation}
Consider $\Id:(M,g_0)\to(M,h_t)$. Since $E_1,E_2$ are
$g_0$-parallel, commute, and $\phi:(M,h_1)\to(N,h)$ is a local
isometry, for $k=1,2$,
\[
 \begin{aligned}
 h_t\bigl(\tau_{g_0,h_t}(\Id),E_k\bigr)
 &=\sum_{i=1}^2\left\{
    E_i\bigl(h_t(E_i,E_k)\bigr)
       -\tfrac12E_k\bigl(h_t(E_i,E_i)\bigr)\right\}\\
 &=t\sum_{i=1}^2h_1(\nabla^{h_1}_{E_i}E_i,E_k)
  =t\,h\bigl(\tau(\phi),d\phi(E_k)\bigr)=0.
 \end{aligned}
\]
Hence
\begin{equation}\label{eq:harmonicidentitypath}
 \tau_{g_0,h_t}(\Id)=0\qquad(0\leq t\leq1).
\end{equation}

Use the target coframe \eqref{eq:torusparallelframe}.
The covering $\phi $ has constant orientation sign
$\operatorname{sgn}(\deg\phi)$.
Naturality of the Hodge star under a local isometry and
\eqref{eq:torusparallelidentities} give
\begin{equation}\label{eq:coverharmonicpullback}
 \begin{gathered}
 d(\phi^*\alpha^i)=\phi^*(d\alpha^i)=0\qquad(i=1,2),\\
 *_{h_1}(\phi^*\alpha^1)
   =\operatorname{sgn}(\deg\phi)\,\phi^*(*_h\alpha^1)
   =\operatorname{sgn}(\deg\phi)\,\phi^*\alpha^2,\\
 |\phi^*\alpha^1|_{h_1}=|\alpha^1|_h\circ\phi
       =e^{-w\circ\phi}>0.
 \end{gathered}
\end{equation}
Thus $\phi^*\alpha^1$ is $h_1$-harmonic and nonzero, so its
de Rham class is nonzero.

Define
\begin{equation}\label{eq:harmoniccoframepath}
 \nu_t:=\Harm_{h_t}(\phi^*\alpha^1),
\end{equation}
where $\Harm_{h_t}$ is the $L^2(h_t)$-orthogonal projection onto
$h_t$-harmonic one-forms. By this definition and the Hodge
decomposition of the closed form $\phi^*\alpha^1$,
\begin{equation}\label{eq:harmoniccoframeproperties}
 \begin{gathered}
 d\nu_t=d(*_{h_t}\nu_t)=0,\qquad
 \nu_t-\phi^*\alpha^1\in dC^\infty(M,\R),\\
 \dRclass{\nu_t}=\dRclass{\phi^*\alpha^1}\neq0.
 \end{gathered}
\end{equation}
Elliptic theory gives smooth dependence on $t$.
For any flat metric conformal to $h_t$, conformal invariance of
the Hodge star makes $\nu_t$ harmonic, hence parallel and nowhere
zero. The orthogonality of $\nu_t,*_{h_t}\nu_t$ gives
\begin{equation}\label{eq:harmoniccoframemetric}
 h_{b,t}:=\nu_t\otimes\nu_t
           +(*_{h_t}\nu_t)\otimes(*_{h_t}\nu_t)
       =|\nu_t|_{h_t}^2h_t.
\end{equation}
The coframe $(\nu_t,*_{h_t}\nu_t)$ is closed by
\eqref{eq:harmoniccoframeproperties}, so $h_{b,t}$ is flat.
Normalize its parallel dual frame as in \eqref{eq:torusframe}:
\begin{equation}\label{eq:metricframepath}
 P_{1,t}=\frac{\nu_t^{\sharp_{h_t}}}{|\nu_t|_{h_t}},\qquad
 P_{2,t}=\frac{(*_{h_t}\nu_t)^{\sharp_{h_t}}}{|\nu_t|_{h_t}},\qquad
 \omega^{(t) }=h_t(\nabla^{h_t}P_{1,t},P_{2,t}).
\end{equation}
The forms $\omega^{(t) }$ depend smoothly on $t$. Apply
Corollary~\ref{cor:harmonicperiodcondition} to the harmonic identity
map in \eqref{eq:harmonicidentitypath}, whose pulled-back
connection form is $\Id^*\omega^{(t) }=\omega^{(t) }$:
\begin{equation}\label{eq:discreteperiodpath}
 \Harm_{g_0}(\omega^{(t) })\in\cH_{0,0}\qquad(0\leq t\leq1).
\end{equation}
At $t=0$, $\nu_0$ and $P_{1,0},P_{2,0}$ are $g_0$-parallel,
so $\omega^{(0) }=0$. Continuity and discreteness of $\cH_{0,0}$ give
\[
 \Harm_{g_0}(\omega^{(t) })=0\qquad(0\leq t\leq1).
\]

At $t=1$, $\nu_1=\phi^*\alpha^1$. Since
$\phi:(M,h_1)\to(N,h)$ is a local isometry,
\eqref{eq:torusparallelidentities}, \eqref{eq:coverharmonicpullback},
and \eqref{eq:metricframepath} give
\begin{equation}\label{eq:metricframeendpoint}
 \begin{aligned}
 d\phi(P_{1,1})
 &=\frac{(\alpha^1)^{\sharp_h}\circ\phi}
          {|\alpha^1|_h\circ\phi}=P_1\circ\phi,\\
 d\phi(P_{2,1})
 &=\operatorname{sgn}(\deg\phi)\,
   \frac{(\alpha^2)^{\sharp_h}\circ\phi}
          {|\alpha^1|_h\circ\phi}
   =\operatorname{sgn}(\deg\phi)\,P_2\circ\phi,\\
 \omega^{(1) }(X)
 &=h\bigl(d\phi(\nabla^{h_1}_XP_{1,1}),d\phi(P_{2,1})\bigr)\\
 &=\operatorname{sgn}(\deg\phi)\,
    h\bigl(\nabla^\phi_X(P_1\circ\phi),P_2\circ\phi\bigr)
 =\operatorname{sgn}(\deg\phi)A_\phi(X).
 \end{aligned}
\end{equation}
Therefore $\cP(\phi)=\Harm_{g_0}(A_\phi)=0$.
\end{proof}

\begin{proof}[Proof of Theorem~\ref{thm:harmonicintro}]
Combine Proposition~\ref{prop:harmoniccover} with
Lemma~\ref{lem:coverperiod}. The initial change from $g$ to
$g_0$ preserves harmonicity by \eqref{eq:conformal}.
\end{proof}

\section{Coupled solutions and the twofold quotient}\label{sec:existence}

\begin{theorem}[Existence at prescribed connection periods]\label{thm:periodexistence}
Fix a nontrivial class into a nonflat closed oriented surface of
genus at least one, and its framed model from Section~\ref{sec:geometry}.
Suppose one of the following holds:
\begin{enumerate}
\item The class is cyclic, the source spin structure is nontrivial,
and $\eta\in\cH_{\sigma,\delta}$ satisfies
\eqref{eq:regularperiod}.
\item The target is a nonflat torus, $\det B\neq0$, the source
spin structure is arbitrary, and
$\eta\in\cH_{\sigma,\delta}\setminus\{0\}$.
\end{enumerate}
There is a smooth coupled Dirac-harmonic pair $(\Phi_\eta,\Psi_\eta)$
from the flat spin torus $(M,g_0)$ to the framed model, with
$[\pi\circ\Phi_\eta]=\mathfrak f$, such that
\begin{equation}\label{eq:prescribedresult}
 \cP(\Phi_\eta)=\eta,\qquad
 E(\Phi_\eta)=m_{\mathfrak f}(\eta)
       \geq\frac{\|\eta\|_{L^2(g_0)}^2}{2C_\omega^2}.
\end{equation}
Its spinor is nowhere zero and nonparallel, and its map has
rank two on an open set where the target curvature is nonzero.
\end{theorem}

\begin{proof}
Proposition~\ref{prop:smoothmin} gives a smooth minimizer
$\Phi_\eta\in\mathcal A_\eta$ and its unique multiplier
$\lambda_\eta$. Since $\mathcal A_\eta\subset\mathcal B_{\mathfrak f}$,
\eqref{eq:smoothclasshomotopy} gives
$[\pi\circ\Phi_\eta]=\mathfrak f$.
Equation~\eqref{eq:energyperiod} gives \eqref{eq:prescribedresult}.
The hypotheses, Corollary~\ref{cor:harmonicperiodcondition}, and
Theorem~\ref{thm:harmonicintro} imply
$\tau(\Phi_\eta)\not\equiv0$.
Equation~\eqref{eq:multiplier} implies that the parallel one-form
$\lambda_\eta$ is nowhere zero. Applying
Corollary~\ref{cor:reduction} with $j_\eta=-\lambda_\eta$
gives a nowhere-zero spinor $\Psi_\eta$ such that
$(\Phi_\eta,\Psi_\eta)$ is coupled.

Proposition~\ref{prop:smoothmin} gives a nonempty open set on which
$\rank d\Phi_\eta=2$ and $K_Y\circ\Phi_\eta\neq0$. Thus
\eqref{eq:pullbackconnectioncurvature} gives
$dA_{\Phi_\eta}\not\equiv0$, and \eqref{eq:diracparallelcriteria}
shows that $\Psi_\eta$ is not parallel.
\end{proof}

For cyclic classes with the trivial spin structure, we use the
following twofold quotient.

\begin{lemma}[A twofold quotient for cyclic classes]\label{lem:doublequotient}
Let $(M,g_0)$ be the marked flat torus of
Section~\ref{sec:conventions}, and let $\mathfrak f$ be a
nontrivial cyclic class into a closed oriented surface $N$ of
genus at least one.
There exist a flat torus $M'$, a two-sheeted local isometry
$q_2:M\to M'$, a nontrivial spin structure $\chi'$ on $M'$, and
a nontrivial cyclic class $\mathfrak f'$ such that
\[
 q_2^*\chi'\text{ is trivial},\qquad
 [\phi'\circ q_2]=\mathfrak f
       \quad\text{whenever }[\phi']=\mathfrak f'.
\]
Pullback by $q_2$ preserves the Dirac-harmonic equations, coupling,
rank two at a point, and the properties that the twisted spinor is
nowhere zero and nonparallel. It multiplies map energy by two.
\end{lemma}

\begin{proof}
\emph{The lattices and the induced homomorphism.}
Keep the identification $\pi_1(M)\cong\Z^2$ of
Section~\ref{sec:conventions}. For a target of genus at least two,
take $c$ and $p$ from \eqref{eq:cyclicclass}. For a torus target,
take $\mathbf c$ and $p$ from \eqref{eq:rankonefactorization},
and let $c\in\pi_1(N)$ correspond to $\mathbf c$ under
$\pi_1(N)\cong\Z^2$. In either case, for the chosen based representative,
\[
 \mathfrak f_*(k)=c^{p_1k_1+p_2k_2},\quad k\in\Z^2,\qquad
 \gcd(|p_1|,|p_2|)=1.
\]
Choose $m_1,m_2\in\Z$ with $p_1m_1+p_2m_2=1$ and set
\[
 \ell_1=2\pi\binom{m_1}{m_2},\qquad
 \ell_2=2\pi\binom{-p_2}{p_1},\qquad
 \det\begin{pmatrix}m_1&-p_2\\m_2&p_1\end{pmatrix}=1.
\]
Thus $\ell_1,\ell_2$ are an oriented basis of $\Lambda=(2\pi\Z)^2$.
Put
\[
 \Lambda'=\Z\ell_1+\Z(\ell_2/2),\qquad
 M'=\R^2/\Lambda',\qquad
 q_2([\widetilde x]_\Lambda)=[\widetilde x]_{\Lambda'}.
\]
Here $[\widetilde x]_\Lambda=\widetilde x+\Lambda$, and similarly
for $\Lambda'$. Since $[\Lambda':\Lambda]=2$, $q_2$ is two-sheeted.
The constant metric descends to $g'_0$ with
\[
 q_2^*g'_0=g_0,\qquad \Vol_{g_0}(M)=2\Vol_{g'_0}(M').
\]
Use the oriented marking $(\ell_1,\ell_2/2)$ on $M'$.
For
$k=(k_1,k_2)^{\mathsf T}\in\Z^2$,
\[
 2\pi k=(p_1k_1+p_2k_2)\ell_1
       +2(-m_2k_1+m_1k_2)\frac{\ell_2}{2}.
\]
Hence, in the respective marked bases of the two fundamental groups,
\begin{equation}\label{eq:sourcequotientmatrix}
 (q_2)_*\binom{k_1}{k_2}
 =\begin{pmatrix}p_1&p_2\\-2m_2&2m_1\end{pmatrix}
       \binom{k_1}{k_2}
 =\binom{p_1k_1+p_2k_2}{2(-m_2k_1+m_1k_2)}.
\end{equation}
Its matrix has determinant $2(p_1m_1+p_2m_2)=2$.

\emph{The spin structure and the map class.}
In the marking of $M'$, take the spin character
\[
 \chi'(k)=(-1)^{k_2},\qquad k\in\Z^2.
\]
It is the nontrivial spin structure with parameters $(0,1/2)$,
and its pullback satisfies
\[
 (q_2^*\chi')(k)
   =\chi'((q_2)_*k)=(-1)^{2(-m_2k_1+m_1k_2)}=1.
\]
Choose $\mathfrak f'$ with based homomorphism
\[
 \mathfrak f'_*(k)=c^{k_1},\qquad k\in\Z^2.
\]
The classification \eqref{eq:homotopyclassification} realizes
this homomorphism by a nontrivial cyclic class, and
\[
 \mathfrak f'_*((q_2)_*k)
 =c^{p_1k_1+p_2k_2}=\mathfrak f_*(k).
\]
By \eqref{eq:homotopyclassification},
$[\phi'\circ q_2]=\mathfrak f$ whenever $[\phi']=\mathfrak f'$.

\emph{Pullback of the equations.}
Let $(\phi',\Psi')$ be a pair on $M'$ and set $\phi=\phi'\circ q_2$.
The orientation-preserving local isometry $q_2$ and the pullback
spin structure give isometric bundle identifications
\[
 \Sigma_{g_0}M\cong q_2^*\Sigma_{g'_0}M',\qquad
 \phi^*TN=q_2^*((\phi')^*TN).
\]
For a local vector field $X'$ on $M'$, let $X$ be its lift, so
$dq_2(X)=X'\circ q_2$. For a local spinor $\xi'$ and a local
section $V'$ of $(\phi')^*TN$, these identifications give
\begin{equation}\label{eq:sourcequotientconnections}
 \begin{aligned}
 X\cdot q_2^*\xi'&=q_2^*(X'\cdot\xi'),\\
 \nabla^{\Sigma_{g_0}}_X(q_2^*\xi')
   &=q_2^*(\nabla^{\Sigma_{g'_0}}_{X'}\xi'),\\
 \nabla^\phi_X(q_2^*V')&=q_2^*(\nabla^{\phi'}_{X'}V').
 \end{aligned}
\end{equation}
Set $\Psi=q_2^*\Psi'$ under
$\cE_\phi\cong q_2^*\cE_{\phi'}$.
Equation~\eqref{eq:sourcequotientconnections} and the local isometry
property imply
\[
 \D^\phi\Psi=q_2^*(\D^{\phi'}\Psi'),\qquad
 \tau(\phi)=q_2^*\tau(\phi'),\qquad
 \cR(\phi,\Psi)=q_2^*\cR(\phi',\Psi').
\]
Moreover, under the same bundle identifications,
\[
 \begin{gathered}
 |\Psi|=|\Psi'|\circ q_2,\qquad
 |\nabla^{\cE_\phi}\Psi|_{g_0}
       =|\nabla^{\cE_{\phi'}}\Psi'|_{g'_0}\circ q_2,\\
 \rank d\phi=(\rank d\phi')\circ q_2 .
 \end{gathered}
\]
Surjectivity preserves the stated nonvanishing properties.
Integration over the two sheets gives
\begin{equation}\label{eq:quotientenergy}
 E_M(\phi'\circ q_2)
 =\tfrac12\int_M|d\phi'|_{g'_0,h}^2\circ q_2\dd V_{g _{0 }}
 =\int_{M'}|d\phi'|_{g'_0,h}^2\dd V_{g' _{0 }}=2E_{M'}(\phi').
\end{equation}
\end{proof}

\begin{proof}[Proof of Theorem~\ref{thm:main}]
Let $g_0$ be a flat metric conformal to $g$. By
Section~\ref{sec:geometry}, every nontrivial class is cyclic or
has nonzero degree. Lemma~\ref{lem:doublequotient} reduces the
cyclic case with trivial spin structure to a nontrivial spin
structure: pullback recovers the prescribed class and the required
solution properties, and multiplies energy by two, as in \eqref{eq:quotientenergy}. Thus it suffices to prove the flat-metric assertion under the
hypotheses of Theorem~\ref{thm:periodexistence}.

Fix $R>0$. The rank-two affine lattice $\cH_{\sigma,\delta}$ in
\eqref{eq:spinlattice} remains unbounded after removing the line
or the origin excluded by Theorem~\ref{thm:periodexistence}.
Choose an admissible $\eta$ with
$\|\eta\|_{L^2(g_0)}^2>2C_\omega^2R$.
Theorem~\ref{thm:periodexistence} gives a smooth coupled pair
$(\Phi,\Psi)$ satisfying
\[
 [\pi\circ\Phi]=\mathfrak f,\qquad
 E(\Phi)\geq\frac{\|\eta\|_{L^2(g_0)}^2}{2C_\omega^2}>R.
\]
Its spinor is nowhere zero, and $\rank d\Phi=2$ on a nonempty
open set where $K_Y\circ\Phi\neq0$.

Set
$\phi=\pi\circ\Phi$ and $\Psi^N=(\Id\otimes d\pi)\Psi$.
The local isometry $\pi$, with the metric
\eqref{eq:targetcovermetric}, gives
\begin{equation}\label{eq:targetprojectionnaturality}
 \begin{aligned}
 \tau(\phi)&=d\pi\,\tau(\Phi),\\
 \D^\phi\Psi^N&=(\Id\otimes d\pi)\D^\Phi\Psi,\\
 \cR(\phi,\Psi^N)&=d\pi\,\cR(\Phi,\Psi),\\
 |d\phi|_{g_0,h}&=|d\Phi|_{g_0,h_Y},\qquad
 |\Psi^N|=|\Psi|.
 \end{aligned}
\end{equation}
Since $d\pi$ is an isometry and $K_Y=K_h\circ\pi$, projection
preserves all the required properties and gives $E_{g_0}(\phi)>R$.
Together with the pullback in Lemma~\ref{lem:doublequotient}, this
proves the assertion for every flat spin torus.

For any resulting pair $(\phi,\Psi_0)$ on the original $(M,g_0)$,
write $g=e^{2\varpi}g_0$ and set
\[
 \Psi_g=e^{-\varpi/2}(\mathcal B\otimes\Id)\Psi_0.
\]
Equation~\eqref{eq:conformal} gives
\[
 \D_g^\phi\Psi_g=0,\qquad
 \tau_g(\phi)=\cR_g(\phi,\Psi_g)
       =e^{-2\varpi}\tau_{g_0}(\phi)\not\equiv0,\qquad
 E_g(\phi)=E_{g_0}(\phi)>R.
\]
The spinor remains nowhere zero, and the map, its homotopy class
and its rank are unchanged. Taking $R=n$ for $n=1,2,\ldots$ gives
the required sequence, with $E_g(\phi_n)>n$.
\end{proof}

\begin{proof}[Proof of Corollary~\ref{cor:criterion}]
Flat targets have zero curvature force. If
$\deg(\mathfrak f)\chi(N)\neq0$, the genus--degree hypothesis of
\cite[Theorem~3.3]{CJSZsurfaces} becomes
$0<|\deg(\mathfrak f)\chi(N)|$ for a torus source, so every
Dirac-harmonic pair in $\mathfrak f$ is uncoupled. This proves necessity.
For sufficiency, Theorem~\ref{thm:main} treats every nontrivial
class into a target of genus at least one, while
\cite[Theorem 1.1]{JZcurves} treats every null-homotopic class.
For a sphere target, \eqref{eq:criterion} gives degree zero,
which is precisely the null-homotopic class because
$[M,S^2]\cong H^2(M;\Z)\cong\Z$ via degree
\cite[Section 4.3]{Hatcher}. Both existence results give unbounded
map energies.
\end{proof}

\appendix
\section{Weak variations and ambient test fields}\label{app:variations}

In Subsections~\ref{app:exponentialbounds}--\ref{app:weakvariationcalculus},
we prove the estimates and weak chain rules used in
Lemma~\ref{lem:firstvariation}. In these subsections, $\mathscr S_{t,z}$ and $\mathbf{f} _{t,z}$
are defined by \eqref{eq:expvariationfamily},
$\partial\in\{\partial_t,\partial_{z_1},\partial_{z_2}\}$, and
$|t|,|z_1|,|z_2|\leq1$. Constants depend on
$\mathbf{f} ,\mathbf{f} ^{(1)},\mathbf{f} ^{(2)}$, the metrics, and the fixed frame; they may
change from line to line.

Subsection~\ref{app:ambienttests} treats the ambient test fields
used in Proposition~\ref{prop:smoothmin}.

\subsection{Uniform derivatives of the exponential family}
\label{app:exponentialbounds}

The Riemannian covering $\pi:(Y,h_Y)\to(N,h)$ makes $Y$ complete,
because geodesics in the closed surface $N$ lift to geodesics in $Y$.
Projection of geodesics gives
\begin{equation}\label{eq:expcoveridentity}
 \pi(\exp_y^Y V)=\exp_{\pi(y)}^N(d\pi_yV),\qquad
 |d\pi_yV|_h=|V|_{h_Y}.
\end{equation}
For fixed $R<\infty$, the disk bundle
$\mathcal K_R=\{V\in TN:|V|_h\leq R\}$ is compact.
The first two derivatives of $\exp^N$ are bounded in finitely many
bundle charts covering $\mathcal K_R$. By \eqref{eq:expcoveridentity},
these bounds lift uniformly to $Y$.

Set
\begin{equation}\label{eq:expframecoordinates}
 \mathscr X(y,\mathbf q)
 :=\exp_y^Y\bigl(q_1P_1(y)+q_2P_2(y)\bigr),\qquad
 \mathbf q\in\R^2.
\end{equation}
Use the product connection on $Y\times\R^2$ and the pullback
connection $\nabla^{\mathscr X}=\mathscr X^*\nabla^Y$ on
$\mathscr X^*TY$. The exponential estimates and
\eqref{eq:framebounds} give, for $|\mathbf q|\leq R$,
\begin{equation}\label{eq:expframebounds}
 \|d_y\mathscr X\|_{\mathrm{op}}
 +\|d_{\mathbf q}\mathscr X\|_{\mathrm{op}}
 +\|\nabla_{\mathbf q}d_{\mathbf q}\mathscr X\|_{\mathrm{op}}
 +\|\nabla_y d_{\mathbf q}\mathscr X\|_{\mathrm{op}}\leq C_R.
\end{equation}
All tangent norms use $h_Y$, and coefficient norms are Euclidean.

Choose $R>\|\mathbf{f} \|_\infty+\|\mathbf{f} ^{(1)}\|_\infty+\|\mathbf{f} ^{(2)}\|_\infty$.
Then
\begin{equation}\label{eq:variationcoefficientbounds}
 |\mathbf{f} _{t,z}|<R,\qquad
 |d\mathbf{f} _{t,z}|+|\partial \mathbf{f} _{t,z}|+|d(\partial \mathbf{f} _{t,z})|\leq C.
\end{equation}
Since $\mathscr S_{t,z}(x,y)=\mathscr X(y,\mathbf{f} _{t,z}(x))$,
\begin{equation}\label{eq:expfamilyparameterderivative}
 \begin{aligned}
 d_x\mathscr S_{t,z}[X]&=d_{\mathbf q}\mathscr X[d\mathbf{f} _{t,z}(X)],\\
 d_y\mathscr S_{t,z}[V]&=d_y\mathscr X[V],\qquad
 \partial\mathscr S_{t,z}=d_{\mathbf q}\mathscr X[\partial \mathbf{f} _{t,z}],
 \end{aligned}
\end{equation}
where $\mathscr X$ is evaluated at $(y,\mathbf{f} _{t,z}(x))$.
Covariant differentiation with $\mathscr S_{t,z}^*\nabla^Y$ yields
\begin{equation}\label{eq:expfamilymixedcalculation}
 \begin{aligned}
 (\nabla_x\partial\mathscr S_{t,z})[X]
 &=(\nabla_{\mathbf q}d_{\mathbf q}\mathscr X)
       [d\mathbf{f} _{t,z}(X),\partial \mathbf{f} _{t,z}]
       +d_{\mathbf q}\mathscr X[d(\partial \mathbf{f} _{t,z})(X)],\\
 (\nabla_y\partial\mathscr S_{t,z})[V]
 &=(\nabla_y d_{\mathbf q}\mathscr X)[V,\partial \mathbf{f} _{t,z}].
 \end{aligned}
\end{equation}
Equations~\eqref{eq:expframebounds}--\eqref{eq:expfamilymixedcalculation}
give
\begin{equation}\label{eq:expfamilybounds}
 \begin{gathered}
 |\partial\mathscr S_{t,z}|_{h_Y}
 +\|d_x\mathscr S_{t,z}\|_{\mathrm{op}}
 +\|d_y\mathscr S_{t,z}\|_{\mathrm{op}}\leq C,\\
 \|\nabla_x\partial\mathscr S_{t,z}\|_{\mathrm{op}}
 +\|\nabla_y\partial\mathscr S_{t,z}\|_{\mathrm{op}}\leq C.
 \end{gathered}
\end{equation}
The input norms are $g_0$ and $h_Y$, respectively, and the output
norm is $h_Y$.

\subsection{Equivariant lifts and displacement}\label{app:displacement}

Equip $\widetilde Y=\R^2$ with $\operatorname{pr}_Y^*h_Y$.
For $y=\operatorname{pr}_Y(\widetilde y)$, define
\begin{equation}\label{eq:expvariationlift}
 \begin{gathered}
 \widetilde{\mathscr S}_{t,z}(x,\widetilde y)
 :=\exp_{\widetilde y}^{\widetilde Y}\!\left(
       \bigl((d\operatorname{pr}_Y)_{\widetilde y}\bigr)^{-1}
       \left[\sum_{j=1}^2(\mathbf{f} _{t,z})_j(x)P_j(y)\right]\right),\\
 \operatorname{pr}_Y\circ\widetilde{\mathscr S}_{t,z}(x,\cdot)
       =\mathscr S_{t,z}(x,\cdot)\circ\operatorname{pr}_Y.
 \end{gathered}
\end{equation}
For $\ell\in\Lambda_Y$ from \eqref{eq:targetlattice}, the translation
$T_\ell(\widetilde y)=\widetilde y+\ell$ satisfies
\[
 \operatorname{pr}_Y\circ T_\ell=\operatorname{pr}_Y,\qquad
 T_\ell^*(\operatorname{pr}_Y^*h_Y)=\operatorname{pr}_Y^*h_Y.
\]
Geodesic uniqueness gives
\begin{equation}\label{eq:exponentialdeckidentity}
 \exp^{\widetilde Y}_{\widetilde y+\ell}((dT_\ell)_{\widetilde y}V)
 =T_\ell(\exp^{\widetilde Y}_{\widetilde y}V).
\end{equation}
By \eqref{eq:exponentialdeckidentity},
\begin{equation}\label{eq:expdisplacement}
 \begin{gathered}
 \widetilde{\mathscr S}_{t,z}(x,\widetilde y+\ell)
 =\widetilde{\mathscr S}_{t,z}(x,\widetilde y)+\ell,\\
 \Delta_{t,z}(x,\operatorname{pr}_Y(\widetilde y))
 :=\widetilde{\mathscr S}_{t,z}(x,\widetilde y)-\widetilde y.
 \end{gathered}
\end{equation}
Thus $\Delta_{t,z}:M\times Y\to\R^2$ is well defined and smooth.
The coordinate covectors satisfy
\begin{equation}\label{eq:coordinatecovectorbounds}
 \begin{aligned}
 |dr|_{h_Y}^2+|ds|_{h_Y}^2
 &=e^{-2w_Y}(1+\sech^2r)\leq C&&\text{(cylinder)},\\
 |dy^1|_h^2+|dy^2|_h^2
 &=e^{-2w}\tr(H^{-1})\leq C&&\text{(torus)}.
 \end{aligned}
\end{equation}
Consequently \eqref{eq:coverinversedifferential} gives
\begin{equation}\label{eq:inversebranchbound}
 |d(\operatorname{pr}_Y^{-1})_y[V]|_{\R^2}\leq C|V|_{h_Y}.
\end{equation}
The curve $\tau\mapsto\mathscr S_{\tau t,\tau z}(x,y)$,
$0\leq\tau\leq1$, is a geodesic of constant speed
$|\mathbf{f} _{t,z}(x)|$. Integrating
\eqref{eq:inversebranchbound} along its lift gives
\begin{equation}\label{eq:expdisplacementvaluebound}
 |\Delta_{t,z}(x,y)|_{\R^2}
 \leq C|\mathbf{f} _{t,z}(x)|\leq C(|t|+|z_1|+|z_2|).
\end{equation}
Choose local inverse branches through $\widetilde y$ and
$\widetilde{\mathscr S}_{t,z}(x,\widetilde y)$.
Differentiating \eqref{eq:expvariationlift} gives
\begin{equation}\label{eq:expdisplacementdifferentials}
 \begin{aligned}
 d_x\Delta_{t,z}
 &=d(\operatorname{pr}_Y^{-1})_{\mathscr S_{t,z}(x,y)}
       \circ d_x\mathscr S_{t,z},\\
 d_y\Delta_{t,z}
 &=d(\operatorname{pr}_Y^{-1})_{\mathscr S_{t,z}(x,y)}
       \circ d_y\mathscr S_{t,z}-d(\operatorname{pr}_Y^{-1})_y.
 \end{aligned}
\end{equation}
By \eqref{eq:coverinversedifferential}, the inverse differentials
do not depend on the branch in real coordinates. Equations
\eqref{eq:expfamilybounds} and \eqref{eq:inversebranchbound} imply
\begin{equation}\label{eq:expdisplacementderivativebounds}
 \|d_x\Delta_{t,z}\|_{\mathrm{op}}
       +\|d_y\Delta_{t,z}\|_{\mathrm{op}}\leq C.
\end{equation}

\subsection{Weak chain rules}\label{app:weakvariationcalculus}

Let
$F:M\times Y\to\R^m$ be smooth, with
$|F|_{\R^m}+\|d_xF\|_{\mathrm{op}}+\|d_yF\|_{\mathrm{op}}\leq C$.
Choose the smooth maps from Lemma~\ref{lem:weakpullback}:
$\Phi_n\to\Phi$ almost everywhere and
$\Theta_j(\Phi_n)\to\Theta_j(\Phi)$ in $L^2$.
For $i=1,2$, the smooth chain rule is
\[
 E_i\{F(\,\cdot\,,\Phi_n)\}
 =d_xF(\,\cdot\,,\Phi_n)[E_i]
  +\sum_{j=1}^2d_yF(\,\cdot\,,\Phi_n)[P_j\circ\Phi_n]
       \Theta_j(\Phi_n)(E_i).
\]
Boundedness of the coefficients, dominated convergence, and
\eqref{eq:strongproductestimate} give
\begin{equation}\label{eq:boundedcompositionconvergence}
 \begin{gathered}
 F(\,\cdot\,,\Phi_n)\longrightarrow F(\,\cdot\,,\Phi)
       \quad\text{in }L^2,\\
 E_i\{F(\,\cdot\,,\Phi_n)\}
 \longrightarrow d_xF(\,\cdot\,,\Phi)[E_i]
       +d_yF(\,\cdot\,,\Phi)[d\Phi(E_i)]
       \quad\text{in }L^2.
 \end{gathered}
\end{equation}
Passing to the limit in the integration-by-parts identity
\eqref{eq:sourceintegration} and using
\eqref{eq:sourceweakderivative}--\eqref{eq:sourceh1criterion} gives
\begin{equation}\label{eq:boundedcompositionchain}
 \begin{gathered}
 F(\,\cdot\,,\Phi)\in H^1(M,\R^m),\\
 E_i\{F(\,\cdot\,,\Phi)\}(x)
 =d_xF(x,\Phi(x))[E_i(x)]
       +d_yF(x,\Phi(x))[(d\Phi)_x(E_i(x))]
       \quad\text{a.e. }x\in M.
 \end{gathered}
\end{equation}

Apply \eqref{eq:boundedcompositionchain} to $F=\Delta_{t,z}$.
For a local real lift $\widetilde\Phi_U$ as in
\eqref{eq:mapliftdomains},
\begin{equation}\label{eq:expreallift}
 \widetilde\Phi_{t,z,U}
 :=\widetilde\Phi_U+\Delta_{t,z}(\,\cdot\,,\Phi)\in H^1(U,\R^2),
 \qquad \operatorname{pr}_Y\circ\widetilde\Phi_{t,z,U}
       =\Phi_{t,z}|_U.
\end{equation}
Taking weak derivatives gives
\[
 d\widetilde\Phi_{t,z,U}(E_i)
 =d\widetilde\Phi_U(E_i)+d_x\Delta_{t,z}[E_i]
       +d_y\Delta_{t,z}[d\Phi(E_i)].
\]
Substitute \eqref{eq:expdisplacementdifferentials} and
\eqref{eq:mapliftdifferentials}. Multiplication by
$(d\operatorname{pr}_Y)_{\widetilde\Phi_{t,z,U}}$ then proves
\eqref{eq:expweakchainrule}. Substituting the coordinate expressions
for $\widetilde\Phi_U$ in \eqref{eq:expreallift} gives
\eqref{eq:exprealcorrections}.

For a smooth section $\mathscr V$ of $\mathscr S_{t,z}^*TY$ with
bounded value and first covariant derivatives, set
$F_j=h_Y(\mathscr V,P_j\circ\mathscr S_{t,z})$.
Equations~\eqref{eq:expfamilybounds} and \eqref{eq:framebounds}
bound $F_j,d_xF_j,d_yF_j$. Applying
\eqref{eq:boundedcompositionchain} to these coefficients gives
\begin{equation}\label{eq:covariantcompositionchain}
 \begin{aligned}
 \nabla^{\Phi_{t,z}}_{E_i}\{\mathscr V(\,\cdot\,,\Phi)\}
 &=\sum_{j=1}^2\bigl(d_xF_j[E_i]+d_yF_j[d\Phi(E_i)]\bigr)
                         P_j\circ\Phi_{t,z}\\
 &\quad+\sum_{j=1}^2F_j(\,\cdot\,,\Phi)
       \nabla^Y_{d_x\mathscr S_{t,z}[E_i]+d_y\mathscr S_{t,z}[d\Phi(E_i)]}P_j\\
 &=(\nabla_x\mathscr V)(x,\Phi(x))[E_i(x)]
       +(\nabla_y\mathscr V)(x,\Phi(x))[d\Phi(E_i)].
 \end{aligned}
\end{equation}
The middle line uses \eqref{eq:expweakchainrule}; all terms lie in
$L^2$. Taking $\mathscr V=\partial\mathscr S_{t,z}$ yields
\begin{equation}\label{eq:expparameterchainrule}
 \nabla^{\Phi_{t,z}}_{E_i}W_\partial
 =(\nabla_x\partial\mathscr S_{t,z})(x,\Phi(x))[E_i(x)]
       +(\nabla_y\partial\mathscr S_{t,z})(x,\Phi(x))[d\Phi(E_i)].
\end{equation}
On the fixed full-measure set where $\Phi$ and $d\Phi$ are defined,
use the right-hand side of \eqref{eq:expweakchainrule} for all
parameters. It is smooth in $(t,z)$. Torsion-freeness gives
\begin{equation}\label{eq:expmixedderivatives}
 \begin{aligned}
 \nabla_\partial d\Phi_{t,z}(E_i)
 &=(\nabla_\partial d_x\mathscr S_{t,z})[E_i]
       +(\nabla_\partial d_y\mathscr S_{t,z})[d\Phi(E_i)]\\
 &=(\nabla_x\partial\mathscr S_{t,z})[E_i]
       +(\nabla_y\partial\mathscr S_{t,z})[d\Phi(E_i)]
 =\nabla^{\Phi_{t,z}}_{E_i}W_\partial.
 \end{aligned}
\end{equation}

\subsection{Ambient test fields}\label{app:ambienttests}

\begin{lemma}[Approximation of ambient test fields]\label{lem:ambienttests}
Let $\iota:(N,h)\hookrightarrow\R^m$ be a smooth isometric embedding,
$\Phi\in\mathcal B_{\mathfrak f}$, and $U=\iota\circ\pi\circ\Phi$.
For $\mathbf b\in C^\infty(M,\R^m)$, define $f_j,W$ by
\eqref{eq:ambienttestcoefficients}, where $\Pi(U)$ is the Euclidean
orthogonal projection onto $T_U\iota(N)$.
Then $f_j\in H^1(M)\cap L^\infty(M)$, and there are smooth
$f_{j,n}$ such that $W_n:=\sum_{j=1}^2f_{j,n}(P_j\circ\Phi)$
satisfies $W_n\to W$ and $\nabla^\Phi W_n\to\nabla^\Phi W$ in $L^2$.
Moreover, for $i=1,2$,
\begin{equation}\label{eq:ambienttestpairing}
 \bigl\langle dU(E_i),E_i\bigl(\Pi(U)\mathbf b\bigr)\bigr\rangle
 =h_Y\bigl(d\Phi(E_i),\nabla^\Phi_{E_i}W\bigr)
 \quad\text{a.e.}
\end{equation}
\end{lemma}

\begin{proof}
Let $\mathrm{II}$ be the second fundamental form of
$\iota(N)\subset\R^m$. The Gauss formula and
\eqref{eq:framebounds} give, for $V\in T_yY$ and $j=1,2$,
\begin{equation}\label{eq:ambientframederivative}
 \begin{aligned}
 d\{d(\iota\circ\pi)(P_j)\}_y[V]
 &=d(\iota\circ\pi)_y(\nabla^Y_VP_j)\\
 &\quad+\mathrm{II}_{\iota(\pi(y))}
       \bigl(d(\iota\circ\pi)_yV,d(\iota\circ\pi)_yP_j\bigr),\\
 |d\{d(\iota\circ\pi)(P_j)\}_y[V]|_{\R^m}
 &\leq\bigl(|\nabla^YP_j|_{h_Y}+\|\mathrm{II}\|_\infty\bigr)|V|_{h_Y}
 \leq C|V|_{h_Y}.
 \end{aligned}
\end{equation}
Thus \eqref{eq:boundedcompositionchain} applies to the unit-length
functions $y\mapsto d(\iota\circ\pi)_yP_j(y)$ and gives
\begin{equation}\label{eq:ambientframechainrule}
 \begin{aligned}
 E_i\{d(\iota\circ\pi)_\Phi(P_j\circ\Phi)\}
 &=d(\iota\circ\pi)_\Phi
             \bigl(\nabla^\Phi_{E_i}(P_j\circ\Phi)\bigr)\\
 &\quad+\mathrm{II}_U\bigl(dU(E_i),
                  d(\iota\circ\pi)_\Phi(P_j\circ\Phi)\bigr),\\
 |E_i\{d(\iota\circ\pi)_\Phi(P_j\circ\Phi)\}|_{\R^m}
 &\leq C|dU(E_i)|_{\R^m}\quad\text{a.e.}
 \end{aligned}
\end{equation}
The weak product rule and \eqref{eq:ambientframechainrule} imply
\begin{equation}\label{eq:ambientcoefficientbounds}
 \begin{gathered}
 |f_1|+|f_2|\leq2|\mathbf b|_{\R^m},\qquad
 f_1,f_2\in H^1(M)\cap L^\infty(M),\\
 |df_1|_{g_0}+|df_2|_{g_0}
 \leq C\bigl(|d\mathbf b|_{g_0,\R^m}
       +|\mathbf b|_{\R^m}|dU|_{g_0,\R^m}\bigr).
 \end{gathered}
\end{equation}
By \eqref{eq:periodicmollification} and \eqref{eq:aeextraction},
choose smooth $f_{j,n}\to f_j$ in $H^1$ and a.e., with
$\|f_{j,n}\|_\infty\leq\|f_j\|_\infty$.
For $W_n:=\sum_j f_{j,n}(P_j\circ\Phi)$,
\eqref{eq:strongproductestimate} and \eqref{eq:framebounds} give
\begin{equation}\label{eq:ambienttestconvergence}
 \begin{gathered}
 \|W_n-W\|_{L^2(h_Y)}\longrightarrow0,\\
 \|\nabla^\Phi(W_n-W)\|_{L^2(g_0,h_Y)}
 \leq\sum_{j=1}^2\|d(f_{j,n}-f_j)\|_2
       +C\sum_{j=1}^2\|(f_{j,n}-f_j)|d\Phi|_{g_0,h_Y}\|_2
 \longrightarrow0.
 \end{gathered}
\end{equation}
Equation~\eqref{eq:ambientframechainrule} and the product rule give
\[
 E_i\bigl(\Pi(U)\mathbf b\bigr)
 =d(\iota\circ\pi)_\Phi(\nabla^\Phi_{E_i}W)
       +\mathrm{II}_U\bigl(dU(E_i),\Pi(U)\mathbf b\bigr).
\]
Pairing with the tangent vector $dU(E_i)$ removes the normal term
and gives \eqref{eq:ambienttestpairing}.
\end{proof}

\subsection*{Statement of usage of artificial intelligence}
The central mathematical ideas and overall approach are due to the authors.
OpenAI ChatGPT was used to assist with proofreading and improving the
clarity of the exposition. The authors take full responsibility for the
paper's content and the accuracy of its references.


\begin{thebibliography}{99}

\bibitem{AGindex}
B.~Ammann and N.~Ginoux,
\emph{Dirac-harmonic maps from index theory},
Calc. Var. Partial Differential Equations \textbf{47} (2013),
739--762.
{arXiv:1110.1309}.

\bibitem{AGexamples}
B.~Ammann and N.~Ginoux,
\emph{Some examples of Dirac-harmonic maps},
Lett. Math. Phys. \textbf{109} (2019), 1205--1218.



\bibitem{Bekesi}
B.~M.~B\'ek\'esi,
\emph{A note on coupled Dirac-harmonic maps},
{arXiv:2610.10330}, 2026.

\bibitem{CJLW}
Q.~Chen, J.~Jost, J.~Li, and G.~Wang,
\emph{Dirac-harmonic maps},
Math. Z. \textbf{254} (2006), 409--432.
{arXiv:math/0411402}.

\bibitem{CJSZsurfaces}
Q.~Chen, J.~Jost, L.~Sun, and M.~Zhu,
\emph{Dirac-harmonic maps between Riemann surfaces},
Asian J. Math. \textbf{23} (2019), 107--126.

\bibitem{ColdingMinicozzi}
T.~H.~Colding and W.~P.~Minicozzi II,
\emph{A Course in Minimal Surfaces},
Graduate Studies in Mathematics, vol.~121,
American Mathematical Society, Providence, RI, 2011.


\bibitem{FarbMargalit}
B.~Farb and D.~Margalit,
\emph{A Primer on Mapping Class Groups},
Princeton Mathematical Series, vol.~49,
Princeton University Press, Princeton, NJ, 2012.

\bibitem{GT}
D.~Gilbarg and N.~S.~Trudinger,
\emph{Elliptic Partial Differential Equations of Second Order},
Classics in Mathematics, Springer, Berlin, 2001.

\bibitem{Hatcher}
A.~Hatcher,
\emph{Algebraic Topology},
Cambridge University Press, Cambridge, 2002.

\bibitem{JostRS}
J.~Jost,
\emph{Compact Riemann Surfaces: An Introduction to Contemporary Mathematics},
third edition, Universitext, Springer, Berlin, 2006.

\bibitem{JostRGA}
J.~Jost,
\emph{Riemannian Geometry and Geometric Analysis},
eigth edition, Universitext, Springer, Cham, 2026.


\bibitem{JMZ}
J.~Jost, X.~Mo, and M.~Zhu,
\emph{Some explicit constructions of Dirac-harmonic maps},
J. Geom. Phys. \textbf{59} (2009), 1512--1527.


\bibitem{JSZ}
J.~Jost, L.~Sun, and J.~Zhu,
\emph{Dirac-harmonic maps with the trivial index},
Sci. China Math. \textbf{68} (2025), 917--938.

\bibitem{JZalpha}
J.~Jost and J.~Zhu,
\emph{$\alpha$-Dirac-harmonic maps from closed surfaces},
Calc. Var. Partial Differential Equations \textbf{60} (2021),
Paper No.~111.


\bibitem{JZflow}
J.~Jost and J.~Zhu,
\emph{Short-time existence of the $\alpha$-Dirac-harmonic map flow
and applications},
Comm. Partial Differential Equations \textbf{46} (2021), 442--469.


\bibitem{JZcurves}
J.~Jost and J.~Zhu,
\emph{Coupled Dirac-harmonic maps from closed spin tori},
{arXiv:2610.09754}, 2026.

\bibitem{LawsonMichelsohn}
H.~B.~Lawson, Jr., and M.-L.~Michelsohn,
\emph{Spin Geometry},
Princeton Mathematical Series, vol.~38,
Princeton University Press, Princeton, NJ, 1989.

\bibitem{LeeSM}
J.~M.~Lee,
\emph{Introduction to Smooth Manifolds},
second edition, Graduate Texts in Mathematics, vol.~218,
Springer, New York, 2013.

\bibitem{Sampson}
J.~H.~Sampson,
\emph{Some properties and applications of harmonic mappings},
Ann. Sci. \'Ecole Norm. Sup. (4) \textbf{11} (1978), no.~2, 211--228.


\bibitem{SchoenYau}
R.~Schoen and S.-T.~Yau,
\emph{On univalent harmonic maps between surfaces},
Invent. Math. \textbf{44} (1978), 265--278.

\bibitem{SharpTopping}
B.~Sharp and P.~Topping,
\emph{Decay estimates for Rivi\`ere's equation, with applications
to regularity and compactness},
Trans. Amer. Math. Soc. \textbf{365} (2013), 2317--2339.



\bibitem{WangXu}
C.~Wang and D.~Xu,
\emph{Regularity of Dirac-harmonic maps},
Int. Math. Res. Not. IMRN (2009), no.~20, 3759--3792.


\bibitem{Yang}
L.~Yang,
\emph{A structure theorem of Dirac-harmonic maps between spheres},
Calc. Var. Partial Differential Equations \textbf{35} (2009),
409--420.

\end{thebibliography}
\end{document}